%% file: neurips_2024.tex
\pdfoutput=1
\documentclass{article}

\PassOptionsToPackage{numbers}{natbib}
\usepackage[final]{neurips_2024}

\usepackage[utf8]{inputenc} % allow utf-8 input
\usepackage[T1]{fontenc}    % use 8-bit T1 fonts
\usepackage{hyperref}       % hyperlinks
\usepackage{url}            % simple URL typesetting
\usepackage{booktabs}       % professional-quality tables
\usepackage{amsfonts}       % blackboard math symbols
\usepackage{amsthm}
\usepackage{nicefrac}       % compact symbols for 1/2, etc.
\usepackage{microtype}      % microtypography
\usepackage{xcolor}         % colors

\title{Faster Accelerated First-order Methods for Convex Optimization with Strongly Convex Function Constraints}

\author{%
  Zhenwei~Lin\\
  Shanghai University of Finance and Economics\\
  \texttt{zhenweilin@163.sufe.edu.cn}
  \And
  Qi~Deng\\
  Shanghai Jiao Tong University\\
  \texttt{qdeng24@sjtu.edu.cn}
}

\newtheorem{assumption}{\protect\assumptionname}
\newtheorem{theorem}{\protect\theoremname}
\newtheorem{definition}{\protect\definitionname}
\newtheorem{remark}{\protect\remarkname}
\newtheorem{corollary}{\protect\corollaryname}
\newtheorem{proposition}{\protect\propositionname}
\newtheorem{lemma}{\protect\lemmaname}

\providecommand{\assumptionname}{Assumption}
\providecommand{\corollaryname}{Corollary}
\providecommand{\definitionname}{Definition}
\providecommand{\lemmaname}{Lemma}
\providecommand{\propositionname}{Proposition}
\providecommand{\remarkname}{Remark}
\providecommand{\theoremname}{Theorem}

\usepackage{amsmath}
\usepackage{array}
\usepackage{verbatim}
\usepackage{float}
\usepackage{calc}
\usepackage{mathtools}
\usepackage{bm}
\usepackage{multicol}
\usepackage{multirow}
\usepackage{amssymb}
\usepackage{graphicx}
\usepackage{wasysym}
\usepackage{url}
\PassOptionsToPackage{normalem}{ulem}
\usepackage{ulem}
\pdfpageheight\paperheight
\pdfpagewidth\paperwidth

\usepackage{titletoc}
\newcommand{\insertseparator}{%
  \par
  \noindent\rule{\linewidth}{1.5pt} % 粗的直线
  \par
}

\renewcommand\frac[2]{\tfrac{#1}{#2}}
\usepackage{algorithm}
\usepackage{algpseudocode}
\usepackage{graphicx}
\usepackage{array}
\usepackage{eqparbox}

\floatstyle{ruled}
\newfloat{algorithm}{tbp}{loa}
\providecommand{\algorithmname}{Algorithm}
\floatname{algorithm}{\protect\algorithmname}

\usepackage{etoolbox}  % patch def of algorithmic environment
\makeatletter
\patchcmd{\algorithmic}{\addtolength{\ALC@tlm}{\leftmargin} }{\addtolength{\ALC@tlm}{\leftmargin}}{}{}
\makeatother
\let\oldsum\sum
\renewcommand{\sum}{\textstyle\oldsum}

\begin{document}
\include{macros}

\maketitle

\begin{abstract}
\input{abstract_neurips}
\end{abstract}

\input{main_text_neurips}

\bibliographystyle{plainnat}
\bibliography{neurips_ref}

%%%%%%%%%%%%%%%%%%%%%%%%%%%%%%%%%%%%%%%%%%%%%%%%%%%%%%%%%%%%

\appendix
\newpage
\part{Appendix} % Start the appendix part
\insertseparator
\startcontents[appendix] 
\printcontents[appendix]{}{1}{}
\insertseparator
% \section{Appendix / supplemental material}

% Optionally include supplemental material (complete proofs, additional experiments and plots) in appendix.
% All such materials \textbf{SHOULD be included in the main submission.}

\input{appendix_neurips}

\end{document}

%% file: macros.tex
\global\long\def\inprod#1#2{\left\langle #1,#2\right\rangle }%

\global\long\def\inner#1#2{\langle#1,#2\rangle}%

\global\long\def\binner#1#2{\big\langle#1,#2\big\rangle}%

\global\long\def\norm#1{\Vert#1\Vert}%
\global\long\def\bnorm#1{\big\Vert#1\big\Vert}%
\global\long\def\Bnorm#1{\Big\Vert#1\Big\Vert}%

\newcommand{\intr}{\mathop{\bf int}}
\newcommand{\ri}{\mathop{\bf ri}}

\global\long\def\brbra#1{\big(#1\big)}%
\global\long\def\Brbra#1{\Big(#1\Big)}%
\global\long\def\rbra#1{(#1)}%
\global\long\def\sbra#1{[#1]}%
\global\long\def\bsbra#1{\big[#1\big]}%
\global\long\def\Bsbra#1{\Big[#1\Big]}%
\global\long\def\lrsbra#1{\left[#1\right]}%
\global\long\def\abs#1{\vert#1\vert}%
\global\long\def\babs#1{\big\vert#1\big\vert}%
\global\long\def\cbra#1{\{#1\}}%
\global\long\def\bcbra#1{\big\{#1\big\}}%
\global\long\def\Bcbra#1{\Big\{#1\Big\}}%

\global\long\def\vertiii#1{\left\vert \kern-0.25ex  \left\vert \kern-0.25ex  \left\vert #1\right\vert \kern-0.25ex  \right\vert \kern-0.25ex  \right\vert }%

\global\long\def\matr#1{\bm{#1}}%

\global\long\def\til#1{\tilde{#1}}%

\global\long\def\wtil#1{\widetilde{#1}}%

\global\long\def\wh#1{\widehat{#1}}%

\global\long\def\mcal#1{\mathcal{#1}}%

\global\long\def\mbb#1{\mathbb{#1}}%

\global\long\def\mtt#1{\mathtt{#1}}%

\global\long\def\ttt#1{\texttt{#1}}%

\global\long\def\dtxt{\textrm{d}}%

\global\long\def\bignorm#1{\bigl\Vert#1\bigr\Vert}%

\global\long\def\Bignorm#1{\Bigl\Vert#1\Bigr\Vert}%

\global\long\def\rmn#1#2{\mathbb{R}^{#1\times#2}}%

\global\long\def\deri#1#2{\frac{d#1}{d#2}}%

\global\long\def\pderi#1#2{\frac{\partial#1}{\partial#2}}%

\global\long\def\limk{\lim_{k\rightarrow\infty}}%

\global\long\def\trans{\textrm{T}}%

\global\long\def\onebf{\mathbf{1}}%

\global\long\def\zerobf{\mathbf{0}}%

\global\long\def\zero{\bm{0}}%

% expectation

\global\long\def\Euc{\mathrm{E}}%

\global\long\def\Expe{\mathbb{E}}%

\global\long\def\rank{\mathrm{rank}}%

\global\long\def\range{\mathrm{range}}%

\global\long\def\diam{\mathrm{diam}}%

\global\long\def\epi{\mathrm{epi} }%

\global\long\def\inte{\operatornamewithlimits{int}}%

\global\long\def\cov{\mathrm{Cov}}%

\global\long\def\argmin{\operatornamewithlimits{argmin}}%

\global\long\def\argmax{\operatornamewithlimits{argmax}}%

\global\long\def\proj{\operatornamewithlimits{proj}}%

\global\long\def\prox{\operatornamewithlimits{prox}}%

\global\long\def\tr{\operatornamewithlimits{tr}}%

\global\long\def\dis{\operatornamewithlimits{dist}}%

\global\long\def\sign{\operatornamewithlimits{sign}}%

\global\long\def\prob{\mathrm{Prob}}%

\global\long\def\st{\operatornamewithlimits{s.t.}}%

\global\long\def\dom{\mathrm{dom}}%

\global\long\def\prox{\mathrm{prox}}%

\global\long\def\diag{\mathrm{diag}}%

\global\long\def\and{\mathrm{and}}%

\global\long\def\st{\mathrm{s.t.}}%

\global\long\def\Var{\operatornamewithlimits{Var}}%

\global\long\def\raw{\rightarrow}%

\global\long\def\law{\leftarrow}%

\global\long\def\Raw{\Rightarrow}%

\global\long\def\Law{\Leftarrow}%

\global\long\def\vep{\varepsilon}%

\global\long\def\dom{\operatornamewithlimits{dom}}%

\global\long\def\tsum{{\textstyle {\sum}}}%

\global\long\def\Cbb{\mathbb{C}}%

\global\long\def\Ebb{\mathbb{E}}%

\global\long\def\Fbb{\mathbb{F}}%

\global\long\def\Nbb{\mathbb{N}}%

\global\long\def\Rbb{\mathbb{R}}%

\global\long\def\extR{\widebar{\mathbb{R}}}%

\global\long\def\Pbb{\mathbb{P}}%

\global\long\def\Mrm{\mathrm{M}}%

\global\long\def\Acal{\mathcal{A}}%

\global\long\def\Bcal{\mathcal{B}}%

\global\long\def\Ccal{\mathcal{C}}%

\global\long\def\Dcal{\mathcal{D}}%

\global\long\def\Ecal{\mathcal{E}}%

\global\long\def\Fcal{\mathcal{F}}%

\global\long\def\Gcal{\mathcal{G}}%

\global\long\def\Hcal{\mathcal{H}}%

\global\long\def\Ical{\mathcal{I}}%

\global\long\def\Kcal{\mathcal{K}}%

\global\long\def\Lcal{\mathcal{L}}%

\global\long\def\Mcal{\mathcal{M}}%

\global\long\def\Ncal{\mathcal{N}}%

\global\long\def\Ocal{\mathcal{O}}%

\global\long\def\Pcal{\mathcal{P}}%

\global\long\def\Scal{\mathcal{S}}%

\global\long\def\Tcal{\mathcal{T}}%

\global\long\def\Xcal{\mathcal{X}}%

\global\long\def\Ycal{\mathcal{Y}}%

\global\long\def\Zcal{\mathcal{Z}}%

\global\long\def\i{i}%
% Use bold text symbol for vector and matrix

\global\long\def\abf{\mathbf{a}}%

\global\long\def\bbf{\mathbf{b}}%

\global\long\def\cbf{\mathbf{c}}%

\global\long\def\fbf{\mathbf{f}}%

\global\long\def\gbf{\mathbf{g}}%

\global\long\def\qbf{\mathbf{q}}%

\global\long\def\lambf{\bm{\lambda}}%

\global\long\def\alphabf{\bm{\alpha}}%

\global\long\def\sigmabf{\bm{\sigma}}%

\global\long\def\thetabf{\bm{\theta}}%

\global\long\def\deltabf{\bm{\delta}}%

\global\long\def\mubf{\bm{\mu}}%

\global\long\def\lbf{\mathbf{l}}%

\global\long\def\ubf{\mathbf{u}}%

\global\long\def\vbf{\mathbf{v}}%

\global\long\def\wbf{\mathbf{w}}%

\global\long\def\xbf{\mathbf{x}}%

\global\long\def\ybf{\mathbf{y}}%

\global\long\def\zbf{\mathbf{z}}%

\global\long\def\Abf{\mathbf{A}}%

\global\long\def\Ubf{\mathbf{U}}%

\global\long\def\Pbf{\mathbf{P}}%

\global\long\def\Ibf{\mathbf{I}}%

\global\long\def\Ebf{\mathbf{E}}%

\global\long\def\Mbf{\mathbf{M}}%

\global\long\def\Qbf{\mathbf{Q}}%

\global\long\def\Lbf{\mathbf{L}}%

\global\long\def\Pbf{\mathbf{P}}%

\global\long\def\Xbf{\mathbf{X}}%
% Use bold symbol for vector and matrix

\global\long\def\abm{\bm{a}}%

\global\long\def\bbm{\bm{b}}%

\global\long\def\cbm{\bm{c}}%

\global\long\def\dbm{\bm{d}}%

\global\long\def\ebm{\bm{e}}%

\global\long\def\fbm{\bm{f}}%

\global\long\def\gbm{\bm{g}}%

\global\long\def\hbm{\bm{h}}%

\global\long\def\pbm{\bm{p}}%

\global\long\def\qbm{\bm{q}}%

\global\long\def\rbm{\bm{r}}%

\global\long\def\sbm{\bm{s}}%

\global\long\def\tbm{\bm{t}}%

\global\long\def\ubm{\bm{u}}%

\global\long\def\vbm{\bm{v}}%

\global\long\def\wbm{\bm{w}}%

\global\long\def\xbm{\bm{x}}%

\global\long\def\ybm{\bm{y}}%

\global\long\def\zbm{\bm{z}}%

\global\long\def\Abm{\bm{A}}%

\global\long\def\Bbm{\bm{B}}%

\global\long\def\Cbm{\bm{C}}%

\global\long\def\Dbm{\bm{D}}%

\global\long\def\Ebm{\bm{E}}%

\global\long\def\Fbm{\bm{F}}%

\global\long\def\Gbm{\bm{G}}%

\global\long\def\Hbm{\bm{H}}%

\global\long\def\Ibm{\bm{I}}%

\global\long\def\Jbm{\bm{J}}%

\global\long\def\Lbm{\bm{L}}%

\global\long\def\Obm{\bm{O}}%

\global\long\def\Pbm{\bm{P}}%

\global\long\def\Qbm{\bm{Q}}%

\global\long\def\Rbm{\bm{R}}%

\global\long\def\Ubm{\bm{U}}%

\global\long\def\Vbm{\bm{V}}%

\global\long\def\Wbm{\bm{W}}%

\global\long\def\Xbm{\bm{X}}%

\global\long\def\Ybm{\bm{Y}}%

\global\long\def\Zbm{\bm{Z}}%

\global\long\def\lambm{\bm{\lambda}}%

\global\long\def\alphabm{\bm{\alpha}}%

\global\long\def\albm{\bm{\alpha}}%

\global\long\def\taubm{\bm{\tau}}%

\global\long\def\mubm{\bm{\mu}}%

\global\long\def\yrm{\mathrm{y}}%

\global\long\def\dist{\mathop{\bf dist}}%

\global\long\def\brbra#1{\big(#1\big)}%
\global\long\def\Brbra#1{\Big(#1\Big)}%
\global\long\def\rbra#1{(#1)}%
\global\long\def\sbra#1{[#1]}%
\global\long\def\bsbra#1{\big[#1\big]}%
\global\long\def\Bsbra#1{\Big[#1\Big]}%
\global\long\def\abs#1{\vert#1\vert}%
\global\long\def\babs#1{\big\vert#1\big\vert}%
\global\long\def\cbra#1{\{#1\}}%
\global\long\def\bcbra#1{\big\{#1\big\}}%
\global\long\def\Bcbra#1{\Big\{#1\Big\}}%
\global\long\def\vrho{\hat{\rho}}
\global\long\def\trho{\tilde{\rho}}
\global\long\def\coloneqq{:=}

\global\long\def\aleq{\overset{(a)}{\leq}}%
\global\long\def\aeq{\overset{(a)}{=}}%
\global\long\def\ageq{\overset{(a)}{\geq}}%
\global\long\def\bleq{\overset{(b)}{\leq}}%
\global\long\def\beq{\overset{(b)}{=}}%
\global\long\def\bgeq{\overset{(b)}{\geq}}%
\global\long\def\cleq{\overset{(c)}{\leq}}%
\global\long\def\ceq{\overset{(c)}{=}}%
\global\long\def\cgeq{\overset{(c)}{\geq}}%
\global\long\def\dleq{\overset{(d)}{\leq}}%
\global\long\def\deq{\overset{(d)}{=}}%
\global\long\def\dgeq{\overset{(d)}{\geq}}%
\global\long\def\eleq{\overset{(e)}{\leq}}%
\global\long\def\eeq{\overset{(e)}{=}}%
\global\long\def\egeq{\overset{(e)}{\geq}}%
% \newcommand{\Halmos}{\ensuremath{\blacksquare}}
% \global\long\def\endproof{\hfill \Halmos}%

\global\long\def\BPu{\textsc{BPu}}%
\global\long\def\APDSCE{\textup{APDPro}}%
\global\long\def\IMPROVE{\textsc{Improve}}%
\global\long\def\UpdateN{\textsc{TerminateIter}}%

\global\long\def\APDPI{\textup{APDPi}}%
\global\long\def\APD{\textup{APD}}%
\global\long\def\resapdsce{%
\textup{rAPDPro}}%
\global\long\def\restpapd{\text{Restart Version of TPAPD}}%
\global\long\def\TPAPD{\textup{rAPDPro}}%
\global\long\def\MSAPD{\textup{msAPD}}%
\global\long\def\mirror{\textup{Mirror-Prox}}%
\global\long\def\nonde{\text{non-degenerated}}%

%% file: abstract_neurips.tex
In this paper, we introduce faster accelerated primal-dual algorithms for minimizing a convex function subject to strongly convex function constraints. 
Prior to our work, the best complexity bound was $\mathcal{O}(1/{\varepsilon})$, regardless of the strong convexity of the constraint function.
It is unclear whether the strong convexity assumption can enable even better convergence results. 
To address this issue, we have developed novel techniques to progressively estimate the strong convexity of the Lagrangian function.
Our approach, for the first time, effectively leverages the constraint strong convexity, obtaining an improved complexity of $\mathcal{O}(1/\sqrt{\varepsilon})$. This rate matches the complexity lower bound for strongly-convex-concave saddle point optimization and is therefore order-optimal.
We show the superior performance of our methods in sparsity-inducing constrained optimization, notably Google's personalized PageRank problem. Furthermore, we show that a restarted version of the proposed methods can effectively identify the optimal solution's sparsity pattern within a finite number of steps, a result that appears to have independent significance.

% Prior to our work, the best complexity bound was obtained as $\mathcal{O}(1/{\sqrt{\varepsilon}})$ if the strong convexity of the objective function was used and $\mathcal{O}(1/{\varepsilon})$ for merely convex objective. However, it remains unclear how to improve $\mathcal{O}(1/{\varepsilon})$ by leveraging the strong convexity assumption of constrained function.  We address this issue by developing novel techniques to estimate the strong convexity of the Lagrangian function progressively.
% Our approach yields an improved complexity of $\mathcal{O}(1/\sqrt{\varepsilon})$,   matching the complexity lower bound for strongly-convex-concave saddle point optimization.
% We show the superior performance of our methods in sparsity-inducing constrained optimization, notably Google's personalized PageRank problem. Furthermore, we show that a restarted version of the proposed methods can effectively identify the optimal solution's sparsity pattern within a finite number of steps, a result that appears to have independent significance.

%% file: main_text_neurips.tex
\section{Introduction}
In this paper, we are interested in the following convex function-constrained
problem:
\begin{equation}
\begin{split}\textstyle\min_{\text{\ensuremath{\xbf}}\in\mathbb{R}^{n}}  \quad f(\text{\ensuremath{\xbf}}) \quad
\st  \quad g_i(\text{\ensuremath{\xbf}})\leq 0,\ 1\le i \le m,
\end{split}
\label{eq:constrainedproblem}
\end{equation}
where $f:\mathbb{R}^{n}\to\mathbb{R}$  is a convex continuous function and bounded from below and  $g_{i}:\mathbb{R}^{n}\to\mathbb{R}$, $i=1,2,\ldots, m$, are {strongly} convex continuous functions. 
% Constrained problems such as \eqref{eq:constrainedproblem} frequently appear in machine learning and operations research. 
An important application of this problem, commonly encountered in statistics and engineering, involves the objective $f(\xbf)$ as a proximal-friendly regularizer and $g_i(\xbf)$ as a data-driven loss function used to gauge model fidelity.

% The field of first-order methods for optimization problems is comprehensive, with significant attention dedicated to unconstrained problems or those within projection-efficient or linearly constrained domains. Recent advancements in this area can be found in works like \cite{nesterov2018lectures, lan2020firstorder}.
To apply first-order methods for the above function-constrained problems, 
a common strategy involves a double-loop procedure that repeatedly employs fast first-order methods, such as Nesterov's accelerated method,  to solve specific strongly convex proximal subproblems.
Popular methods among this category include Augmented Lagrangian methods~\citep{lan2016iterationcomplexity, xu2021iterationa}, level-set methods~\cite{lin2018levelset}, penalty  methods~\cite{lan2013iterationcomplexitya}. 
When both $f(\xbf)$ and $g_i(\xbf)$ are convex and smooth (or composite), it has been found that these double-loop algorithms can attain an iteration complexity of $\Ocal(1/\vep)$ to achieve an $\vep$-error in both the optimality gap and constraint violation.  When the objective is strongly convex, the complexity can be further improved to $\Ocal(1/\sqrt{\vep})$ (\cite{xu2021iterationa, lin2018levelset}).
 % For example, \cite{boob2022level} proposed a level constrained proximal gradient (LCPG) was proposed, which can further reduce the complexity associated with evaluating function and gradient values. Notably, when $f(\xbf)$ is strongly convex, LCPG  only requires an $\Ocal(\log(1/\vep))$ number of function and gradient evaluations. However, it is worth noting that LCPG still necessitates solving a structured strongly convex constrained subproblem to a predetermined accuracy.

In contrast to these double-loop algorithms, single-loop algorithms remain popular due to their simplicity in implementation. Along this research line, \cite{xu2021firstordera} developed a first-order algorithm based on linearizing the augmented Lagrangian function, which obtains an iteration complexity of $\Ocal(1/\vep)$.
\cite{zhang2022solving} extended the augmented Lagrangian method to stochastic function-constrained problems where both objective and constraint exhibit an expectation form.  
Viewing \eqref{eq:constrainedproblem} as a special case of the min-max problem:
\begin{equation}\label{eq:min-max}
    \textstyle\min_{\xbf\in\Rbb^n} \max_{\ybf\in\Rbb^m}\ \Lcal(\xbf,\ybf)\coloneqq f(\xbf) + \sum_{i=1}^m y_i {g_i(\xbf)}, \quad\text{s.t.}\ y_i\ge0,\ i=1,2,\ldots, m,
\end{equation}
\cite{hamedani2021primal} proposed to solve~\eqref{eq:constrainedproblem} and~\eqref{eq:min-max} by an accelerated primal-dual method (\APD{}), which generalizes the primal-dual hybrid gradient method~\citep{chambolle2016ergodic} initially developed for saddle point optimization with bilinear coupling term. 
Under mild conditions, \APD{}  achieves the best iteration complexity of $\Ocal(1/\vep)$   for general convex constrained problem and a further improved complexity of $\Ocal(1/\sqrt{\vep})$ when $f(\xbf)$ is strongly convex.
\cite{boob2019stochastic} proposed a unified constrained extrapolation method that can be applied to both deterministic and stochastic constrained optimization problems.

Despite these recent progresses, to the best of our knowledge, all available algorithms are suboptimal in the presence of strongly convex function constraints~\eqref{eq:constrainedproblem}. Specifically,  direct applications of previously discussed algorithms yield an $\Ocal(1/\vep)$ complexity, which is inferior to the   $\Ocal(1/\sqrt{\vep})$ optimal bound for the strongly-convex-concave saddle point problem~\citep{lin2020near}. 
It is somewhat unsatisfactory that the strong convexity of $g(\xbf)$ has not been found helpful in further algorithmic acceleration. 
The core underlying issue arises from the dynamics of saddle point optimization: it is the strong convexity of $\Lcal(\cdot,\ybf)$ that offers more potential acceleration advantages, yet the strong convexity of $\Lcal(\cdot,\ybf)$  is substantially harder to estimate than that of $g(\xbf)$.  This difficulty is compounded by the interplay between  $g(\xbf)$ and the varying dual sequence $\{\ybf_k\}$.
The challenge naturally leads us to question: \emph{Is it possible to further improve the convergence rate of first-order methods for solving the strongly convex constrained problem~\eqref{eq:constrainedproblem}?}

\paragraph{Key intuitions}
We make an assumption that the minimizer of $f(\xbf)$ is infeasible for the function constraint $g_i(\xbf)\le 0$, $1\le i \le m$. If this assumption were not made, we would be dealing with an unconstrained optimization problem that would not depend on $g(\xbf)$. This assumption also implies that the optimal dual variables are non-zero, and as a result, the Lagrangian function is strongly convex with respect to $\xbf$. By leveraging the strong convexity, we can use more aggressive step sizes and achieve faster convergence rates compared to other state-of-the-art algorithms.

% We assume that at least one of the constraints $\{g_i(\xbf)\}$ is active at the optimal solution. Without it, we would essentially be dealing with an unconstrained optimization that does not rely on $g(\xbf)$.  Additionally, this assumption ensures strict separation between the origin and the subdifferential at optimal solution $\partial f(\xbf^*)$.  By leveraging the structures of such an objective function and the strongly convex constraints,  we develop new lower bounds for the norm of optimal dual variables $\norm{\ybf^*}$. This derived bound,  valid under mild assumptions and easily computable using the generated iterates, allows us to harness the strong convexity of the Lagrangian function, apply more aggressive stepsizes, and hence achieve faster convergence rates compared to state-of-the-art algorithms.

\paragraph{Applications in sparsity-constrained optimization}
We consider the constrained Lasso-type problem, which minimizes a sparsity-inducing regularizer while explicitly ensuring data-driven error remains controlled:
\begin{equation}
\textstyle\min_{\xbf\in\mathbb{R}^{n}}  \ \|\xbf\|_{1} \quad  \text{s.t.} \ g(\xbf)\le 0,
\label{eq:constrained-lasso}
\end{equation}
where $g(\cdot)$ is a convex smooth loss term. A motivating application is the approximate personalized PageRank problem~\citep{fountoulakis2019variational}, where $g(\xbf)=\frac{1}{2}\inner{\xbf}{Q\xbf}-\inner{\bbf}{\xbf}$ is strongly convex quadratic and $Q$ integrates the graph Laplacian with an identity matrix. 
Compared to the standard Lasso problem~\citep{tibshirani1996regression}, 
$\textstyle\min_{\xbf\in\Rbb^n} g(\xbf) + \lambda\norm{\xbf}_1,$
the constrained problem~\eqref{eq:constrained-lasso} offers enhanced control over the data fitting error. This advantage, however, is counterbalanced by the  challenge of dealing with a nonlinear constraint.
Besides concerns about the efficiency in solving \eqref{eq:constrained-lasso},  it is often desired to show the active set (or sparsity) identification, namely, the nonzero patterns of the optimal solution $\xbf^*$ can be identified by the solution sequence $\{\xbf_k\}$ in a finite number of iterations. 
% Leverage the low-dimensional solution structure, one for further acceleration~\cite{iutzeler2020nonsmoothness}.
Identifying the embedded solution structure within a broader context is referred to as the manifold identification problem~\citep{wright1993identifiablea, hare2004identifying}. Exploiting the sparsity pattern is particularly desirable in large-scale PageRank problems, as it could result in significant runtime savings.
For the regularized Lasso-type problem, it has been known that proximal gradient methods~(e.g. \cite{iutzeler2020nonsmoothness, lee2012manifold, nutini2019activeset}) possess the finite active-set identification property. 
Specifically, \cite{nutini2019activeset} introduced ``active set complexity'', which is defined as the number of iterations required before an algorithm is guaranteed to have reached the optimal manifold,
and they proved the proximal gradient method with constant stepsize
can identify the optimal manifold in a finite number of iterations.
However,  for the problem~\eqref{eq:constrained-lasso}, it remains unclear whether first-order methods can identify the sparsity pattern in finite time. 

% and a regularized dual-averaging method can identify the sparsity pattern with high probability in stochastic setting \cite{lee2012manifold}.
% \cite{nutini2019activeset} introduced 
% ``active set complexity'', which is defined as the number of iterations
% required before an algorithm is guaranteed to have reached the optimal manifold,
% and they proved the proximal gradient method with constant stepsize
% can identify the optimal manifold in a finite number of iterations.

% \subsection{Contributions}
\paragraph{Contributions}
We address the theoretical questions about strongly convex constrained optimization and the application of sparse optimization. Our contributions are summarized as follows.

First, we present a new accelerated primal-dual algorithm with progressive strong convexity estimation ({\APDSCE}) for solving problem~\eqref{eq:constrainedproblem}.
\APDSCE{} employs a novel strategy to estimate the lower bound of the dual variables, which leads to a gradually refined estimated strong convexity modulus of $\Lcal(\cdot, \ybf)$. With additional cut constraints on the dual update, \APDSCE{} is able to separate the dual search space from the origin point, which is critical for maintaining the desired strong convexity over the entire solution path.  With these two important ingredients,  \APDSCE{} exhibits an $\mcal O\brbra{{(\|\xbf_0-\xbf^*\|+D_Y)}/{\sqrt{\vep}}}$ complexity bound to obtain an $\vep$-error on the  function value gap and constraint violation, where $D_Y$ is a known upper-bound of $\|\ybf_0-\ybf^*\|$.
Moreover, we show that for the last iterate to have an $\vep$ error (i.e.,  $\|\xbf_K-\xbf^*\|^2\le \vep$), \APDSCE{} requires a total iteration of $\Ocal\brbra{(\|\xbf_0-\xbf^*\|+\|\ybf_0-\ybf^*\|)/\sqrt{\vep}}$. Both complexity results appear new in the literature for strongly convex-constrained optimization.

Second, we present a new restart algorithm (\resapdsce{})  which calls {\APDSCE} repeatedly with the input parameters properly changing over time. 
 Different from \APDSCE{}, \resapdsce{} dynamically adjusts the iteration number of \APDSCE{} in each epoch based on the progressive strong convexity estimation. 
 We show that \resapdsce{} exhibits a complexity of $\Ocal\brbra{\log(D_X/\sqrt{\vep})+D_Y/\sqrt{\vep}}$ to ensure $\vep$-error in the last iterate convergence where $D_X$ is the estimated diameter of the primal feasible domain.  
 While it is difficult to improve  the overall $\Ocal(1/\sqrt{\vep})$ bound, \resapdsce{} appears to be more advantageous when $D_X$ and $D_Y$ are the same order of $\|\xbf_0-\xbf^*\|$ and $\|\ybf_0-\ybf^*\|$, respectively, and  $D_X\gg D_Y$. 
 In addition, we show that a similar restart strategy can further accelerate the standard \APD{}.
 % that does not utilize the strong convexity assumption. 
 The 
 % resulting 
 multistage-accelerated primal dual method~(\MSAPD{}) obtains a comparable $\Ocal(1/\sqrt{\vep})$ complexity of \APDSCE{} without introducing additional cut constraint.

Third, we apply our proposed methods to the sparse learning problem~\eqref{eq:constrained-lasso}. In view of the theoretical analysis, all our methods converge at an $\Ocal(1/\sqrt{\vep})$ rate, which is substantially better than the rates of state-of-the-art first-order algorithms. Moreover, we conduct a new analysis to show that the restart algorithm~\resapdsce{} has the favorable feature of identifying the optimal sparsity pattern.  Note that such active-set/manifold identification is substantially more challenging to prove due to the coupling of dual variables and constraint functions. To establish the desired property, we develop asymptotic convergence of the dual sequence to the optimal solution, which can be of independent interest.

% \paragraph{Comparison with  Frank-Wolfe}
% We note that the strongly convex function constraint in \eqref{eq:constrainedproblem} is a special case of a strongly convex set constraint, as demonstrated in~\cite{journee2010generalized}.
% Over the strongly convex set, it has been shown that Frank-Wolfe Algorithm (FW) can obtain convergence rates substantially better than the worst-case $\Ocal(1/\vep)$ rate. Under the bounded gradient assumption, \cite{dunn1979rates, levitin1966constrained} show that FW obtains linear convergence over a strongly convex set. \cite{garber2015faster} shows that FW obtains an $\Ocal(1/\sqrt{\vep})$ rate when the gradient is the order of the square root of the function value gap.  For more recent progress, please refer to~\cite{braun2022conditional}. Despite the attractive convergence property,  FW exhibits certain limitations when applied to the general function constraints~\eqref{eq:constrainedproblem} addressed in this paper.  Specifically, FW involves a sequence of linear optimization problems throughout the iterations. While linear optimization over certain strongly convex sets, such as $\ell_p$-ball, admits a closed-form solution, there exists no efficient routine to handle general function constraints explored in this paper. 

% \subsection{Outline}
\paragraph{Outline}
Section~\ref{sec:Preliminaries} sets notations and assumptions for the later analysis. Section~\ref{sec:APDSCE} presents the {\APDSCE} algorithm and develops its stepsize rule and complexity rate. 
Section~\ref{sec:restarted} presents the restart \APDSCE{} (\resapdsce{}) algorithm. 
% Section~\ref{sec:MSAPD} presents the multistage algorithm \MSAPD{}, which does not rely on adding a dual cut constraint. 
Section~\ref{sec:Manifold} applies our proposed methods for sparsity-inducing optimization
and shows the sparsity identification result for {\resapdsce}. Section~\ref{sec:numerical} empirically examine the convergence performance and sparsity identification of our proposed algorithms. Finally, we draw the conclusion in Section~\ref{sec:Conclusion}. All the missing proofs are provided in the appendix sections.

\section{\label{sec:Preliminaries}Preliminaries}
We use bold letters like $\xbf$ to represent vectors.
Suppose $\xbf\in\mathbb{R}^{n}$, $q\ge 1$, we use $\|\xbf\|_{q}=(\sum_{i=1}^{n}|\xbf_{(i)}|^{q})^{1/q}$
to represent the $l_{q}$-norm, where $\xbf_{(i)}$ is the $i$-th
element of $\xbf$. For brevity, $\|\xbf\|$ stands for $l_{2}$-norm. 
For a matrix $A$, we denote the matrix norm induced by 2-norm as
$\|A\|=\sup_{\|\xbf\|\leq1}\|A\xbf\|$.
The normal cone of $\mathcal{U}$ at $\ubf$ is denoted as $\mathcal{N}_{\mathcal{U}}(\ubf):=\left\{ \vbf\mid\left\langle \vbf,\xbf-\ubf\right\rangle \leq 0,\forall\xbf\in\mathcal{U}\right\} $.  Let $\mathcal{B}(\xbf,r)$
be the closed ball centered at $\xbf$ with radius $r > 0$,
i.e., $\mathcal{B}(\xbf,r)=\left\{ \ybf\mid\|\ybf-\xbf\|\leq r\right\} $.
We denote the set of feasible solutions by  $\mcal{X}_{G}\coloneqq\left\{ \xbf \mid g_i(\xbf) \leq 0, \forall i\in[m] \right\}$ and write the constraint function  as $G(\xbf)\coloneqq[g_{1}(\xbf),\ldots,g_{m}(\xbf)]^{\top}$. We assume 
each $g_{i}(\xbf)$ is a $\mu_i$ strongly convex function, and denote $\mubm:=[\mu_{1},\ldots,\mu_{m}]^{\top}$.
Let $[m]:=\{1,\ldots,m\}$ for integer $m$. We denote minimum and maximum strongly convexity $\underline{\mu}:=\min_{j\in[m]}\{\mu_{j}\}$,
and $\bar{\mu}:=\max_{j\in[m]}\{\mu_{j}\}$ and the vector of elements 0 by $\mathbf{0}$. 
The Lagrangian function
of problem~\eqref{eq:constrainedproblem} is given by $\mcal L(\xbf,\ybf):=f(\xbf)+\left\langle \ybf,G(\xbf)\right\rangle$ where $\ybf \in \Rbb^m_+$.
\begin{definition}[KKT condition]
    We say that $\xbf^*$  satisfies the \emph{KKT condition} of~\eqref{eq:constrainedproblem} if there exists a Lagrangian multiplier vector $\ybf^*\in\mathbb{R}^m_+$ such that $\zerobf \in \partial_x \Lcal(\xbf^*, \ybf^*)$ and $\langle \ybf^*, G(\xbf^*)\rangle =0 $.
\end{definition}

% \begin{equation*}\label{eq:KKT-cond}
% \zerobf \in \partial \Lcal(\xbf^*, \ybf^*),\  \zerobf \le \ybf^* \perp -G(\xbf^*) \ge \zerobf.
% \end{equation*}
The KKT condition is necessary for optimality when a constraint qualification (CQ) holds at $\xbf^*$. We assume Slater's CQ (Assumption~\ref{assu:Slater's}) holds, which guarantees that an optimal solution is also a KKT point~\citep{bertsekasnonlinear}. 
\begin{assumption}
\label{assu:Slater's}There exists a strictly feasible point $\widetilde{\xbf}\in \mathbb{R}^n$ such that $G(\widetilde{\xbf})<\mathbf{0}$.
\end{assumption}
We use $\tilde{\xbf}$ to denote a strictly feasible point throughout the paper.
Moreover, we require Assumption~\ref{assu:dual-optim} to circumvent any trivial solution.
\begin{assumption}
\label{assu:dual-optim}For any $\xbf_{0}^{*}\in\argmin_{\xbf\in\mbb R^{n}}f(\xbf)$,
there exists an $i\in[m]$ such that $g_{i}(\xbf_{0}^{*})>0$.
\end{assumption}
\begin{remark}
Assumption~\ref{assu:dual-optim} is essential for our analysis. While verifying Assumption~\ref{assu:dual-optim} can be indeed challenging,  it is achievable for the sparsity-inducing problem considered in our paper. In this example, the solution $\xbf_0^*=\boldsymbol{0}$ is the single minimizer of the sparsity penalty.
\end{remark}
% \begin{remark}
%     It is important to further highlight that Assumption~\ref{assu:dual-optim} is generally difficult to verify in most cases. Unconstrained convex problems may have an infinite number of optimal solutions. Therefore, it is impossible to verify that all $\xbf^*_0$ satisfy $g_i(\xbf_0^*)>0$. However, we also need to point out that it is possible when projection onto the level set of $f(\xbf)$ is easy. First, we compute an $\vep$-solution $\hat{\xbf}_0^*$ of $\min f(\xbf)$ efficiently, e.g., $\mcal O(1/\sqrt{\vep})$ by using accelerated method~\cite{allen2016optimal}. Next, we can switch the objective and constraint and consider problem $\phi^*=\min_{\xbf}\max_{1\le i\le m}g_i(\xbf)\ \ \st f(\xbf)\leq f(\hat{x}_0^*)$. If the projection is easy to do, we can find an approximate solution $\bar{\xbf}$ and value $\bar{\phi}$ satisfying $\bar{\phi}\le \phi^*+\vep$ in $\mcal O(1/\sqrt{\vep})$ by using accelerated gradient method~\cite{lin2018levelset,allen2016optimal}. If $\bar{\phi}>\vep$, then we have $\phi^*>0$ and hence Assumption~\ref{assu:dual-optim} holds. Otherwise, we have $\phi^*\leq \bar{\phi}\leq\vep$, then the solution $\bar{\xbf}$ naturally becomes an $\vep$-solution of the original problem. (\textcolor{red}{check here})
% \end{remark}

Next, we give several useful properties about the optimal solutions of problem~\eqref{eq:constrainedproblem}. Please refer to Appendix~\ref{sec:proof_of_prop_bounded_set_y} for the proof of Proposition~\ref{prop:bounded_set_y} and Appendix~\ref{sec:proof_of_proposition_uniquex} for the proof of Proposition~\ref{prop:uniquex}.
\begin{proposition}
\label{prop:bounded_set_y}Suppose Assumption~\ref{assu:Slater's}
holds. Then, for any optimal solution $\xbf^*$ of problem~\eqref{eq:constrainedproblem}, there exists $\ybf^*\in\mbb R^m$ such that KKT condition holds. Moreover,  $\ybf^*$  falls into set $\mcal Y:=\bcbra{\ybf\mid \norm{\ybf}_1\leq \bar{c}}$,  where $\bar{c}\coloneqq \frac{f(\widetilde{\xbf})-\min_{\xbf\in \mbb R^n}f(\xbf)}{\min_{i\in[m]}\{-g_{i}(\widetilde{\xbf})\}}$.
\end{proposition}
% \begin{proof}
%     For complete proof details, please refer to Appendix~\ref{sec:proof_of_prop_bounded_set_y}.
% \end{proof}
% \begin{proof}
% \proof 
% Under Slater's CQ, the existence of the KKT condition is a standard result in nonlinear programming. For example, one can refer to \cite{bertsekasnonlinear}.
% For any  $\xbf\in\Xcal_G$,  we have 
% \[
% f(\xbf)+\left\langle \ybf^{*},G(\xbf)\right\rangle \geq f(\xbf^{*})+\left\langle \ybf^{*},G(\xbf^{*})\right\rangle =f({\xbf}^{*}),
% \]
% where the equality is from the complementary slackness. In view of the above result and the Slater's condition (i.e., $G(\widetilde{\xbf})<\boldsymbol{0}$),
% we have
% \begin{equation}
% \begin{split} & f(\widetilde{\xbf})\geq f(\widetilde{\xbf})+\left\langle \ybf^{*},G(\widetilde{\xbf})\right\rangle \geq f(\xbf^{*}),\\
% \text{\ensuremath{\Rightarrow}} & \|\ybf^{*}\|_{1}\min_{i\in[m]}\{-g_{i}(\widetilde{\xbf})\}\leq-\left\langle \ybf^{*},G(\widetilde{\xbf})\right\rangle \leq f(\widetilde{\xbf})-f(\xbf^{*}),\\
% \Rightarrow & \|\ybf^{*}\|\leq\|\ybf^{*}\|_{1}\leq\frac{f(\widetilde{\xbf})-f(\xbf^{*})}{\min_{i\in[m]}\{-g_{i}(\widetilde{\xbf})\}}\leq \bar{c},
% \end{split}
% \label{eq:tighter}
% \end{equation}
% where the last inequality is by $f(\xbf^{*})\ge\min_{\xbf\in\mathbb{R}^{n}}f(\xbf)$.
% % , which is derived by $\xbf_{0}^{*}\in\argmin_{\xbf\in\mathbb{R}^{n}}f(\xbf)$.
% \endproof

\begin{proposition}\label{prop:uniquex}
Under Assumption~\ref{assu:dual-optim}, $\xbf^{*}$ is the unique solution of~\eqref{eq:constrainedproblem}. Furthermore, set  $\mathcal{Y}^{*}=\argmax_{\ybf\in \mbb R_+^m}\mcal L(\xbf^{*},\ybf)$
% denote the set containing all the optimal dual solutions,
% then $\mathcal{Y}^{*}$
is convex and bounded. 
\end{proposition}

In view of Assumption~\ref{assu:dual-optim}, Proposition~\ref{prop:uniquex}, and closedness of the subdifferential set of proper convex functions~\citep[Theorem 3.9]{beck2017first},~\cite[Chapter 23]{rockafellar1970convex}, we know that $\dist(\partial f(\xbf^{*}),\mathbf{0})>0,$ where $\dist(\partial f(\xbf^{*}),\mathbf{0})\coloneqq\min_{\mathbf{\xi}\in\partial f(\xbf^{*})}\|\xi\|$. Furthermore, we make the following assumption:
\begin{assumption}\label{assu:r}
    Throughout the paper, suppose that a constant $r$ satisfying
    \begin{equation}\label{eq:subdiff-lb}
    \dist(\partial f(\xbf^{*}),\mathbf{0})\geq r>0,
    \end{equation}
is known.
\end{assumption}
% derive a subdifferential separation result which will be used repeatedly in our algorithm development. Specifically, for the optimal solution $\xbf^*$, we have 
% \begin{equation}\label{eq:subdiff-lb}
% \dist(\partial f(\xbf^{*}),\mathbf{0})\geq r>0,
% \end{equation}
% for some real value $r$,

We give some important examples for which the lower bound $r$ can be estimated.
Suppose  $f(\xbf)$ is a Lasso regularizer, i.e., $f(\xbf)=\|\xbf\|_{1}$, then $r=1$ satisfies \eqref{eq:subdiff-lb}. More general, consider the   group Lasso regularizer, i.e., $f(\xbf)=\sum_{i=1}^{B}p_{i}\|\xbf_{(i)}\|$, where $\xbf_{(i)}\in\mathbb{R}^{b_{i}}$ and $\sum_{i=1}^{B}b_{i}=n$,
$B$ is the number of blocks, then $r=\min_{i\in [B]}\cbra{p_{i}}$ when $\xbf^{*}\neq\boldsymbol{0}$.
Another example is  $f(\xbf)=\mathbf{c}^{\top}\xbf$,  then we have $r=\|\mathbf{c}\|$.

\begin{remark}
Condition~\eqref{eq:subdiff-lb} is similar to the bounded gradient assumption that has been used for accelerating the convergence of the Frank-Wolfe algorithm. See Appendix~\ref{sec:frank} for more discussions.
% Specifically, assuming that $f(\xbf)$ is smooth and its gradient is uniformly bounded by $\|\nabla f(\xbf)\|>c >0$ for any feasible solution $\xbf$, \cite{dunn1979rates, levitin1966constrained} further improve the convergence rate of Frank-Wolfe algorithm from $\Ocal(1/\vep)$ to $\Ocal\brbra{\log(1/\vep)}$. Nevertheless, the uniform bounded gradient assumption appears to be stronger than ours, as we only impose the lower boundedness assumption on the optimal solution $\xbf^*$ and allow the objective to be non-differentiable. 
\end{remark}

When considering the Lipschitz continuity of function in $\mbb R^n$, even quadratic functions are not Lipschitz continuous. However, the Lipschitz continuity of $g_i(x)$ is crucial for algorithm convergence. Therefore, we define the bounded feasible region in the following proposition, with its proof provided in Appendix~\ref{sec:proof_of_proposition_bounded_set_x}.
\begin{proposition}
\label{prop:bounded_set_x}Let $\mcal X:=\mcal B\brbra{\widetilde{\xbf},\min_{i\in[m]}2\sqrt{\frac{-2g_{i}(\xbf_{i}^{*})}{\mu_{i}}}}$, where $\xbf_{i}^{*}=\argmin_{\xbf\in\Rbb^n}g_{i}(\xbf)$. Then under Assumptions~\ref{assu:Slater's}
and~\ref{assu:dual-optim},  we have $\xbf^{*}\in\intr\mathcal{X}$.
\end{proposition}
% \begin{proof}
%     For complete proof details, please refer to Appendix~\ref{sec:proof_of_proposition_bounded_set_x}.
% \end{proof}
% \begin{proof}
% \proof 
% From the strong convexity of $g_{i}(\xbf)$, we have $g_{i}(\xbf)\geq g_{i}(\xbf_{i}^{*})+\frac{\mu_{i}}{2}\|\xbf-\xbf_{i}^{*}\|^{2},$ which implies
% $$
% \|\xbf-\xbf_{i}^{*}\|^{2}\leq(g_{i}(\xbf)-g_{i}(\xbf_{i}^{*}))\frac{2}{\mu_{i}}\leq\frac{-2g_{i}(\xbf_{i}^{*})}{\mu_{i}}, \forall \xbf \in \mcal{X}_G.
% $$
% In view of the triangle inequality and the above result, for any $\xbf_{1},\xbf_{2}\in \mcal{X}_G$,
% we have
% \[
% \|\xbf_{1}-\xbf_{2}\|  \leq\|\xbf_{1}-\xbf_{i}^{*}\|+\|\xbf_{2}-\xbf_{i}^{*}\|\leq2\sqrt{\frac{-2g_{i}(\xbf_{i}^{*})}{\mu_{i}}}.
% \]
% Hence, $\xbf^*\in \mathcal{B}\Brbra{\widetilde{\xbf},\min_{i\in[m]}2\sqrt{\frac{-2g_{i}(\xbf_{i}^{*})}{\mu_{i}}}}$. Since  $\zeta$ is a positive constant,   we have $\xbf^*\in \textup{int }\mcal{X}$.
% % \endproof
% \endproof

\begin{assumption}
\label{assu:lipschitz}There exist $L_{X},L_{G}>0$ such that
\begin{align}
\|\nabla G(\xbf)-\nabla G(\bar{\xbf})\| & \leq L_{X}\|\xbf-\bar{\xbf}\|,\ \ \forall\xbf,\bar{\xbf}\in\mcal X,\label{eq:lipschitz_x}\\
\|G(\xbf)-G(\bar{\xbf})\| & \le L_{G}\|\xbf-\bar{\xbf}\|,\ \ \forall\xbf,\bar{\xbf}\in\mathcal{X},\label{eq:lipschitz_func}
\end{align}
where $\nabla G(\xbf):=[\nabla g_{1}(\xbf),\cdots,\nabla g_{m}(\xbf)]\in\mathbb{R}^{n\times m}$ and $\mcal X$ is defined in Proposition~\ref{prop:bounded_set_x}.
\end{assumption}

% Proposition~\ref{prop:bounded_set_x} and Proposition~\ref{prop:bounded_set_y} guarantee that the primal and dual optimal solutions are in an set $\Xcal$ and $\Ycal$, respectively, thereby ensuring the  Lipschitzness of Lagrangian function in the domain.
The Lipschitz smoothness of the Lagrangian function with respect to the primal variable $\xbf$ is crucial for the convergence of algorithms. Given that the dual variable $\ybf$ is bounded from above, and considering the smoothness of the constraint functions, we can derive the smoothness of the Lagrangian function.
Combining~\eqref{eq:lipschitz_x} and the
fact $\|\ybf\|\leq\|\ybf\|_{1}\leq \bar{c},\forall \ybf \in \mcal{Y}$,
we obtain that
\begin{equation}
\|\nabla G(\xbf)\ybf-\nabla G(\bar{\xbf})\ybf\|\leq L_{XY}\|\xbf-\bar{\xbf}\|\ \ \forall\xbf,\bar{\xbf}\in\mathcal{X},\ \forall\ybf\in \mathcal{Y},\label{eq:lipschitz_xy}
\end{equation}
 where $L_{XY}=\bar{c} L_{X}$.  For set $\mcal {X}$, $\mcal{Y}$, we use $D_X$ and $D_Y$ to denote their diameters, respectively, i.e., $D_{X}:=\max_{\xbf_1,\xbf_2\in \mcal X}\|\xbf_1 - \xbf_2\|$ and $D_{Y}:= \max_{\ybf_1, \ybf_2\in \mcal{Y}} \|\ybf_1 - \ybf_2\|$.

{\footnotesize
\begin{algorithm}[h]
\caption{\label{alg:APDSCE}\uline{A}ccelerated \uline{P}rimal-\uline{D}ual
Algorithm with \uline{P}rogressive St\uline{ro}ng Convexity
Estimation (APDPro)}
    \begin{algorithmic}[1]
        \Require{$\tau_0>0, \sigma_0>0$, $\mathbf{x}_0\in\mathcal{X}, \mathbf{y}_0\in \mathcal{Y},  \rho_0\geq 0, N>0$}
        \State{\textbf{Initialize:} $\left(\mathbf{x}_{-1}, \mathbf{y}_{-1}\right) \leftarrow\left(\mathbf{x}_{0}, \mathbf{y}_{0}\right), \bar{\mathbf{x}}_{0}\leftarrow \mathbf{x}_0$, $ \sigma_{-1} \leftarrow \sigma_{0},T_0 = 0$}
        \State{Set $\Delta_{XY}=\frac{1}{2\tau_{0}}D_{X}^{2}+\frac{1}{2\sigma_{0}}D_{Y}^{2}$}
        \For{$k = 0,1,\ldots, N$ }
        \State{$\mathcal{Y}_k\leftarrow\left\{\mathbf{y}\in \mathbb{R}_{+}^m\mid \norm{\mathbf{y}}_1 \cdot \underline{\mu} \geq \rho_k\right\}\bigcap \mathcal{Y}$,\label{lst:line:cutting_plane}}
        \State{$\zbf_k\leftarrow {(1+{\sigma_{k-1}}/{\sigma_{k}})G(\mathbf{x}_k)- ({\sigma_{k-1}}/{\sigma_{k}})G(\mathbf{x}_{k-1})}$}
\State{$\ybf_{k+1}\leftarrow \argmin_{\ybf\in\mcal Y_k}\|\ybf - (\ybf_k + \sigma_k\zbf_k)\|^2 $\label{lst:line:updateyintpapd}}
\State{$\mathbf{x}_{k+1}\leftarrow \prox_{f,\mcal{X}} (\mathbf{x}_k-\tau_{k}\nabla G(\mathbf{x}_k)\mathbf{y}_{k+1},\tau_{k})$\label{lst:line:x}}
\State{Compute $t_k$, $\quad \bar{\mathbf{x}}_{k+1}\leftarrow (T_k\bar{\mathbf{x}}_k+t_k \mathbf{x}_{k+1})/(T_k + t_k)$,$\quad T_{k+1}\leftarrow T_k + t_k$\label{lst:line:weighted_average}}

\State{Update $\rho_{k+1}\leftarrow$ \Call{Improve}{$\xbf_k$, $\bar{\xbf}_k$, $\frac{\sigma_0\tau_{k-1}\Delta_{XY}}{\sigma_{k-1}}$, $\frac{\Delta_{XY}}{T_k}$, $\rho_k$}\label{lst:line:maxrho}}
% \State{$\tau_{k+1} \leftarrow \tau_k/\sqrt{1+\rho_{k+1}\tau_k}$, $\quad \sigma_{k+1} \leftarrow \sigma_k\cdot \tau_k / \tau_{k+1}$\label{lst:line:upgamma}}
\State{Update $\tau_{k+1}$ and $\sigma_{k+1}$ depending on $\rho_{k+1}$\label{lst:line:upgamma}}
\EndFor
    \State{\textbf{Output:}\quad $\mathbf{x}_{N+1},\mathbf{y}_{N+1}$}
    \Procedure{Improve}{$\xbf$, $\bar\xbf$, $\beta$, $\bar\beta$,  $\rho_{\text{old}}$}
    \State{{Compute} 
    $\rho  =\underline{\mu}\cdot\max\bigg\{r\bsbra{\norm{\nabla G({\xbf})}+L_X\sqrt{2\beta}}^{-1},  \Bsbra{\frac{L_X}{r}\sqrt{\frac{\bar\beta}{2\underline{\mu}}}
    +\sqrt{\frac{L_X^2\bar\beta }{2\underline{\mu}r^2}+\frac{\norm{\nabla G(\bar\xbf)}}{r}}}^{-2}\bigg\}$
    \label{compute-rho-hat}
    }
    \State{{Set} $\rho_\text{new} = \max\{\rho_\text{old}, \rho\}$}
    \State{\Return $\rho_{\text{new}}$}
    \EndProcedure \label{improve-end}
    \end{algorithmic}
\end{algorithm}
}

\section{\label{sec:APDSCE}\APD{} with progressive
strong convexity estimation}
We present the Accelerated Primal-Dual
Algorithm with Progressive Strong Convexity Estimation ({\APDSCE}) to solve problem~\eqref{eq:constrainedproblem}.
% In the algorithm running process,
% we apply a new technique to improve the estimated strong convexity modulus of the Lagrangian function.
For problem~\eqref{eq:constrainedproblem}, {\APDSCE}  achieves the improved convergence rate $\mcal O(1/\sqrt{\vep})$ without relying on the uniform strong convexity assumption~\citep{hamedani2021primal,lin2020near}. 
For the rest of this paper, we denote $\prox_{f,\mcal{X}}(\xbf-\eta\zbf,\eta)\coloneqq\argmin_{\hat{\xbf}\in\mcal X}f(\hat{\xbf})+\left\langle \zbf,\hat{\xbf}\right\rangle +\tfrac{1}{2\eta}\|\hat{\xbf}-\xbf\|^{2}$
as the proximal mapping.
% associated with $f(\cdot) + \inner{\zbf}{\cdot}$ on $\Xcal$.

 We describe {\APDSCE} in Algorithm~\ref{alg:APDSCE}.
The main component of {\APDSCE} contains a dual ascent step to update $\ybf_k$ based on the extrapolated gradient, followed by a primal proximal step to update $\xbf_k$.
Compared with standard \APD{}~\citep{hamedani2021primal}, {\APDSCE} has two more steps. 
First, line~\ref{lst:line:cutting_plane} of Algorithm~\ref{alg:APDSCE} applies a novel cut constraint to separate the dual sequence $\{\ybf_k\}$ from the origin, which allows us to leverage the strong convexity of the Lagrangian function and hence obtain a faster rate of convergence than \APD{}. Second, to use the strong convexity more effectively, in line~\ref{lst:line:maxrho}, we perform a progressive estimation of the strong convexity by using the latest iterates $\xbf_k$ and  $\bar{\xbf}_k$. Throughout the algorithm process, we use a routine  {\IMPROVE} to construct a non-decreasing sequence $\{\rho_k\}$, which provides increasingly refined lower bounds of the strong convexity of the Lagrangian function.

\noindent \textbf{The \IMPROVE{} step}
In order to estimate the strong convexity of the Lagrangian function, we rely on the subdifferential separation  (eq.~\eqref{eq:subdiff-lb}) to bound the dual variables. 
From the first-order optimality condition in minimizing $\Lcal(\xbf, \ybf^*)$ 
and the fact that $\xbf^{*}\in\intr\mcal X$ (Proposition~\ref{prop:bounded_set_x}), we have
$
\boldsymbol{0}\in\partial f(\xbf^{*})+\nabla G(\xbf^{*})\ybf^{*}+\mcal N_\Xcal(\xbf^*)=\partial f(\xbf^{*})+\nabla G(\xbf^{*})\ybf^{*}.
$ 
It follows from \eqref{eq:subdiff-lb} that
\begin{equation}\label{eq:r_lower}
    r\leq\|\nabla G(\xbf^{*})\ybf^{*}\| \le \|\nabla G(\xbf^*)\|\cdot\|\ybf^*\| \leq \norm{\ybf^*}_1\|\nabla G(\xbf^{*})\|,
\end{equation}
where the last inequality use the fact that $\norm{\cdot} \le \norm{\cdot}_1$. Note that the bound $\norm{\ybf^*}_1\ge {r}/{\norm{\nabla G(\xbf^*)}}$ can not be readily used in the algorithm implementation because  $\xbf^*$ is generally unknown. To resolve this issue, we develop more concrete dual lower bounds by using the generated solution $\hat\xbf$ in the proximity of $\xbf^*$. As we will show in the analysis, {{\APDSCE}} keeps track of two primal sequences $\{\xbf_k\}$ and $\{\bar\xbf_k\}$, for which we can establish bounds on $\norm{\xbf_k-\xbf^*}^2$ and ${(\ybf^*)^{\top}\mubm}\cdot \norm{\hat\xbf-\xbf^*}^2/2$, respectively. This drives us to develop the following lower bound property, with the proof provided in Appendix~\ref{sec:proof_of_propostion_local_estimation}.
\begin{proposition}\label{prop:local_estimation}Suppose Assumption~\ref{assu:lipschitz} holds.
Let  $\ybf^*\in\Ycal^*$  be a dual optimal solution. 
\begin{enumerate}
\item Suppose that $\|\hat\xbf-\xbf^*\|^2 \le 2\beta$, then we have 
\begin{equation}\label{eq:y-lb-1}
\|\ybf^*\|_1  \ge h_1(\hat{\xbf}, \beta)\coloneqq r\bsbra{\norm{\nabla G(\hat{\xbf})}+L_X\sqrt{2\beta}}^{-1}.
\end{equation}
\item Suppose $(\ybf^*)^{\top}\mubm\cdot \|\hat\xbf-\xbf^*\|^2 \le 2 \beta$,  then we have
\begin{equation}\label{eq:y-lb-2}
    \|\ybf^*\|_1 \ge h_2(\hat{\xbf}, \beta) := \lrsbra{\tfrac{L_X}{r}\sqrt{\tfrac{\beta}{2\underline{\mu}}}
    +\sqrt{\tfrac{L_X^2\beta }{2\underline{\mu}r^2}+\tfrac{\norm{\nabla G(\hat\xbf)}}{r}}}^{-2}.
\end{equation}
\end{enumerate}
\end{proposition}
% \begin{proof}
%     For complete proof details, please refer to Appendix~\ref{sec:proof_of_propostion_local_estimation}.
% \end{proof}
% \proof{}
% Moreover, using
% the triangle inequality and~\eqref{eq:lipschitz_x}, we have
% \[
% \|\nabla G(\xbf^{*})\|-\|\nabla G(\hat\xbf)\|\le\|\nabla G(\xbf^{*})-\nabla G(\hat\xbf)\|\leq L_{X}\|\hat\xbf-\xbf^{*}\|.
% \]
% Combining the above inequality and~\eqref{eq:r_lower}, we obtain
% \begin{equation}\label{eq:mid-04}
% \tfrac{r}{\norm{\ybf^*}_1}\le L_X \|\hat\xbf-\xbf^{*}\|+\|\nabla G(\hat\xbf)\|.
% \end{equation}
% Next,  we develop more specific lower bounds on $\norm{\ybf}_1$.  i). Inequality~\eqref{eq:y-lb-1} can be easily verified since  we have $\|\hat\xbf-\xbf^*\| \le \sqrt{2\beta}$.  
% ii). Suppose ${(\ybf^*)^{\top}\mubm}\cdot \|\hat\xbf-\xbf^*\|^2 \le 2\beta$, then together with \eqref{eq:mid-04}  we have 
% $$
% \tfrac{r}{\norm{\ybf^{*}}_1}\le L_X \sqrt{\tfrac{2\beta}{(\ybf^*)^{\top}\mubm}}+\norm{\nabla G(\hat\xbf)}\le L_X \sqrt{\tfrac{2\beta}{\underline{\mu}\norm{\ybf^*}_1}}+\norm{\nabla G(\hat\xbf)}.
% $$
% Note that the above inequality can be expressed as $at^2-bt-c \le 0$ with $t=\norm{\ybf^*}_1^{-1/2}$, $a=r, b=L_X\sqrt{{2\beta}/{\underline{\mu}}}$ and $c=\norm{\nabla G(\hat\xbf)}$. Standard analysis implies that $t\le (b+\sqrt{b^2+4ac})/{2a}$, which gives the desired bound~\eqref{eq:y-lb-2}.
% \endproof

Our next goal is to conduct the convergence analysis for {\APDSCE} in Theorem~\ref{thm:gapdiminishing} and Corollary~\ref{cor:mainthm}. Complete proof details are provided in Appendix~\ref{sec:proof_of_thm_gapdiminishing} and~\ref{proof:cor}. 
% which can be divided into two parts. First, under the two preconditions~\eqref{eq:param-01}
% and~\eqref{eq:param-03}, we present the convergence conclusion for the primal-dual gap and optimality gap of problem~\eqref{eq:constrainedproblem}
% (see Theorem~\ref{thm:gapdiminishing}). 
% Next, we derived
% the convergence rate of {\APDSCE} for a specific 
%  stepsize selection (see Corollary~\ref{cor:mainthm}). 

\begin{theorem}
\label{thm:gapdiminishing}Suppose for any $\ybf^{*}\in\mcal Y^{*}$, $(\ybf^{*})^{\top}\boldsymbol{\mu}\geq\rho_{0}$ holds, and let the sequence $\{\tau_k,\sigma_k,t_k,\rho_{k+1}\}$ generated by Algorithm~\ref{alg:APDSCE} satisfy:
\begin{align}
t_{k+1}(\tau_{k+1}^{-1}-\rho_{k+1})\leq t_{k}\tau_{k}^{-1},\quad t_{k+1}\sigma_{k+1}^{-1}\leq t_{k}\sigma_{k}^{-1},\quad L_{XY}+L_{G}^{2} \sigma_{k} \leq\tau_{k}^{-1}.\label{eq:param-03}
\end{align}
Then, the set $\Ycal_k$ is nonempty and $\Ycal^*\subseteq\Ycal_k$. Let  $\Delta(\xbf, \ybf)\coloneqq \tfrac{1}{2\tau_{0}} \|\xbf - \xbf_0\|^2 + \tfrac{1}{2\sigma_{0}} \|\ybf - \ybf_0\|^2,\bar{\ybf}_{K}=T_{K}^{-1}\sum_{s=0}^{K-1}t_{s}\ybf_{s}$.  The sequence $\{\bar{\xbf}_{k},\xbf_k,\bar{\ybf}_k \}$ generated by~{\APDSCE} satisfies
\begin{equation}
\label{eq:apd_primal_dual_converge}
    \tfrac{t_{K-1}\tau_{K-1}^{-1}}{2T_{K}}\|\xbf^{*}-\xbf_{K}\|^{2}+\mathcal{L}(\bar{\xbf}_{K},\ybf^{*})-\mathcal{L}(\xbf^{*},\bar{\ybf}_{K})\leq\tfrac{1}{T_{K}}\Delta(\xbf^*, \ybf^*).
\end{equation}
\end{theorem}
% \begin{proof}
%     For complete proof details, please refer to Appendix~\ref{sec:proof_of_thm_gapdiminishing}.
% \end{proof}

Next, we develop more concrete complexity results in Corollary~\ref{cor:mainthm}.

\begin{corollary}
\label{cor:mainthm}
Suppose that $\sigma_k,\tau_k,t_k$ satisfy:
\begin{equation}\label{eq:cor_param}
\begin{split}\tau_{0}^{-1}\geq L_{XY}+L_{G}^{2}\sigma_{0}, &\ \  t_{k}=\sigma_{k}/\sigma_{0},\\
\tau_{k+1}=\tau_{k}/\sqrt{1+\rho_{k+1}\tau_k},& \ \ \sigma_{k+1}=\sigma_{k}\tau_{k}/\tau_{k+1}
\end{split}    
\end{equation}
Then we have
\begin{equation}\label{eq:cor_complexity}
\begin{aligned}
f(\bar{\xbf}_{K})-f(\xbf^{*})
&\leq  \tfrac{6}{6+\tau_{0}\trho_{K} (K+1)K}\Brbra{\tfrac{1}{2\tau_0}\norm{\xbf_0-\xbf^*}^2+\tfrac{D_Y^2}{2\sigma_0}},\\
\|[G(\bar{\xbf}_{K})]_{+}\| & \leq \tfrac{6}{c^{*}(6+\tau_{0}\trho_{K} (K+1)K)}\Brbra{\tfrac{1}{2\tau_0}\norm{\xbf_0-\xbf^*}^2+\tfrac{D_Y^2}{2\sigma_0}},\\
\tfrac{1}{2}\|\xbf_{K}-\xbf^{*}\|^{2}&\leq\tfrac{3\sigma_{0}}{\vrho_{K}^{2}\tau_{0}^{2}K^{2}+9(\sigma_0/\tau_0)}\Delta(\xbf^*,\ybf^*).
\end{aligned}
\end{equation}
where $c^{*}:=\brbra{f(\xbf^{*})-\min_{\xbf}f(\xbf)}/{\min_{i\in[m]}\{-g_{i}(\widetilde{\xbf})\}}>0$, $\tilde{\rho}_k = 2\sum_{s=0}^{k}\hat{\rho}_s s/\brbra{k(k+1)}$ and $\tilde{\rho}_k$ satisfy the following condition, $\hat{\rho}_{k+1}:=\sqrt{\hat{\rho}_k^2 k^2 + (3\rho_{k+1}\hat{\rho}_k)k}/(k+1), \hat{\rho}_1 = 3\sqrt{\rho_1/\tau_0}$.
\end{corollary}
% \begin{proof}
%     For complete proof details, please refer to Appendix~\ref{proof:cor}.
% \end{proof}

\begin{remark}\label{rem:rate}
In view of Corollary~\ref{cor:mainthm}, \APDSCE{} obtains an iteration complexity of  $\Ocal({1}/{\sqrt{\trho_{K}\vep}})$, which is substantially better than the $\Ocal({1}/{\vep})$ bound of APD~\citep{hamedani2021primal} and ConEx~\citep{boob2019stochastic} when the strong convexity parameter $\tilde{\rho}_K$ is relatively large compared with $\vep$. 
\end{remark}
\begin{remark}
Additionally, we argue that even when $\trho_{K} =O(\vep)$, \APDSCE{} can obtain the matching $\Ocal({1}/{\vep})$ bound of the state-of-the-art algorithms. Specifically, using the definition of $\sigma_k,\tau_k$, we can easily derive the monotonicity of $\{\sigma_k\}$.
It follows from 
$
\sigma_{k+1}=\tau_k\sigma_k/\tau_{k+1}=\tau_k\sigma_k/\brbra{\tau_k/\sqrt{1+\rho_{k+1}\tau_k}}\geq\sigma_k,
$
that $T_k = \sum_{s=0}^{k-1} t_k  = \sigma_0^{-1}\sum_{s=0}^{k-1}\sigma_k \ge k$. Using a similar argument to that of Corollary~\ref{cor:mainthm}, we obtain the bound $f(\bar{\xbf}_{K})-f(\xbf^{*})\le \Ocal({1}/{K})$ and $\|[G(\bar{\xbf}_{K})]_{+}\|\le \Ocal({1}/{K})$.

\end{remark}

\begin{remark}
    The  implementation of \APDSCE{} requires knowing an upper bound on $\|\ybf^*\|$. When the bound is unavailable, \cite{hamedani2021primal} developed an adaptive \APD{} which still ensures the boundedness of dual sequence via line search. Since our main goal of this paper is to exploit the \emph{lower-bound} rather than the \emph{upper bound} of $\|\ybf^*\|$, we leave the extension for the future work.
\end{remark}

\begin{algorithm}[h]
\caption{\label{alg:Restart Version of Stage Two TPAPD}Restarted \APDSCE{} (\resapdsce)}
\begin{algorithmic}[1]
\Require{$\rho^{-1}_{N_{-1}}\geq 0,\bar{\sigma}>0$,  {$\nu_{0}\in(0,1)$}, $\delta\in(0,1)$, $\xbf^{-1}_{N_{-1}}, \ybf^{-1}_{N_{-1}}, S$ }
\State{Compute  {$\bar{\tau}=(1-\nu_{0})\brbra{L_{XY}+L_{G}^{2}\bar{\sigma}/\delta}^{-1}$}}
        \For{$s=0,1,\ldots,S$}
        \State{$\tau_{0}^s = \bar{\tau}, \sigma_{0}^s = \bar{\sigma}, (\xbf_{-1}^s, \ybf_{-1}^s)\leftarrow (\xbf_{N_{s-1}}^{s-1}, \ybf_{N_{s-1}}^{s-1}),(\xbf_{0}^s, \ybf_{0}^s)\leftarrow (\xbf_{N_{s-1}}^{s-1}, \ybf_{N_{s-1}}^{s-1}), \rho_{0}^{s} = \rho_{N_{s-1}}^{s-1}$}
        \State{Set $\Delta_{XY} = \frac{1}{\tau_{0}^s}D_{X}^2 + \frac{1}{2\sigma_{0}^s}D_{Y}^2,\sigma_{-1}^{s}\leftarrow \sigma_{0}^{s}, T_0^s = 0, k=0, \hat{\rho}_{0}^s = 1, N_s=\infty$}
        \While{$k<N_s$}
        \State{Run line~\ref{lst:line:cutting_plane}-\ref{lst:line:upgamma} of {\APDSCE}  with index set $(s,k)$}

 \State{Update $N_s,\vrho_{k+1}^{s}\leftarrow$\Call{TerminateIter}{$\vrho_{k}^{s}, \rho_{k+1}^{s}, s, k$}, $k\leftarrow k + 1$
 
 }
        \EndWhile
        \EndFor
\State{\textbf{Output: }$\mathbf{x}^{S}_{N_{S}},\mathbf{y}^{S}_{N_{S}}$}
    \Procedure{TerminateIter}{$\vrho_{\text{old}},\rho,s, k$}
    \State{
    Compute $\vrho_{\text{new}} =
    \begin{cases}
    \frac{1}{k+1}\sqrt{\vrho_{\text{old}}^{2}k^{2}+3\rho\vrho_{\text{old}} k} & k >1 \\
    3\sqrt{{\rho}/{\tau_{0}}} & k = 1
    \end{cases}
    $
    } 
    \State{Compute $N=\lceil\max\{6(\hat{\rho}_{\text{new}}\tau_{0}^{s})^{-1},\sqrt{2}^{s}\cdot3\sqrt{2}D_{Y}/\brbra{\hat{\rho}_{\text{new}}D_{X}\sqrt{\tau_{0}^{s}\sigma_{0}^{s}}}\}\rceil 
    % \lceil\max\{\frac{2\sqrt{3{\sigma_0}}}{\vrho_{\text{new}}\bar{\tau}^{1.5}},    \frac{\sqrt{6}D_{Y}}{\vrho_{\text{new}}\tau_{0}D_{X}}2^{\frac{s}{2}}\}\rceil
    $}
    \State{\Return $N,\vrho_{\text{new}}$}
    \EndProcedure
    \end{algorithmic}
\end{algorithm}

\section{\label{sec:restarted}{\APDSCE} with a restart scheme}
Note that in the worst case,  \APDSCE{} exhibits an iteration complexity of  $\Ocal\brbra{(D_X+D_Y)/\sqrt{\vep}}$, which has a linear dependence on the diameter. 
 While the $\Ocal({1}/{\sqrt{\vep}})$ is optimal~\citep{ouyang2021lower}, it is possible to improve the complexity with respect to the primal part from $\Ocal\brbra{{D_X}/{\sqrt{\varepsilon}}}$ to $\Ocal\brbra{\log\brbra{{D_X}/{\sqrt{\varepsilon}}}}$.
To achieve this goal, we propose  a restart scheme ({\resapdsce})
that calls {\APDSCE} repeatedly and present the details in  Algorithm~\ref{alg:Restart Version of Stage Two TPAPD}. 
% The total iteration number of \APDSCE{} in each epoch is an important parameter of the restart algorithm.
Inspired by~\cite{lan2020firstorder}, we set the iteration number as a function of the estimated strong convexity, detailed in the {\UpdateN} procedure. For convenience in describing a double-loop algorithm,  we use superscripts for the number of
epochs and subscripts for the number of sub-iterations in parameters
$\xbf,\ybf,\tau,\sigma$, e.g., $\xbf_{1}^{S}$
meaning the $\xbf$ output of first iterations at $S$-th epoch. {To avoid redundancy in the Algorithm~\ref{alg:Restart Version of Stage Two TPAPD}, we call the {\APDSCE} iteration directly. Note that the notation system here is identical to that of {\APDSCE}, with the only difference being the use of superscripts to distinguish the number of epochs.}

In Theorem~\ref{thm:rapdpro}, we show the overall convergence complexity of {\resapdsce} with the proof provided in Appendix~\ref{sec:proof_thm_rapdpro}.

\begin{theorem}
\label{thm:rapdpro}
Let $\{\xbf_{0}^{s}\}_{s\geq0}$
be the sequence generated by {\resapdsce}, then we have
\begin{equation}
\|\xbf_{0}^{s}-\xbf^{*}\|^{2}\leq\Delta_{s}\equiv D_{X}^{2}\cdot 2^{-s},\ \ \ \forall s\geq0.\label{eq:exp decay}
\end{equation}
As a consequence, {\resapdsce}
will find a solution $\xbf_{0}^{S}$ such that $\|\xbf_{0}^{S}-\xbf^{*}\|^{2}\leq\vep$
for any $\vep\in(0, D_{X}^{2})$ in at most $S:=\big\lceil \log_{2}\rbra{ D_{X}^{2}/\vep}\big\rceil $
epochs. Moreover, 
The iteration number of {\resapdsce}
to find $\xbf_{0}^{S}$ such that $\|\xbf_{0}^{S}-\xbf^{*}\|^{2}\leq\vep$ is bounded by
\begin{equation}\label{T_e}
    T_{\vep}:= \brbra{\tfrac{12}{\varpi_{1}\tau_{0}^{s}}+2}\left\lceil \log_{2}\tfrac{D_{X}}{\sqrt{\vep}}+1\right\rceil +\brbra{\tfrac{6(\sqrt{2}+2)}{\varpi_{2}\sqrt{\tau_{0}^{s}\sigma_{0}^{s}}}}\cdot\brbra{\tfrac{D_{Y}}{\sqrt{\vep}}},
\end{equation}
where $\varpi_1$ and $\varpi_2$ satisfy $\sum_{s=0}^{S}(\hat{\rho}_{N_{s}}^s)^{-1} = (\varpi_{1})^{-1}(S+1)$ and $\sum_{s=0}^{S} {\sqrt{2}^s}/{\hat{\rho}_{N_{s}}^s}= (\varpi_2)^{-1}\sum_{s=0}^{S}\sqrt{2}^s$, respectively.
\end{theorem}
% \begin{proof}
%     For complete proof details, please refer to Appendix~\ref{sec:proof_thm_rapdpro}.
% \end{proof}

\begin{remark}
The bound $T_{\vep}$ depends on $\vep$, $\varpi_1$ and $\varpi_2$. If $\varpi_1 = O\big((-\log_{2}\sqrt{\vep})^{-1}\big)$ or $\varpi_2 = O(\sqrt{\vep})$, then we have $T_\vep = \infty$, which implies that we can not guarantee $\norm{\xbf_0^s-\xbf^*} \le \vep$ at finite iterations. $T_{\vep}=\infty$ implies that there exists an epoch with infinite sub-iterations. Hence, {\resapdsce} is reduced to {\APDSCE} if we only consider that epoch.
\end{remark}
\begin{remark}Comparison of \resapdsce{} and \APDSCE{} involves a number of factors. 
In particular, \resapdsce{} compares favorably against \APDSCE{} if  $\|\xbf_0-\xbf^*\|=\widetilde\Omega (\sqrt{\vep}\log D_X)$. Moreover,  the complexity~\eqref{T_e} can be slightly improved if $D_X$ is replaced by any tighter upper bound of $\|\xbf^s_0-\xbf^*\|$. However, it is still unknown whether we can directly replace $D_X$ with $\|\xbf^s_0-\xbf^*\|$ in \eqref{T_e}.
\end{remark}
\noindent \textbf{Dual Convergence}
For dual variables,  we establish asymptotic convergence to the optimal solution, a key condition for developing the active-set identification in the later section. For ease in notation, it is more convenient to label the generated solution as a whole sequence using a single subscript index: $\xbf_1,\xbf_2,\ldots, \xbf_N; \ybf_1, \ybf_2, \ldots, \ybf_N$. Hence, we use the index system $j$ and $(s,k)$ interchangeably.  Note that  $\{\xbf_{0}^{s+1},\ybf_{0}^{s+1}\}$ and $\{\xbf_{N_{s}+1}^{s},\ybf_{N_{s}+1}^{s}\}$ correspond to  the same pair of points.  We present the dual asymptotic result in the following theorem, with the proof provided in Appendix~\ref{proof_of_thm_asymptotic}.

\begin{theorem}\label{thm:asymptotic}
Assume $\bar{\tau}^{-1}>\overline{\rho}$
and  choose  $\nu_{0}>0$ such that 
$1>\inf_{j\ge 0}\{\sigma_{j-1}/\sigma_{j}\}\geq \delta + \nu_{0}.$
We have $(\xbf^{*},\ybf^{*})$
satisfy the KKT condition, where $\ybf^{*}$ is any limit point of $\{\ybf_{j}\}$
generated by {\resapdsce}. 
% (\textcolor{red}{For dual asymptotic convergence, we need to additional introduce $\delta\in(0,1)$ here. However, $\delta=1$ in APDPro.})
\end{theorem}
% \begin{proof}
% For complete proof details, please refer to Appendix~\ref{proof_of_thm_asymptotic}.
% \end{proof}
\begin{remark}
To establish the asymptotic convergence of the dual variable, we introduce an additional constant $\delta\in(0,1)$, which implies that the initial step size must meet a stricter requirement than the convergence condition specified in Corollary~\ref{cor:mainthm}.
Since ${\sigma_{k}^s}/{\sigma_{k-1}^s}=\sqrt{1+\rho_{k}^{s}\tau_{k}^{s}}$,
$\{\rho_{k}^{s}\}$ is bounded due to the boundedness of the dual variable, $\{\tau_{k}^{s}\}$ 
is monotonically decreasing, then $\inf_{0 \le k \le N_s}\{\sigma_{k-1}^{s}/\sigma_{k}^{s}\}\geq (1+\overline{\rho}\bar{\tau})^{-1/2}$. Hence, inequality, $1>\inf_{j\ge 0}\{\sigma_{j-1}/\sigma_{j}\}\geq\delta+\nu_{0}$, is always satisfiable if we choose proper $\delta, \nu_{0}$ such that $(1+\overline{\rho}\bar{\tau})^{-1/2}\geq \delta+\nu_{0}$. Furthermore, Assumption $(\bar{\tau})^{-1}>\overline{\rho}$ is mild. Since we always choose $\bar{\sigma}$ large enough in {\resapdsce}, $\bar{\tau}$ can be sufficiently small.
\end{remark}
\begin{remark}
    Both algorithms proposed previously require solving quadratic optimization with linear constraints when updating dual variables, which may introduce implementation overheads when the constraint number is high. Inspired by the multi-stage algorithm, we additionally propose an algorithm (Multi-Stage APD, msAPD) that uses different step sizes in different stages and dynamically adjusts the number of iterations in each stage by leveraging strong convexity, as detailed in Appendix~\ref{sec:MSAPD}.
\end{remark}

\section{\label{sec:Manifold}Active-set identification in sparsity-inducing optimization}

In this section, we apply our proposed algorithms to  the  aforementioned sparse learning problem:
\begin{equation}\label{eq:sparselearn}
\begin{aligned}
\textstyle \min f(\xbf),\,
 \textrm{s.t.}\,\, g(\xbf)\le 0, \  \xbf=\xbf_{(1)}\times \ldots\times\xbf_{(B)},\  \xbf_{(i)}\in\mathbb{R}^{n_{i}}, 1\le i \le B,
\end{aligned}
\end{equation}
where $f(\xbf)=\sum_{i=1}^{B}p_{i}\|\xbf_{(i)}\|$ is the group Lasso regularizer and $g(\xbf)$ is a strongly convex function.  We use $\xbf_{(i)}$ to express  the $i$-th block coordinates of $\xbf$. 
% For notation convenience, we write $f_{i}(\xbf_{(i)})=p_{i}\|\xbf_{(i)}\|$. 
% A notable example is the Lasso problem where $B=n,p_{i}=1,\forall i$, and $\xbf_{(i)}$ is a real value and we write $f(\xbf)=\sum_{i}|\xbf_{(i)}|$. 
% As has been shown in the earlier sections, our proposed algorithms \APDSCE{}, \resapdsce{} and \MSAPD{} exhibit an $\Ocal({1}/{\sqrt{\vep}})$ rate of convergence. 
The goal of this section is to show that \resapdsce{} can identify the sparsity pattern of the optimal solution of \eqref{eq:sparselearn} in a finite number of iterations. 
% Identification property of other algorithms can be established analogously.

In general, suppose that  $f(\xbf)$ has a separable structure 
$f(\xbf)=\sum_{i=1}^B f_i(\xbf_{(i)})$, 
we define the active set $\mathcal{A}(\xbf)$ for  $f(\xbf)$  by
$
\mathcal{A}(\xbf):=\{i:\partial f_{i}(\xbf_{(i)})\text{ is not a singleton}\}.
$
For $f(\xbf)=\sum_{i=1}^{B}p_{i}\|\xbf_{(i)}\|$, it is easy to see that $\mathcal{A}(\xbf)$ is the index set of the zero blocks: $\mathcal{A}(\xbf^{*})=\big\{ i:\xbf_{(i)}^{*}=\mathbf{0}\big\}$.
Next, we describe one property for the optimal solution of~\eqref{eq:sparselearn} in Proposition~\ref{prop:uniquey} with the proof provided in Appendix~\ref{proof_prop_uniquey}.
\begin{proposition}\label{prop:uniquey}
Under Assumptions~\ref{assu:Slater's} and~\ref{assu:dual-optim}, the KKT point for~\eqref{eq:sparselearn} is unique. 
\end{proposition}
% \begin{proof}
%     For complete proof details, please refer to Appendix~\ref{proof_prop_uniquey}.
% \end{proof}
% \proof
% The uniqueness of primal
% optimal solution $\xbf^{*}$ follows from Proposition~\ref{prop:uniquex}.
% The KKT condition (ensured by Slater's CQ) implies 
% \begin{equation}\label{eq:kkt-2}
%     \zerobf\in \partial f(\xbf^*) + \nabla g(\xbf^*) \ybf^*.
% \end{equation}
% According to Assumption~\ref{assu:dual-optim}, we have $\xbf^* \neq \boldsymbol{0}$, hence $\Acal^c(\xbf^*)=\{1,2,\ldots,B\}\setminus\Acal(\xbf^*)\neq\emptyset$. In view of \eqref{eq:kkt-2}, for any $i\in\Acal^c(\xbf)$, we have 
% $p_i {\xbf^*_{(i)}}/{\norm{\xbf^*_{(i)}}}=-\nabla_{(i)} g(\xbf^*)\ybf^*$, which gives a unique $\ybf^*$.
% \endproof

To identify the sparsity pattern (active set) of the optimal solution, it is common to assume the existence of a non-degenerate optimal solution, which is stronger than the standard optimality condition~\citep{nutini2019activeset, sun2019we}. We say that $\xbf^*$ is non-degenerate if 
$
{\bf 0} \in \ri \partial\Lcal(\xbf^*, \ybf^*)= \ri(\partial f(\xbf^*) +  \nabla g(\xbf^*) \ybf^*)
$ for the Lagrangian multiplier $\ybf^*$, where $\ri$ stands for the relative interior. 
More specifically, $(\xbf^*,\ybf^*)$ satisfies the block-wise optimality condition
\begin{equation*}
\begin{cases}
-[\nabla g(\xbf^{*})\ybf^{*}]_{(i)}=\nabla f_{i}(\xbf_{(i)}^{*}), & \text{if }\, i\notin\mathcal{A}(\xbf^{*}),\\
-[\nabla g(\xbf^{*})\ybf^{*}]_{(i)}\in\intr \brbra{\partial f_{i}(\xbf_{(i)}^{*})}, & \text{if }\, i\in\mathcal{A}(\xbf^{*}).
\end{cases}
\end{equation*}
Inspired by \cite{nutini2019activeset}, we use the radius $\eta:=\textstyle\min_{i\in\mathcal{A}(\xbf^{*})}\left\{ p_{i}-\norm{[\nabla g(\xbf^*)\ybf^*]_{(i)}} \right\}$, which describes the certain distance between the gradient and "subdifferential boundary" of the active set.  
% The proof strategy of active-set identification is similar to that of unconstrained optimization \citep{nutini2019activeset}. 
We demonstrate in the following theorem that the optimal sparsity pattern is identified when the iterates fall in a neighborhood dependent on $\eta$, with the proof provided in Appendix~\ref{sec:proof_thm_active_set_ide}.

\begin{theorem}\label{thm:active_set_ide}
    Set $\mcal X:=\mcal B\Brbra{\tilde{\xbf},\min_{i\in[m]}2\sqrt{\frac{-2g_i(\xbf_i^*)}{\mu_i}}+\zeta}$ with $\zeta>0$ and $3L_{XY}\cdot(\bar{\tau}+(2L_{XY})^{-1})\cdot \zeta>{\eta}{\bar{\tau}}$ in {\resapdsce}, then we have there exists a epoch $\hat{S}_0$ such that $\xbf^{*}_{(i)} = \xbf_{k(i)}^{s}, s \ge  \hat{S}_0, \ \forall k \in [N_{s}],\ \forall i \in \mcal{A}(\xbf^*).$ 
    % (\textcolor{red}{For active set identification specifically $\xbf^{k+1}_i\in \text{int} \mcal X$, we need to introduce $\zeta$.})
\end{theorem}
\begin{remark}
    The active-set identification result is achieved using the optimality condition at the next iterate $\xbf^{k+1}_i$. To ensure $\xbf^{k+1}_i\in \intr \mcal X$, we define an expanded region, which prevents cases where the normal cone differs from $\{\boldsymbol{0}\}$.
\end{remark}

\begin{figure}[h]
\begin{centering}
\begin{minipage}[t]{0.33\columnwidth}%
      \includegraphics[width=4.5cm]{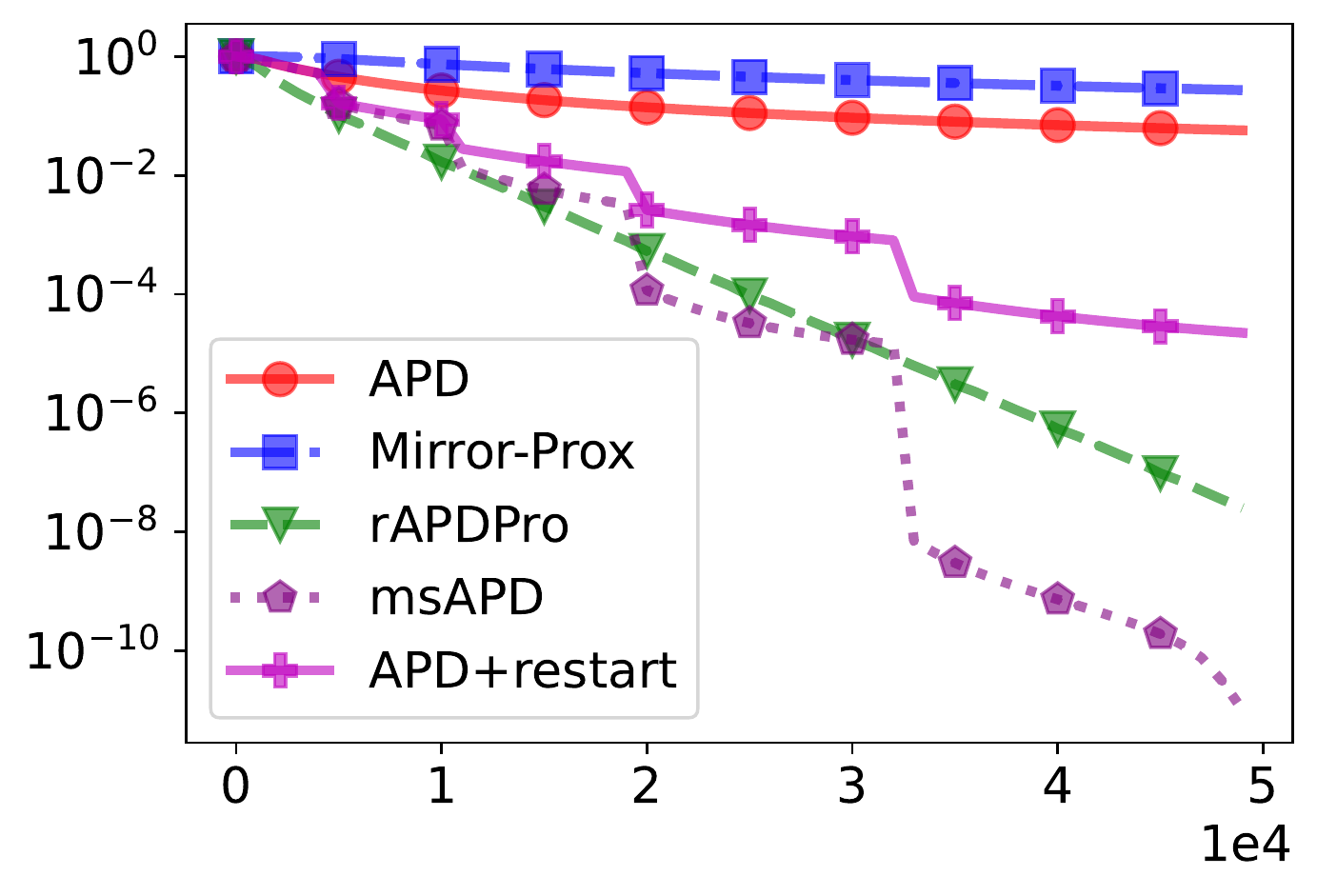}%
\end{minipage}%
\begin{minipage}[t]{0.33\columnwidth}%
\includegraphics[width=4.5cm]{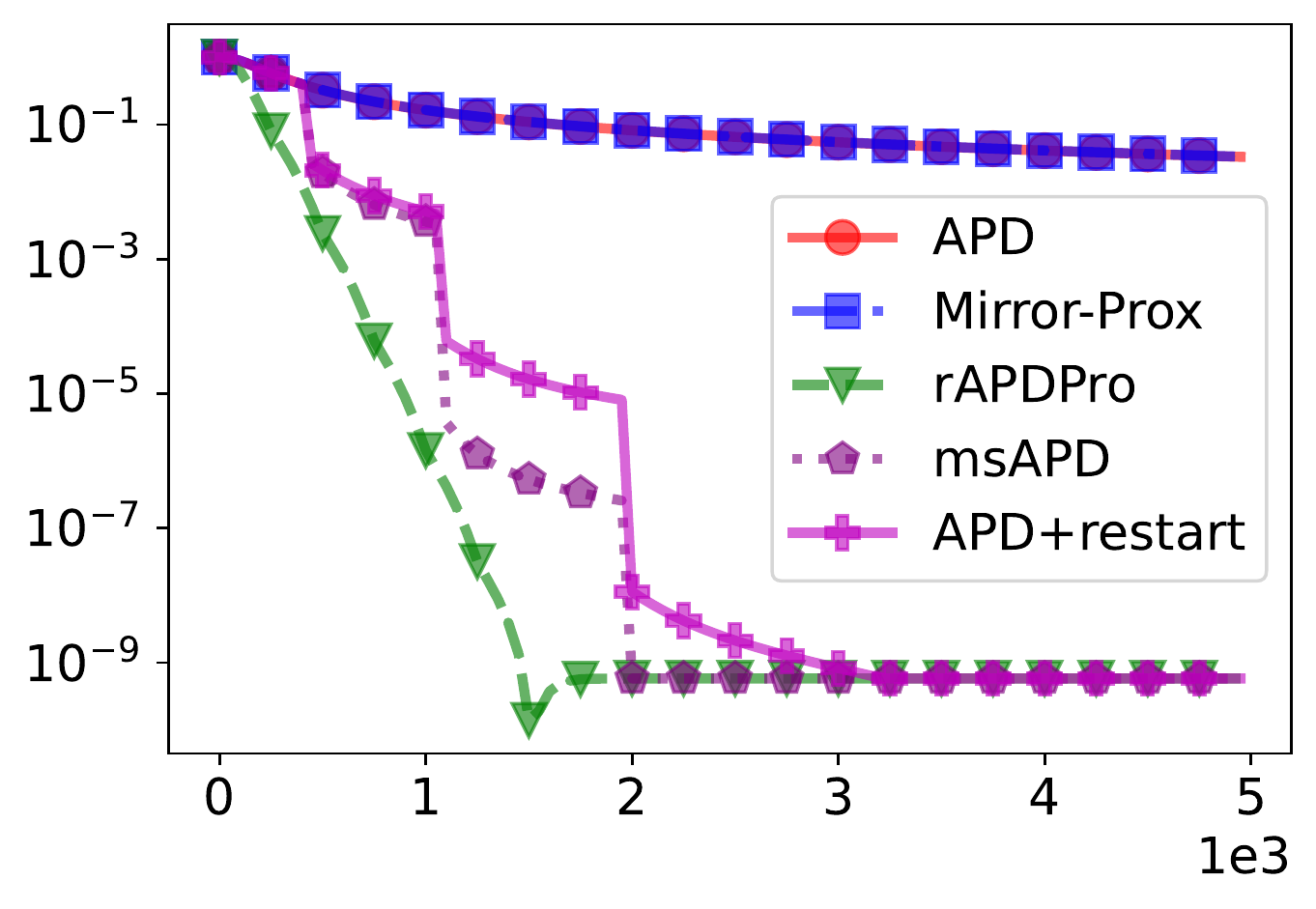}%
\end{minipage}%
\begin{minipage}[t]{0.33\columnwidth}%
\includegraphics[width=4.5cm]{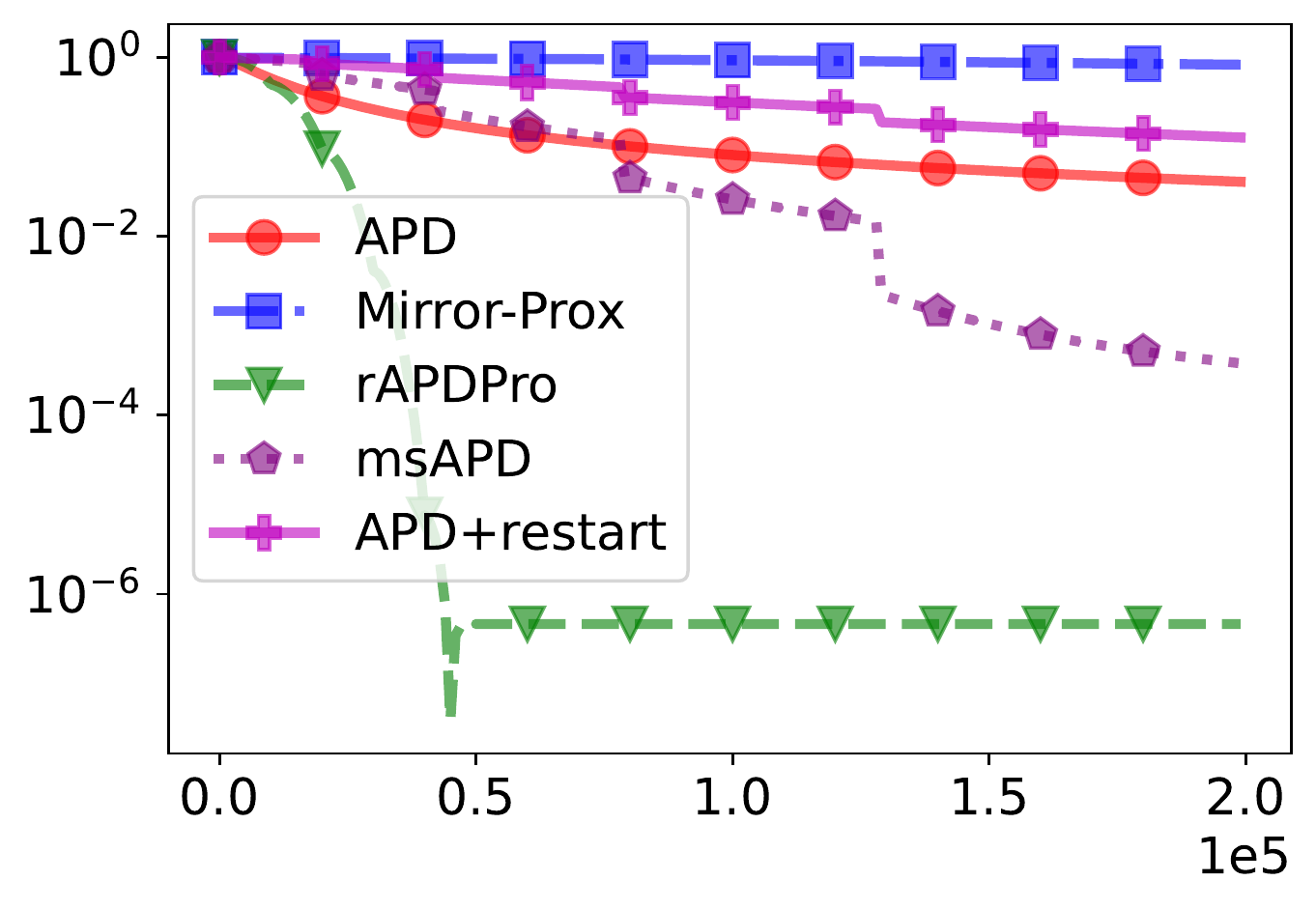}%
\end{minipage}
\\
\begin{minipage}[t]{0.33\columnwidth}%
      \includegraphics[width=4.5cm]{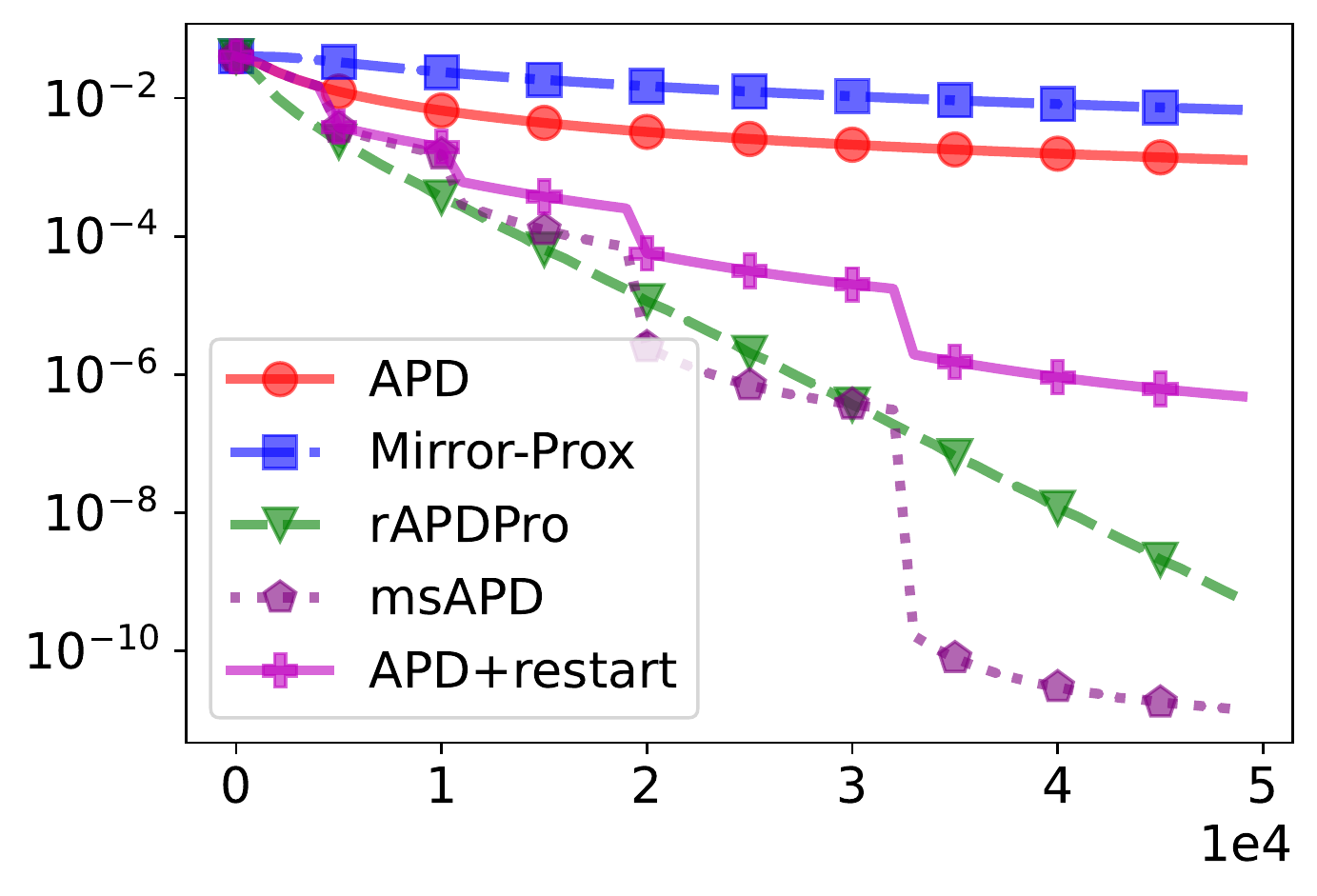}%
\end{minipage}%
\begin{minipage}[t]{0.33\columnwidth}%
\includegraphics[width=4.5cm]{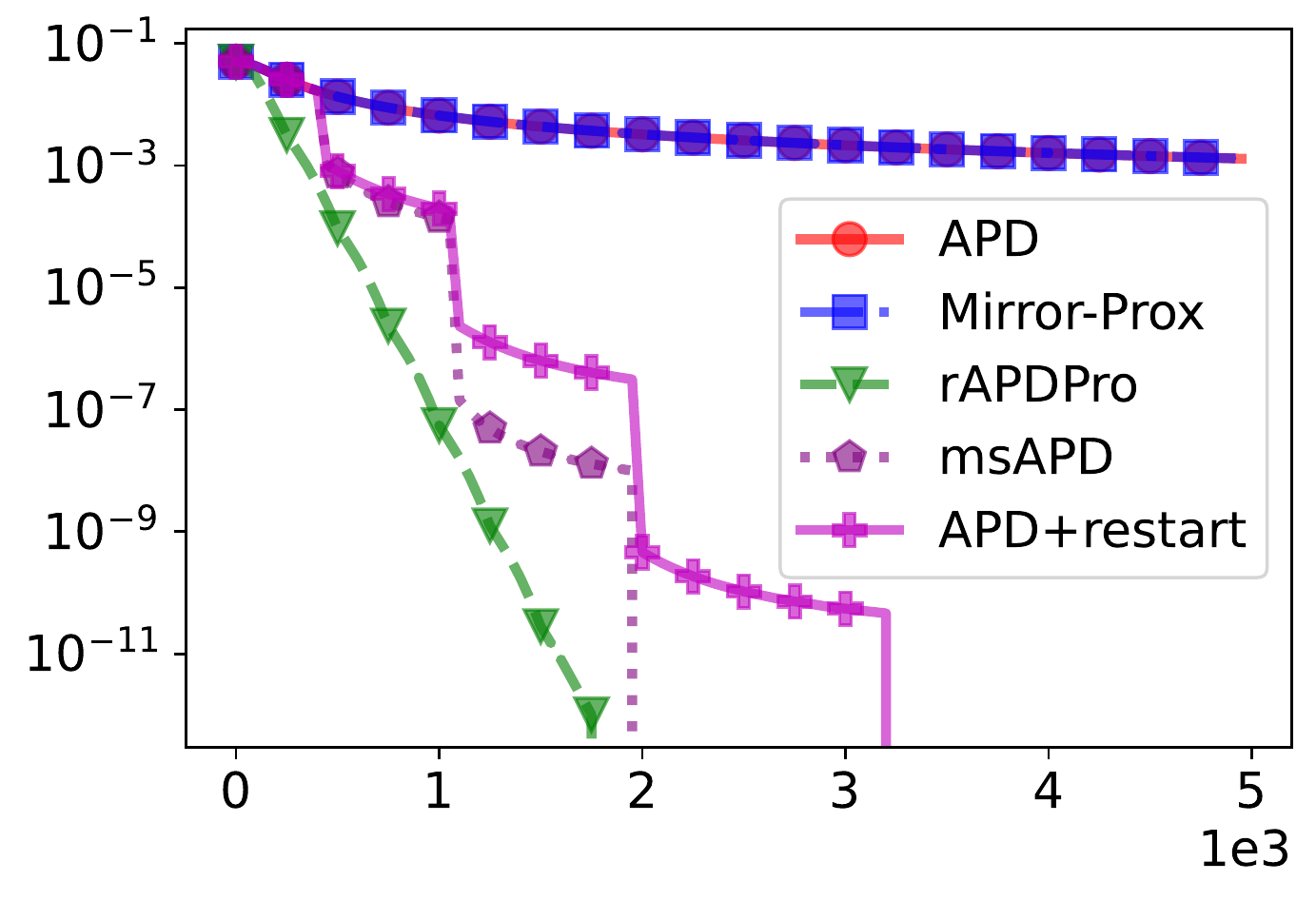}%
\end{minipage}%
\begin{minipage}[t]{0.33\columnwidth}%
\includegraphics[width=4.5cm]{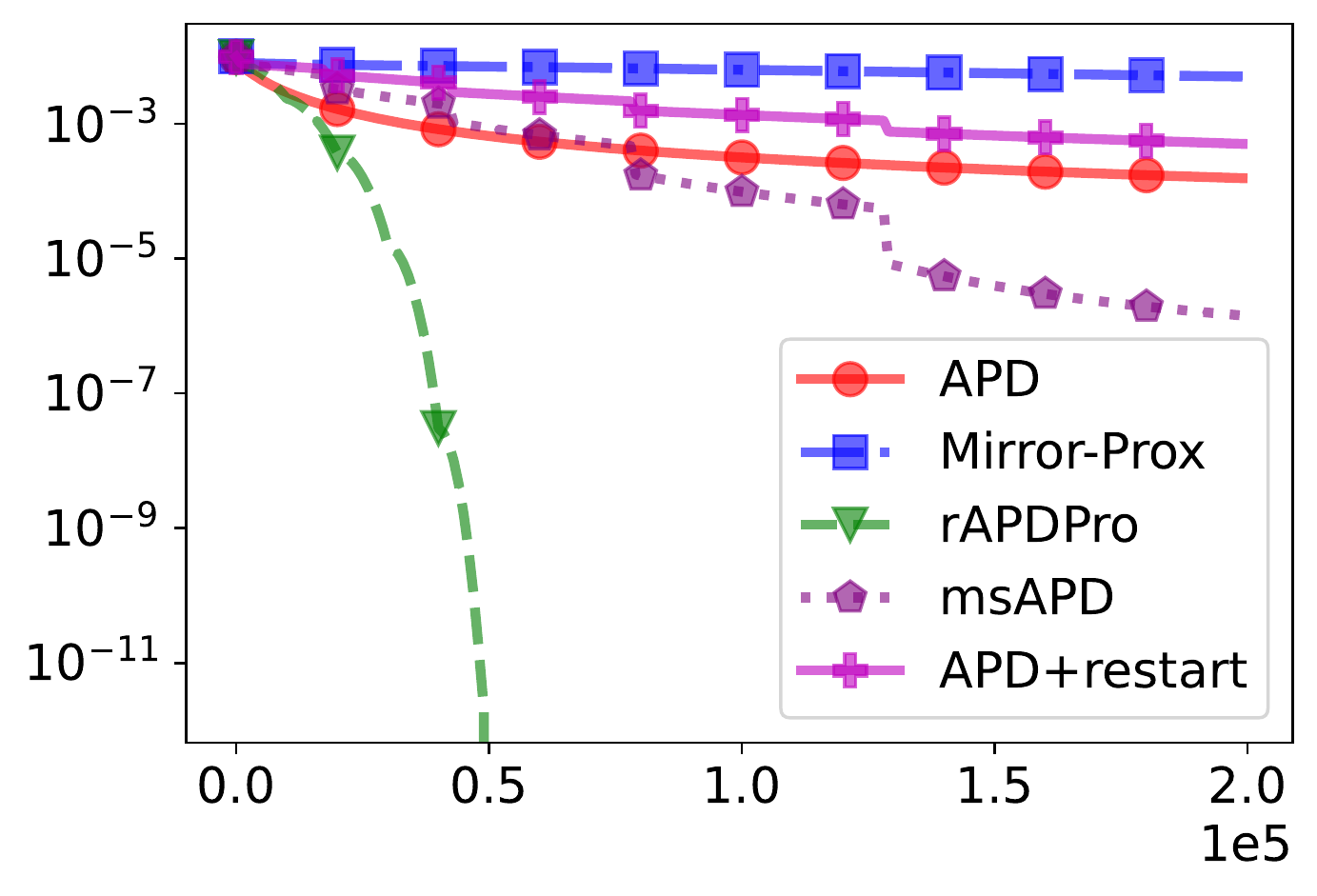}%
\end{minipage}\par
\end{centering}
\caption{\label{fig:Objective-Function-Value} The first row describes the convergence to optimum, where the $y$-axis reports $\log_{10}((\|D^{1/2}\xbf_{k}\|_{1}-\|D^{1/2}\xbf^{*}\|_{1})/\|D^{1/2}\xbf^{*}\|_{1})$
for {\TPAPD}, and  $\log_{10}((\|D^{1/2}\bar{\xbf}_{k}\|_{1}-\|D^{1/2}\xbf^{*}\|_{1})/\| D^{1/2}\xbf^{*}\|_{1})$
for APD, APD+restart, {\MSAPD} and {\mirror} ($\protect\xbf^{*}$ is computed by MOSEK~\citep{aps2019mosek}). The second row describes feasibility violation, where $y$-axis reports the feasibility gap $\log_{10}(\max\{0,G(\protect\xbf_k)\})$
for {\TPAPD}, and $\log_{10}(\max\{0,G(\bar{\protect\xbf}_k)\})$ for
APD, {\MSAPD} and {\mirror}.  Datasets (Left-Right order) correspond to bio-CE-HT, bio-CE-LC and econ-beaflw.}
\end{figure}

\section{\label{sec:numerical}Numerical study}
In this section, we examine the empirical performance of our proposed algorithms for solving
the sparse Personalized PageRank~\citep{fountoulakis2019variational, fountoulakis2022open,martinez2023accelerated}.
% Let $G=(V,E)$ be a connected undirected graph with $n$ vertices. Denote the adjacency matrix of $G$ by $A$, that is, $A_{i,j}=1$ if $i\sim j$ and $0$ otherwise. Let $D=\diag(d_1,\ldots,d_n)$ be the matrix with the degrees $\{d_i\}_{i=1}^n$ in its diagonal.
The constrained form of Personalized PageRank can be written as follows:
$\min_{\xbf\in \mbb R^n}\ \ \norm{D^{1/2}\xbf}_1\ \  \st\,\, \tfrac{1}{2}\inprod{\xbf}{Q\xbf}-\alpha \inner{\mathbf{s}}{D^{-1/2}\xbf}\leq b,$
% \begin{equation}
%     \begin{aligned}
%         \min_{\xbf\in \mbb R^n}\,\, &\norm{D^{1/2}\xbf}_1\\
%         \st\,\, & \tfrac{1}{2}\inprod{\xbf}{Q\xbf}-\alpha \inprod{\mathbf{s}}{D^{-1/2}\xbf}\leq b,
%     \end{aligned}
% \end{equation}
% where $Q=D^{-1/2}\brbra{D-\tfrac{1-\alpha}{2}(D+A)}D^{-1/2}$, $\alpha \in (0,1)$, $\mathbf{s}\in \Delta^n$ is a teleportation distribution over the nodes of the graph $G$ and $b$ is a pre-specific target level.
where $Q,D$ and $\mathbf{s}$ are generated by graph.
We implement both  \resapdsce{} and {\MSAPD}. 
We skip \APDSCE{} as we observe that the restart strategy consistently improves the algorithm performance. 
For comparison, we consider the state-of-the-art accelerated primal-dual (APD) method~\citep{hamedani2021primal},  APD with restart mechanism at fixed iterations (APD+restart) and   {\mirror}~\citep{he2015mirror}. 
6 small to medium-scale datasets from various domains in the Network Datasets~\citep{networkdata} are selected in our experiments. All experiments are implemented on  Mac mini M2 Pro, 32GB.  Due to the page limit, we only report results on three datasets and leave more details in the last Appendix~\ref{sec:appendix_exp_details}.

We plot the  relative  function value gap $|f(\xbf)-f(\xbf^{*})|/|f(\xbf^{*})|$ and the feasibility violation $\max\{G(\xbf), 0\}$ over the iteration number in Figure~\ref{fig:Objective-Function-Value}, respectively.
% Due to space limitations, additional experimental results are included in the appendix.
Firstly, in terms of both optimality gap and constraint violation, the performance of {\resapdsce} and {\MSAPD} is significantly better than that of {\APD}, APD+restart and {\mirror}. Additionally, {\resapdsce} and {\MSAPD} often converge to high-precision solutions. Secondly, based on the experimental results, it is indeed observed that {\MSAPD} exhibits a periodic variation in convergence performance, which aligns with our algorithm theory.

% First, we observe that  \APD{} is able to generate feasible solutions in all the datasets. However, it appears to have inferior performance in reducing the objective value. In contrast, \mirror{} often exhibits a large infeasibility error.
% Second, it is interesting to compare the performance of \MSAPD{} and \APD{}. Notice that both \MSAPD{} and \APD{} employ the same stepsize rules while \MSAPD{} can further leverage the strong convexity to improve the convergence rate.  Our  results  confirm the empirical advantage of the multi-stage strategy.
% Third, we notice that {\TPAPD} outperforms  the other algorithms in all the datasets.
% The experimental results indeed confirm the theoretical advantage of the restart strategy in \TPAPD{}.

\begin{figure}[h]
\begin{centering}
\begin{minipage}[t]{0.33\columnwidth}%
\includegraphics[width=4.5cm]{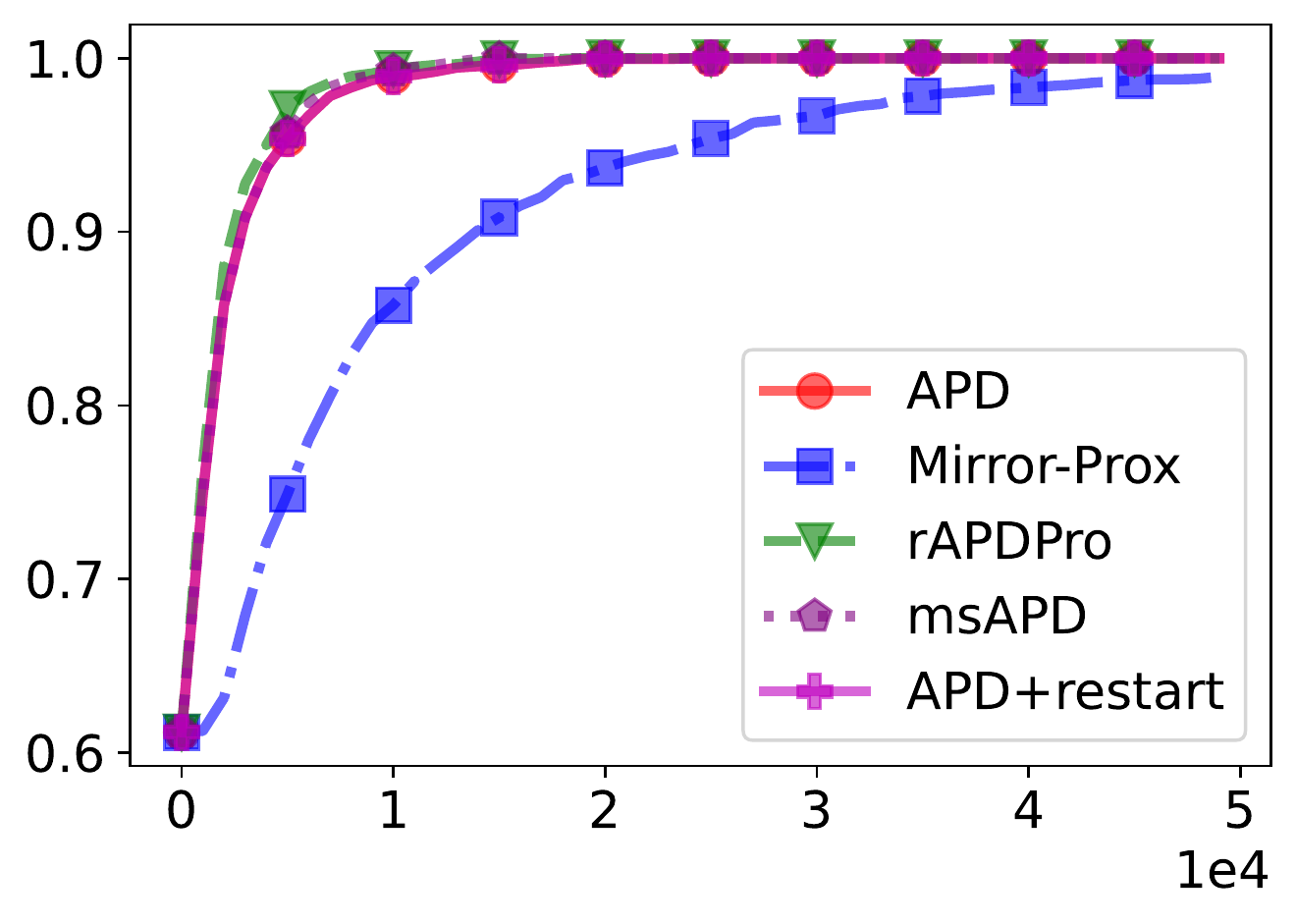}%
\end{minipage}%
\begin{minipage}[t]{0.33\columnwidth}%
\includegraphics[width=4.5cm]{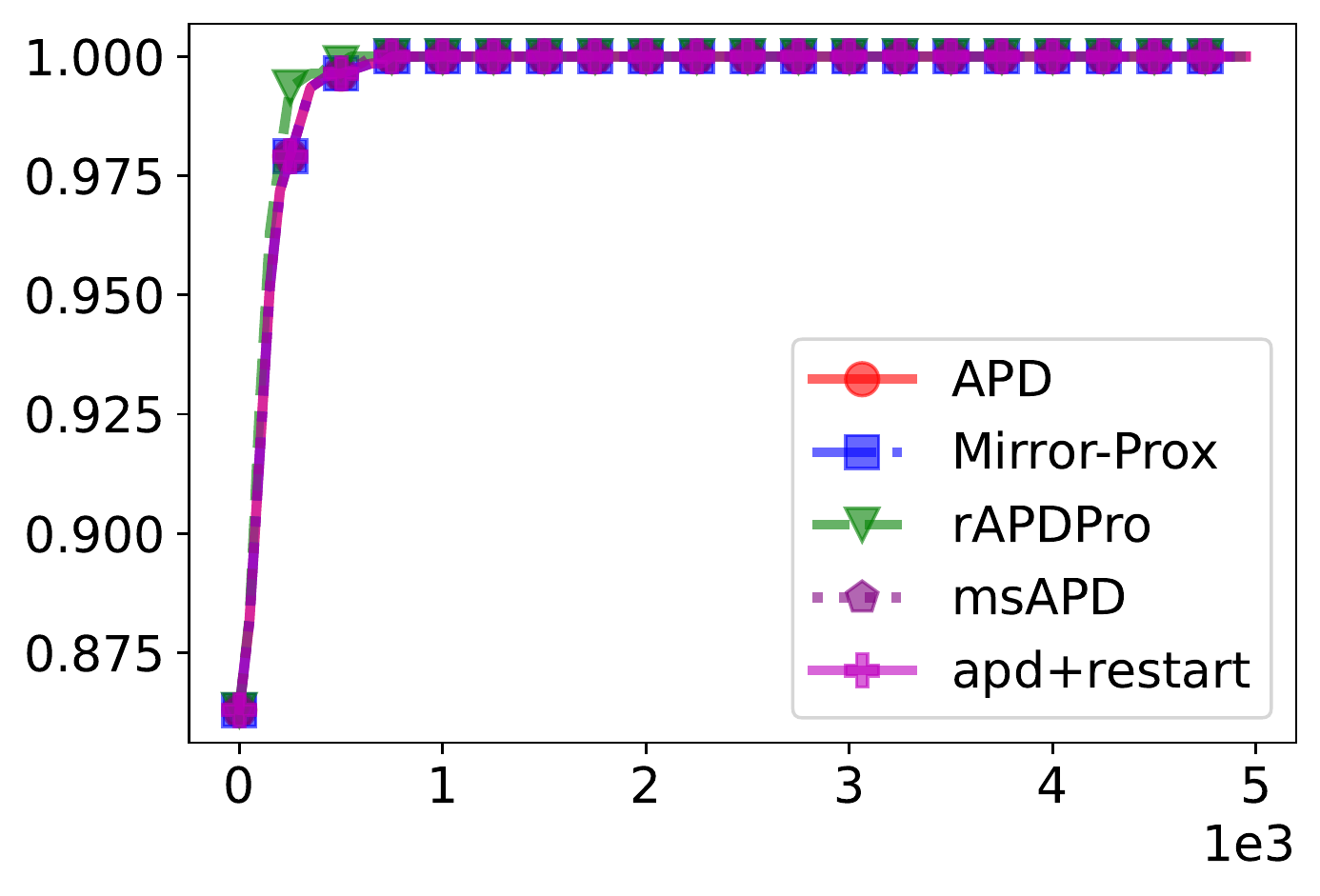}%
\end{minipage}%
\begin{minipage}[t]{0.33\columnwidth}%
\includegraphics[width=4.5cm]{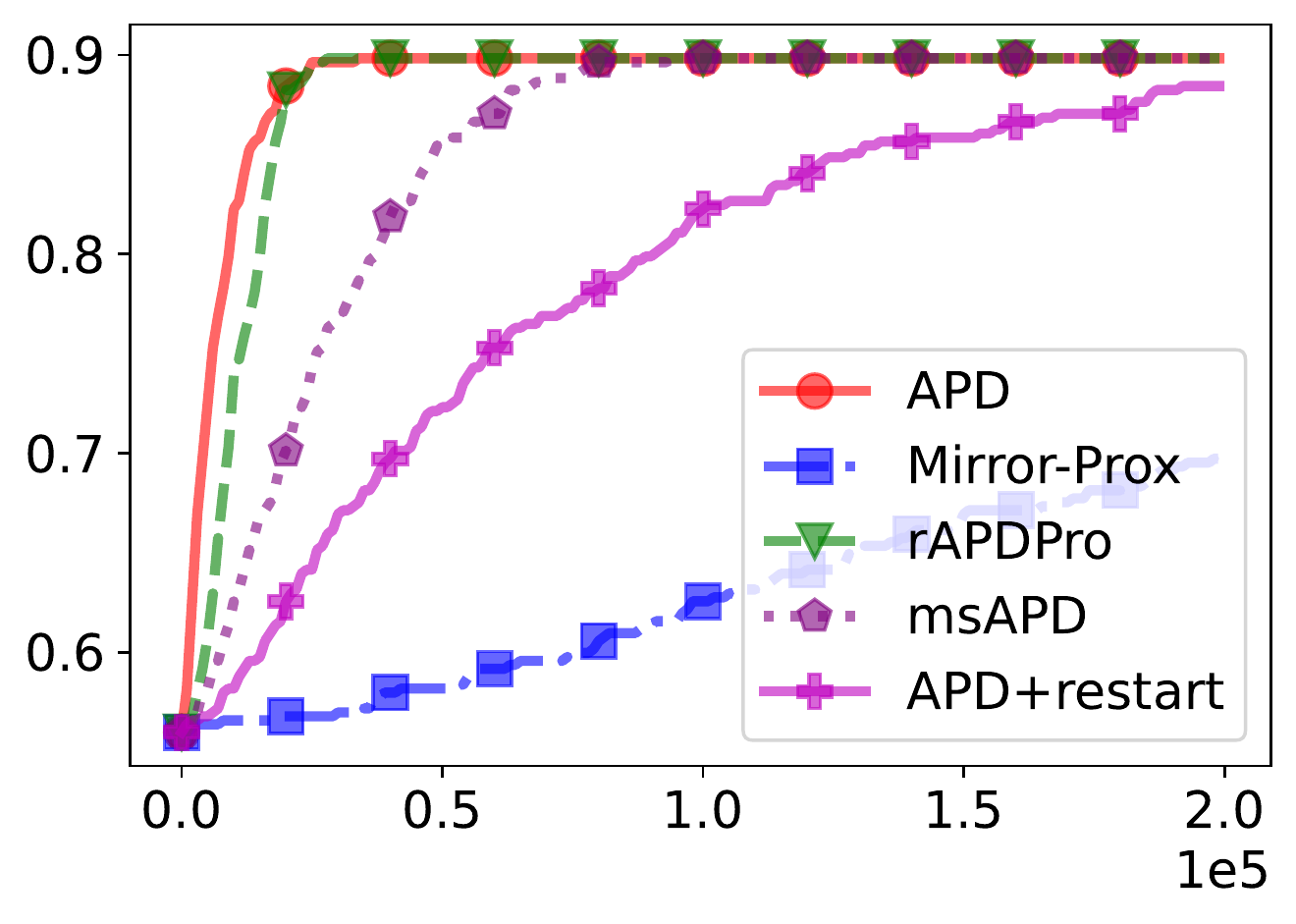}%
\end{minipage}
\par\end{centering}
\caption{\label{fig:identification}The experimental results on active-set identification. Datasets (Left-Right order) correspond to bio-CE-HT, bio-CE-LC and econ-beaflw.
The $x$-axis reports the iteration number and  the $y$-axis reports accuracy in active-set identification.}
\end{figure}

% \paragraph{Sparsity identification}
Next, we examine the algorithm's effectiveness in identifying sparsity patterns. We computed a nearly optimal solution $\xbf^*$ from MOSEK. Note that $\xbf^*$ is a dense vector.
For numerical consideration, we truncate the coordinate values of $\xbf^*$ to zero if the absolute value is below $10^{-8}$ and perform the same truncation to all the generated solutions of the compared algorithms.  
Then we use  $(|\Acal(\xbf)\cap\Acal(\xbf ^*)| + |\Acal^c(\xbf)\cap \Acal^c(\xbf^*)|)/n$ to measure the accuracy of identifying the active set, where $|\cdot|$  denotes the set cardinality. 
For {\TPAPD}, we consider the last iterate $\xbf_k$ while for \APD{}, {\MSAPD} and {\mirror}, we plot the result on $\bar{\xbf}_k$, as these are the solutions where the convergence rates are established.
Figure~\ref{fig:identification} plots the experiment result, from which we observe that 
{\TPAPD} and {\MSAPD} are highly effective in identifying the active set. Often, they are able to recognize the structure of the active set within a small number of iterations.  
% While {\APD} can identify the active set in many cases,  it fails in peking-1.
Overall, the experimental results show the great potential of our proposed algorithms in identifying the sparsity structure and are consistent with our theoretical analysis.

\section{\label{sec:Conclusion}Conclusion}
The key contribution of this paper is that we develop several new first-order primal-dual algorithms for convex optimization with strongly convex constraints. Using some novel strategies to exploit the strong convexity of the Lagrangian function, we substantially improve the best convergence rate from $\mathcal{O}(1/\vep)$ to $\mathcal{O}(1/\sqrt{\vep})$.  
In the application of constrained sparse learning problems, the experimental study confirms the advantage of our proposed algorithms against state-of-the-art first-order methods for constrained optimization. Moreover, we show that one of our proposed algorithms {\resapdsce} has the favorable feature of identifying the sparsity pattern in the optimal solution.
For future work, one direction is to apply the adaptive strategy, such as line search, to our framework to deal with cases when the dual bound is unavailable. Another interesting direction is to further exploit the active set identification property in a general setting. For example, it would be interesting to incorporate our algorithm with active constraint identification, which could be highly desirable when there are a large number of constraints. It would also be interesting to consider a more general convex objective when the proximal operator is not easy to compute.  

\section*{Acknowledgement}
This research is partially supported by the Major Program of National Natural Science Foundation of China (Grant 72394360, 72394364). We sincerely thank all the reviewers for their valuable suggestions, which have significantly improved the quality of our article.

%% file: appendix_neurips.tex
\section*{Structure of the Appendix}
The appendix is structured as follows: Appendix~\ref{sec:limitation} introduces some limitations of our methods, primarily concerning the application scenarios of our algorithm. Appendix~\ref{sec:frank} includes comparisons between ours and some related Frank-Wolfe methods. We give some auxiliary lemmas in Appendix~\ref{sec:aux-Lemmas}, which are very important for the proofs presented later. Appendix~\ref{sec:appendix:preproof},~\ref{sec:Proof-of-Tpapd},~\ref{sec:Proof Details in restart} and~\ref{sec:Proof-of-manifold} present the proof of conclusion in Section~\ref{sec:Preliminaries},~\ref{sec:APDSCE},~\ref{sec:restarted} and~\ref{sec:Manifold}, respectively. Furthermore, Appendix~\ref{sec:MSAPD} introduces a new algorithm to obtain a convergence rate without complicated dual updating. Finally, Appendix~\ref{sec:appendix_exp_details} offers more extensive details on our experiments.

\section{Limitations\label{sec:limitation}}
In this paper, we focus on the theoretical analysis of convex optimization. Although our proposed algorithms for the convex optimization with strongly convex constraints can theoretically improve the existing results from $\mcal O(1/\vep)$ to $\mcal O(1/\sqrt{\vep})$. However, we still need to point out that our optimization algorithm has the following limitations. One is the algorithm needs a lower bound on the norm of sub-gradients of the objective function in the optimal solution, which may not be satisfied for all functions. On the other hand, we require consistent smoothness of the constraints to ensure convergence, and how to use the line search method to ensure convergence is a future direction.

\section{Comparison with  Frank-Wolfe~\label{sec:frank}}
We note that the strongly convex function constraint in \eqref{eq:constrainedproblem} is a special case of a strongly convex set constraint, as demonstrated in~\cite{journee2010generalized}.
Over the strongly convex set, it has been shown that Frank-Wolfe Algorithm (FW) can obtain convergence rates substantially better than the worst-case $\Ocal(1/\vep)$ rate. Under the bounded gradient assumption, \cite{dunn1979rates, levitin1966constrained} show that FW obtains linear convergence over a strongly convex set.
 Nevertheless, the uniform bounded gradient assumption appears to be stronger than ours, as we only impose the lower boundedness assumption on the optimal solution $\xbf^*$ and allow the objective to be non-differentiable. 
 More recently, \cite{garber2015faster} shows that FW obtains an $\Ocal(1/\sqrt{\vep})$ rate when the gradient is the order of the square root of the function value gap.  For more recent progress, please refer to~\cite{braun2022conditional}. Despite the attractive convergence property,  FW exhibits certain limitations when applied to the general function constraints~\eqref{eq:constrainedproblem} addressed in this paper.  Specifically, FW involves a sequence of linear optimization problems throughout the iterations. While linear optimization over certain strongly convex sets, such as $\ell_p$-ball, admits a closed-form solution, there exists no efficient routine to handle general function constraints explored in this paper.

\section{\label{sec:aux-Lemmas}Auxiliary lemmas}
The following three-point property is important in the convergence analysis. 
\begin{lemma} 
\label{lem:descent}Let $f:\mathbb{R}^{n}\to\mathbb{R\cup}\{+\infty\}$
be a closed strongly convex function with modulus $\mu\geq0$.
Give $\bar{\xbf}\in\mathcal{X}$,  where $\mathcal{X}$ is a compact
convex set and $t\geq0$, let
$\xbf^{+}=\argmin_{x\in\mathcal{X}}f(\xbf)+\frac{t}{2}\|\xbf-\bar{\xbf}\|^{2},$
then for all $\xbf \in \mcal{X}$, we have
\[
f(\xbf)+\frac{t}{2}\|\xbf-\bar{\xbf}\|^{2}\geq f(\xbf^{+})+\frac{t+\mu}{2}\|\xbf^{+}-\xbf\|^{2}+\frac{t}{2}\|\xbf^{+}-\bar{\xbf}\|^{2}.
\]
\end{lemma}
\proof
   Since $\Xcal$ is a convex compact set,  $\phi(x)\coloneqq I_{\Xcal}(\xbf)+f(\xbf)+\frac{t}{2}\norm{\xbf - \bar{\xbf}}^2$ is lower-semi-continuous and $(\mu + t)$-strongly convex, where $I_{\Xcal}(\xbf)=\begin{cases}
       0 & \xbf \in \Xcal \\
        \infty & \xbf \notin \Xcal 
   \end{cases}$. 
Using the optimality ($\zerobf\in \phi(\xbf^+)$) and strong convexity, we have $\phi(\xbf)\ge \phi(\xbf^+) + \inner{\zerobf}{\xbf-\xbf^+} + \frac{t+\mu}{2}\|\xbf^{+}-\xbf\|^{2}$, for any $\xbf\in\Xcal$. This immediately gives the desired relation.
\endproof

% \begin{lemma}[Theorem 1 in~\cite{robbins1971convergence}]
% \label{lem:limit_converge}
%     Let $\left\{a_j\right\},\left\{b_j\right\}$ and $\left\{c_j\right\}$ be non-negative real sequences such that $a_{j+1}\le a_j - b_j + c_j$ for all $j\ge 0$, and $\sum_{j=0}^{\infty} c_j <\infty$. Then $a=\lim_{j\to \infty} a_j$ exists, and $\sum_{j=0}^{\infty} b_j < \infty$. 
% \end{lemma}
% \proof
%     Complete proof can be found in~\cite{robbins1971convergence}.
% \endproof
The following result is adjusted from the classic supermartingale convergence theorem~\citep[Theorem 1]{robbins1971convergence}.
We give proof for completeness.
\begin{lemma}\label{lem:limit_converge}
    Let $(\Omega,\mcal F,\mbb P)$ be a probability space and $\mcal F_1\subset \mcal F_2 \subset \cdots $ be a sequence of sub-$\sigma$-algebras of $\mcal F$. For each $j=1,2,\cdots$, let $a_j,b_j$ and $c_j$ be non-negative $\mcal F_n$-measure random variables such $\mbb E[a_{j+1}\mid \mcal F_j]\leq a_j - b_j + c_j$, then we have $\lim_{j\to \infty}a_j<\infty$ exists and $\sum_{j=1}^{\infty}b_j <\infty$ a.s. when $\sum_{j=1}^{\infty}c_j<\infty$.
\end{lemma}
\proof
    Define $d_{j}=a_{j}-\sum_{l=1}^{j-1}(c_{l}-b_{l})$ and for any $\bar{a}>0$, define
$t=\inf\{t:\sum_{l=1}^{t}c_{l}>\bar{a}\}$.
If $j<t$, we have 
\begin{equation}\label{eq:decreasing}
    \mbb E[d_{j+1}\mid\mcal F_{j}]=\mbb E[a_{j+1}-\sum_{l=1}^{j}(c_{l}-b_{l})\mid\mcal F_{j}]\aleq a_{j}-\sum_{l=1}^{j-1}(c_{l}-b_{l})=:d_{j},
\end{equation}
where $(a)$ holds by $\mbb E[a_{j+1}\mid\mcal F_{j}]\leq a_{j}+c_{j}-b_{j}$,
and hence 
\[
\mbb E[d_{\min\{t,(j+1)\}}\mid\mcal F_{j}]=d_{t}\mbb I_{\{t\leq j\}}+\mbb E[d_{j+1}\mid\mcal F_{j}]\mbb I_{\{t>j\}}\aleq d_{\min\{t,j\}},
\]
where $(a)$ holds by~\eqref{eq:decreasing}.
Therefore, we have $\{d_{\min\{t,(j+1)\}},\mcal F_{j},1\leq j\leq\infty\}$
is a supermartingale. Since 
\begin{equation*}
d_{\min\{t,j\}} = a_{\min\{t,j\}} - \sum_{l=1}^{\min\{t,j\} - 1}(c_l - b_l)\ageq-\sum_{l=1}^{\min\{t,(j-1)\}}c_{l}\geq-\bar{a},
\end{equation*}
holds for all $j$, where $(a)$ holds by $a_{\min\{t,j\}},b_l \geq 0$. Then it follows from the martingale convergence
theorem that $\lim_{j\to\infty}d_{\min\{t,j\}}$ exists and is finite
a.s., i.e., $\lim_{j\to\infty}d_{j}$ exists and is finite on $\{t=\infty\}=\{\sum_{j=1}^{\infty}c_{j}\leq\bar{a}\}$.
Since $\bar{a}$ is arbitrary, we see that $\lim_{j\to\infty}d_{j}$
exists and is finite a.s. on $\{\sum_{j=1}^{\infty}c_{j}<\infty\}$.
By $d_{j}=a_{j}-\sum_{l=1}^{j-1}(c_{l}-b_{l})$, we have $\lim_{j\to\infty}a_{j}$
exists and is finite and $\sum_{j=1}^{\infty}b_{j}<\infty$ when $\{\sum_{j=1}^{\infty}c_{j}<\infty\}$.
\endproof

\section{Proof details in Section~\ref{sec:Preliminaries}\label{sec:appendix:preproof}}
% \begin{proposition}
% \label{prop:bounded_set_y}Suppose Assumption~\ref{assu:Slater's}
% holds. Then, for any optimal solution $\xbf^*$ of problem~\eqref{eq:constrainedproblem}, there exists $\ybf^*\in\Rbb^m$ such that the KKT condition holds. Moreover,  $\ybf^*$ is
% bounded within the set $\mathcal{Y}:=\big\{\ybf\in\mathbb{R}_{+}^{m}\mid\|\ybf\|\leq  \bar{c}\big\}$, and 
% % any optimal dual variable
% satisfies $\norm{\ybf^*}_1\leq \bar{c}$, where $\bar{c}\coloneqq \frac{f(\tilde{\xbf})-\min_{\xbf\in \mbb R^n}f(\xbf)}{\min_{i\in[m]}\{-g_{i}(\tilde{\xbf})\}}$.
% \end{proposition}
% \begin{proof}
\subsection{Proof of Proposition~\ref{prop:bounded_set_y}\label{sec:proof_of_prop_bounded_set_y}}
\proof 
% Under Slater's CQ, the existence of the KKT condition  is a standard result in nonlinear programming. 
Under Slater's CQ, it is standard to show that any optimal solution $\xbf^*$ will also satisfy the KKT condition.
For example, one can refer to~\cite{bertsekasnonlinear}.
For any  $\xbf\in\Xcal_G$,  we have 
\[
f(\xbf)+\left\langle \ybf^{*},G(\xbf)\right\rangle \geq f(\xbf^{*})+\left\langle \ybf^{*},G(\xbf^{*})\right\rangle =f({\xbf}^{*}),
\]
where the equality is from the complementary slackness. In view of the above result and the Slater's condition (i.e., $G(\tilde{\xbf})<\boldsymbol{0}$),
we have
\begin{equation}
    f(\tilde{\xbf})> f(\tilde{\xbf})+\left\langle \ybf^{*},G(\tilde{\xbf})\right\rangle \geq f(\xbf^{*}).
\end{equation}
Combining with fact $\norm{\ybf^*}_1\min_{i\in[m]}\bcbra{-g_i(\tilde{\xbf})}\leq -\inner{\ybf^*}{G(\tilde{\xbf})}$, then we have
\begin{equation}
\begin{split} \|\ybf^{*}\|\leq\|\ybf^{*}\|_{1}\leq\frac{f(\tilde{\xbf})-f(\xbf^{*})}{\min_{i\in[m]}\{-g_{i}(\tilde{\xbf})\}} = \bar{c},
\end{split}
\label{eq:tighter}
\end{equation}
where the last inequality is by $f(\xbf^{*})\ge\min_{\xbf\in\mathbb{R}^{n}}f(\xbf)$.
% , which is derived by $\xbf_{0}^{*}\in\argmin_{\xbf\in\mathbb{R}^{n}}f(\xbf)$.
\endproof

% \begin{proposition}
% \label{prop:uniquex}
% Under Assumption~\ref{assu:dual-optim}, problem~\eqref{eq:constrainedproblem} has a unique optimal solution $\xbf^{*}$. Let $\mathcal{Y}^{*}:=\argmax_{\ybf\in \mbb R_+^m}\mcal L(\xbf^{*},\ybf)$
% denote the set containing all the optimal dual variables, then $\mathcal{Y}^{*}$
% is a convex set. 
% \end{proposition}
\subsection{Proof of Proposition~\ref{prop:uniquex}\label{sec:proof_of_proposition_uniquex}}
\proof 
We prove the uniqueness property by contradiction. Suppose that there exist
$(\xbf^{*},\ybf^{*})$, $(\tilde{\xbf}^{*},\tilde{\ybf}^{*})$ satisfying
the KKT condition, then from the complementary slackness, optimality of $\xbf^*$ and $\tilde{\xbf}^*$, we have 
\[\mathcal{L}(\xbf^{*},\ybf^{*})=f(\xbf^{*})=f(\tilde{\xbf}^{*})=\mathcal{L}(\tilde{\xbf}^{*},\tilde{\ybf}^{*}).\]
Moreover,  we have
$\mathcal{L}(\tilde{\xbf}^{*},\tilde{\ybf}^{*})\leq\mathcal{L}(\xbf^{*},\tilde{\ybf}^{*})\leq\mathcal{L}(\xbf^{*},\ybf^{*}).$
Hence, we must have $\mathcal{L}(\tilde{\xbf}^{*},\tilde{\ybf}^{*})=\mathcal{L}(\xbf^{*},\tilde{\ybf}^{*})$.
However, since Assumption~\ref{assu:dual-optim} implies $\tilde{\ybf}^{*}\neq\boldsymbol{0}$, the strongly convex function $\mathcal{L}(\cdot,\tilde{\ybf}^{*})$
has a unique optimizer. Therefore, we conclude that $\xbf^{*}=\tilde{\xbf}^{*}$.

Next, we show that the set of optimal dual variables for problem~\eqref{eq:constrainedproblem}
is convex. Suppose that there exist two optimal dual variables $\ybf_{1}^{*}$
and $\ybf_{2}^{*}$ for the unique primal variable $\xbf^{*}$, both satisfying the KKT condition, then we have
$\left\langle \ybf_{1}^{*},G(\xbf^{*})\right\rangle =\left\langle \ybf_{2}^{*},G(\xbf^{*})\right\rangle =0$.
This implies that any linear combination of $\ybf_{1}^{*}$ and $\ybf_{2}^{*}$
satisfy KKT condition, i.e., $\left\langle a\ybf_{1}^{*}+b\ybf_{2}^{*},G(\xbf^{*})\right\rangle =0,\forall a,b$.
From Proposition~\ref{prop:bounded_set_y}, we know any optimal dual
variable falls into a bounded convex set $\mathcal{Y}$. The intersection of two convex sets is also a convex set. Hence, we complete our proof.
\endproof

% \begin{proposition}
% \label{prop:bounded_set_x}Under Assumptions~\ref{assu:Slater's}
% and~\ref{assu:dual-optim},  we have
%  $\max_{\xbf_1, \xbf_2 \in \mcal{X}_{G} } \|\xbf_1-\xbf_2\| \leq \min_{i\in[m]}2\sqrt{\frac{-2g_{i}(\xbf_{i}^{*})}{\mu_{i}}}$,
% where $\xbf_{i}^{*}=\argmin_{\xbf\in\Rbb^n}g_{i}(\xbf)$. Let $\mcal X:=\mcal B\brbra{\tilde{\xbf},\min_{i\in[m]}2\sqrt{{-2g_{i}(\xbf_{i}^{*})}/{\mu_{i}}}+\zeta}$, 
% where $\zeta$ is a positive constant, then $\xbf^{*}\in\intr\mathcal{X}$.
% \end{proposition}
% \begin{proof}
\subsection{Proof of Proposition~\ref{prop:bounded_set_x}\label{sec:proof_of_proposition_bounded_set_x}}
\proof 
From the strong convexity of $g_{i}(\xbf)$, we have $g_{i}(\xbf)\geq g_{i}(\xbf_{i}^{*})+\frac{\mu_{i}}{2}\|\xbf-\xbf_{i}^{*}\|^{2},$ which implies
\begin{equation}
    \begin{aligned}
        \|\tilde{\xbf}-\xbf_{i}^{*}\|^{2}\leq(g_{i}(\tilde{\xbf})-g_{i}(\xbf_{i}^{*}))\frac{2}{\mu_{i}}\overset{(a)}{<} \frac{-2g_{i}(\xbf_{i}^{*})}{\mu_{i}},\\
        \|\xbf^*-\xbf_{i}^{*}\|^{2}\leq(g_{i}(\xbf^*)-g_{i}(\xbf_{i}^{*}))\frac{2}{\mu_{i}} \leq \frac{-2g_{i}(\xbf_{i}^{*})}{\mu_{i}},
    \end{aligned}
\end{equation}
where $(a)$ holds by $g_i(\tilde{\xbf})<0$.
In view of the triangle inequality and the above result, we have
\[
\|\tilde{\xbf}-\xbf^*\|  \leq\|\xbf_i^*-\xbf^*\|+\|\tilde{\xbf}-\xbf_i^{*}\| < 2\sqrt{\frac{-2g_{i}(\xbf_{i}^{*})}{\mu_{i}}}.
\]
Hence, $\xbf^*\in \textup{int }\mathcal{B}\Brbra{\tilde{\xbf},\min_{i\in[m]}2\sqrt{\frac{-2g_{i}(\xbf_{i}^{*})}{\mu_{i}}}}$.
\endproof

\section{\label{sec:Proof-of-Tpapd}Convergence analysis of 
\APDSCE{}
% Proof details in Section~\ref{sec:APDSCE}
}

% The main proof of Theorem~\ref{thm:gapdiminishing} is based on the modification of~\cite{hamedani2021primal}.

\subsection{Proof of Proposition~\ref{prop:local_estimation}\label{sec:proof_of_propostion_local_estimation}}

\proof
Using
the triangle inequality and~\eqref{eq:lipschitz_x}, we have
\[
\|\nabla G(\xbf^{*})\|-\|\nabla G(\hat\xbf)\|\le\|\nabla G(\xbf^{*})-\nabla G(\hat\xbf)\|\leq L_{X}\|\hat\xbf-\xbf^{*}\|.
\]
Combining the above inequality and~\eqref{eq:r_lower}, we obtain
\begin{equation}\label{eq:mid-04}
\frac{r}{\norm{\ybf^*}_1}\le L_X \|\hat\xbf-\xbf^{*}\|+\|\nabla G(\hat\xbf)\|.
\end{equation}
Next,  we develop more specific lower bounds on $\norm{\ybf}_1$.  i). Inequality~\eqref{eq:y-lb-1} can be easily verified since  we have $\|\hat\xbf-\xbf^*\| \le \sqrt{2\beta}$.  
ii). Suppose ${(\ybf^*)^{\top}\mubm}\cdot \|\hat\xbf-\xbf^*\|^2 \le 2\beta$, then together with \eqref{eq:mid-04}  we have 
$$
\frac{r}{\norm{\ybf^{*}}_1}\le L_X \sqrt{\frac{2\beta}{(\ybf^*)^{\top}\mubm}}+\norm{\nabla G(\hat\xbf)}\le L_X \sqrt{\frac{2\beta}{\underline{\mu}\norm{\ybf^*}_1}}+\norm{\nabla G(\hat\xbf)}.
$$
Note that the above inequality can be expressed as $at^2-bt-c \le 0$ with $t=\norm{\ybf^*}_1^{-1/2}$, $a=r, b=L_X\sqrt{{2\beta}/{\underline{\mu}}}$ and $c=\norm{\nabla G(\hat\xbf)}$. Standard analysis implies that $t\le (b+\sqrt{b^2+4ac})/{2a}$, which gives the desired bound~\eqref{eq:y-lb-2}.
\endproof

\subsection{\label{subsec:Proof_of_thmgap}Proof of Theorem~\ref{thm:gapdiminishing}\label{sec:proof_of_thm_gapdiminishing}}
\proof
First,  it is easy to verify by our construction that $\{\Ycal_k\}$ is a monotone sequence: $\Ycal_1\supseteq\Ycal_2\supseteq\ldots\supseteq \Ycal_k\ldots$.  Our goal is to show $\Ycal^{*}\subseteq\mcal Y_{k}$ holds
for any $k\geq0$ by induction. Note that $\Ycal^{*}\subseteq\mcal Y_{0}$ immediately follows from our assumption that $(\ybf^{*})^{\top}\boldsymbol{\mu}\geq\rho_{0}$, for any $\ybf^{*}\in\mcal Y^{*}$. Suppose that $\Ycal^{*}\subseteq\mcal Y_{k}$
holds for $k=0,\ldots,K-1$, we claim:
\begin{enumerate}
\item For any $\xbf \in \mcal{X}$ and $\ybf \in \mcal{Y}^*$, we have
\begin{equation}
\mathcal{L}(\bar{\xbf}_{K},\ybf)-\mathcal{L}(\xbf,\bar{\ybf}_{K})\leq\frac{1}{T_{K}}\Delta(\xbf,\ybf)-\frac{t_{K-1}\tau_{K-1}^{-1}}{2T_{K}}\|\xbf-\xbf_{K}\|^{2}. \label{eq:thm2}
\end{equation}
\item  $\Ycal^*\subseteq \Ycal_K$.
\end{enumerate}

Part 1. For $k=0,1,2,\ldots, K-1$, taking $-\left\langle \zbf_{k},\cdot\right\rangle $
and $f(\cdot)+\left\langle \nabla G(\xbf_{k})\ybf_{k+1},\cdot\right\rangle $
in Lemma~\ref{lem:descent}, the following relations
\begin{align}
-\langle\ybf_{k+1}-\ybf,\zbf_{k}\rangle & \leq A_{k+1},\label{eq:indicator}\\
f(\xbf_{k+1})+\left\langle \ybf_{k+1},\nabla G(\xbf_{k})^{\top}(\xbf_{k+1}-\xbf)\right\rangle  & \leq f(\xbf)+B_{k+1},\label{eq:finequality}
\end{align}
where
\begin{align}A_{k+1} & \triangleq\frac{1}{2\sigma_{k}}\left(\|\ybf-\ybf{}_{k}\|^{2}-\|\ybf-\ybf_{k+1}\|^{2}-\|\ybf_{k+1}-\ybf_{k}\|^{2}\right),\label{eq:defAk+1}\\
B_{k+1} & \triangleq\frac{1}{2\tau_{k}}\left(\|\xbf-\xbf_{k}\|^{2}-\|\xbf-\xbf_{k+1}\|^{2}-\|\xbf_{k+1}-\xbf_{k}\|^{2}\right)\label{eq:defBk+1},
\end{align}
hold for any $\xbf\in\Xcal$ and $\ybf\in \bigcap_{0\le s\le k}\mcal{Y}_s$. The existence of such $\ybf$ follows from our induction hypothesis. 
Since $\ybf_{k+1}^{\top}G(\cdot)$ is $\rho_{k}$-strongly convex, we have
\begin{align*}
 & \left\langle \ybf_{k+1},\nabla G(\xbf_{k})^{\top}(\xbf_{k+1}-\xbf)\right\rangle \\
% ={} & \left\langle \ybf_{k+1},\nabla G(\xbf_{k})^{\top}(\xbf_{k+1}-\xbf_{k})\right\rangle +\left\langle \ybf_{k+1},\nabla G(\xbf_{k})^{\top}(\xbf_{k}-\xbf)\right\rangle \\
& \ge{} \left\langle \ybf_{k+1},\nabla G(\xbf_{k})^{\top}(\xbf_{k+1}-\xbf_{k})\right\rangle\\ 
&\ \ +\left\langle \ybf_{k+1},G(\xbf_{k+1})-G(\xbf)\right\rangle -\left\langle \ybf_{k+1},G(\xbf_{k+1})-G(\xbf_{k})\right\rangle +\frac{\rho_{k}}{2}\|\xbf-\xbf_{k}\|^{2}.
\end{align*}
Combining this result and~\eqref{eq:finequality}, we have
\begin{equation}
\begin{split} & f(\xbf_{k+1})-f(\xbf)+\left\langle \ybf_{k+1},G(\xbf_{k+1})-G(\xbf)\right\rangle \\
 & \le{}B_{k+1}-\left\langle \ybf_{k+1},\nabla G(\xbf_{k})^{\top}(\xbf_{k+1}-\xbf_{k})\right\rangle +\left\langle \ybf_{k+1},G(\xbf_{k+1})-G(\xbf_{k})\right\rangle -\frac{\rho_{k}}{2}\|\xbf-\xbf_{k}\|^{2}.
\end{split}
\label{eq:temp-01}
\end{equation}
On the other hand, by the definition of $\zbf_{k}$, we have
\begin{align}
 & \langle\ybf-\ybf_{k+1},\zbf_{k}\rangle\nonumber \\
% ={} & \langle\ybf-\ybf_{k+1},G(\xbf_{k})\rangle+\theta_{k}\langle\ybf-\ybf_{k+1},G(\xbf_{k})-G(\xbf_{k-1})\rangle\nonumber \\
&  ={} \langle\ybf-\ybf_{k+1},G(\xbf_{k})-G(\xbf_{k+1})\rangle+\langle\ybf-\ybf_{k+1},G(\xbf_{k+1})\rangle\label{eq:temp-02}\\
 &\ \ +(\sigma_{k-1}/\sigma_k)\langle\ybf-\ybf_{k},G(\xbf_{k})-G(\xbf_{k-1})\rangle+(\sigma_{k-1}/\sigma_k)\langle\ybf_{k}-\ybf_{k+1},G(\xbf_{k})-G(\xbf_{k-1})\rangle.\nonumber 
\end{align}
Let us denote $\qbf_{k}=G(\xbf_{k})-G(\xbf_{k-1})$ for brevity. Combining~\eqref{eq:indicator}
and~\eqref{eq:temp-02} yields
\begin{equation}\label{eq:temp-03}
\begin{split} & \langle\ybf-\ybf_{k+1},G(\xbf_{k+1})\rangle\\
&  \le A_{k+1}+\langle\ybf-\ybf_{k+1},G(\xbf_{k+1})-G(\xbf_{k})\rangle-(\sigma_{k-1}/\sigma_k)\langle\ybf-\ybf_{k},\qbf_{k}\rangle-(\sigma_{k-1}/\sigma_k)\langle\ybf_{k}-\ybf_{k+1},\qbf_{k}\rangle.
\end{split}
\end{equation}
Putting~\eqref{eq:temp-01} and~\eqref{eq:temp-03} together, we
have
\begin{align*}
 & \mcal L(\xbf_{k+1},\ybf)-\mcal L(\xbf,\ybf_{k+1})\\
 % ={}& f(\xbf_{k+1})-f(\xbf)+\left\langle \ybf,G(\xbf_{k+1})\right\rangle -\left\langle \ybf_{k+1},G(\xbf)\right\rangle \\
 & \le{} A_{k+1}+B_{k+1}-\left\langle \ybf_{k+1},\nabla G(\xbf_{k})^{\top}(\xbf_{k+1}-\xbf_{k})\right\rangle +\left\langle \ybf_{k+1},G(\xbf_{k+1})-G(\xbf_{k})\right\rangle \\
 & \ \  +\langle\ybf-\ybf_{k+1},\qbf_{k+1}\rangle-(\sigma_{k-1}/\sigma_k)\langle\ybf-\ybf_{k},\qbf_{k}\rangle+(\sigma_{k-1}/\sigma_k)\langle\ybf_{k+1}-\ybf_{k},\qbf_{k}\rangle-\frac{\rho_{k}}{2}\|\xbf-\xbf_{k}\|^{2}\\
&  \le{}  A_{k+1}+B_{k+1}+\frac{L_{XY}}{2}\|\xbf_{k+1}-\xbf_{k}\|^{2}-\frac{\rho_{k}}{2}\|\xbf-\xbf_{k}\|^{2}\\
 & \ \  +\langle\ybf-\ybf_{k+1},\qbf_{k+1}\rangle-(\sigma_{k-1}/\sigma_k)\langle\ybf-\ybf_{k},\qbf_{k}\rangle+(\sigma_{k-1}/\sigma_k)\langle\ybf_{k+1}-\ybf_{k},\qbf_{k}\rangle,
\end{align*}
where the last inequality is by Lipschitz {smoothness} of $\left\langle \ybf_{k+1},G(\cdot)\right\rangle $.

Next, we bound the term $\langle\qbf_{k},\ybf_{k+1}-\ybf_{k}\rangle$
by Young's inequality, which gives
\begin{equation}
\begin{split}\inner{\ybf_{k+1}-\ybf_{k}}{\qbf_{k}} & \leq\frac{1}{2\sigma_{k-1}}\|\ybf_{k+1}-\ybf_{k}\|^{2}+\frac{\sigma_{k-1}}{2}\|\qbf_{k}\|^{2},\end{split}
\label{eq:upperbound_q}
\end{equation}
% where $1$ is a positive constant. 
It follows from~\eqref{eq:upperbound_q} and $\frac{\sigma_{k}}{2}\|\qbf_{k+1}\|^{2}\leq\frac{L_{G}^{2}\sigma_{k}}{2}\|\xbf_{k+1}-\xbf_{k}\|^{2}$
that
\begin{equation}
\label{eq:bounded_lagran}
\begin{aligned}
 & \mcal L(\xbf_{k+1},\ybf)-\mcal L(\xbf,\ybf_{k+1})\\
& \le  \frac{\tau_{k}^{-1}-\rho_{k}}{2}\|\xbf-\xbf_{k}\|^{2}-\frac{\tau_{k}^{-1}}{2}\|\xbf-\xbf_{k+1}\|^{2}+\frac{(\sigma_{k-1}/\sigma_k)\sigma_{k-1}}{2}\|\qbf_{k}\|^{2}-\frac{\sigma_{k}}{2}\|\qbf_{k+1}\|^{2}\\
 &\ \  +\frac{1}{2\sigma_{k}}\left(\|\ybf-\ybf{}_{k}\|^{2}-\|\ybf-\ybf_{k+1}\|^{2}\right)+\langle\ybf-\ybf_{k+1},\qbf_{k+1}\rangle-(\sigma_{k-1}/\sigma_k)\langle\ybf-\ybf_{k},\qbf_{k}\rangle\\
 &\ \  -\frac{\sigma_{k}^{-1}-(\sigma_{k-1}/\sigma_k)/\sigma_{k-1}}{2}\|\ybf_{k+1}-\ybf_{k}\|^{2}+\frac{L_{G}^{2}\sigma_{k}}{2}\|\xbf_{k+1}-\xbf_{k}\|^{2}-\frac{\tau_{k}^{-1}-L_{X Y}}{2}\|\xbf_{k+1}-\xbf_{k}\|^{2}.
\end{aligned}
\end{equation}
Multiply both sides of the above relation by $t_{k}$ and sum up the
result for $k=0,1,\ldots,K-1$. In view of the parameter relation~\eqref{eq:param-03}, we have
\begin{equation}\label{eq:temp-04}
    \begin{aligned}
 & \sum_{k=0}^{K-1}t_{k}\bsbra{\mathcal{L}(\xbf_{k+1},\ybf)-\mathcal{L}(\xbf,\ybf_{k+1})} \\
 % & \le\frac{t_{0}\rbra{\tau_{0}^{-1}-\rho_{0}}}{2}\|\xbf-\xbf_{0}\|^{2}-\frac{t_{K-1}\tau_{K-1}^{-1}}{2}\|\xbf-\xbf_{K}\|^{2}+\frac{t_{0}\theta_{0}}{2\alpha_{0}}\|\qbf_{0}\|^{2}-\frac{t_{K-1}}{2/\sigma_{k-1}}\|\qbf_{K}\|^{2}\\
 % & \quad+\frac{t_{0}\sigma_{0}^{-1}}{2}\|\ybf-\ybf_{0}\|^{2}-\frac{t_{K-1}\sigma_{K-1}^{-1}}{2}\|\ybf-\ybf_{K}\|^{2}+t_{K-1}\langle\ybf-\ybf_{K},\qbf_{K}\rangle-t_{0}\theta_{0}\langle\ybf-\ybf_{0},\qbf_{0}\rangle \\
 & \overset{(a)}{\leq}\frac{t_{0}\rbra{\tau_{0}^{-1}-\rho_{0}}}{2}\|\xbf-\xbf_{0}\|^{2}-\frac{t_{K-1}\tau_{K-1}^{-1}}{2}\|\xbf-\xbf_{K}\|^{2}
 -\frac{t_{K-1}\sigma_{K-1}}{2}\|\qbf_{K}\|^{2} \\
 & \quad+\frac{t_{0}\sigma_{0}^{-1}}{2}\|\ybf-\ybf_{0}\|^{2}-\frac{t_{K-1}\sigma_{K-1}^{-1}}{2}\|\ybf-\ybf_{K}\|^{2}+t_{K-1}\langle\ybf-\ybf_{K},\qbf_{K}\rangle-t_{0}\langle\ybf-\ybf_{0},\qbf_{0}\rangle \\
 & \overset{(b)}{\leq}\frac{1}{2\tau_{0}}\|\xbf-\xbf_{0}\|^{2}+\frac{1}{2\sigma_{0}}\|\ybf-\ybf_{0}\|^{2}-\frac{t_{K-1}\tau_{K-1}^{-1}}{2}\|\xbf-\xbf_{K}\|^{2}
 % -\frac{(1-1)t_{K-1}\sigma_{K-1}^{-1}}{2}\|\ybf-\ybf_{K}\|^{2},
\end{aligned}
\end{equation}
where $(a)$ uses $\qbf_{0}=\boldsymbol{0}$ and 
$\xbf_{-1}=\xbf_{0}$, and $(b)$ holds by $\rho_0= 0$,  $t_{0}=1$ and
\[t_{K-1}\left\langle \ybf-\ybf_{K},\qbf_{K}\right\rangle \le\frac{t_{K-1}}{2\sigma_{K-1}}\norm{\ybf-\ybf_{K}}^{2}+\frac{t_{K-1}}{2/\sigma_{K-1}}\norm{\qbf_{K}}^{2}.\] Since
$\mathcal{L}(\xbf,\ybf)$ is convex in $\xbf$ and linear in $\ybf$,
we have
\begin{equation}
T_{K}\bsbra{\mathcal{L}(\bar{\xbf}_{K},\ybf)-\mathcal{L}(\xbf,\bar{\ybf}_{K})}\le\sum_{k=0}^{K-1}t_{k}\bsbra{\mathcal{L}(\xbf_{k+1},\ybf)-\mathcal{L}(\xbf,\ybf_{k+1})},\label{eq:temp-05}
\end{equation}
% where $\bar{\xbf}_{k}=T_{k}^{-1}\sum_{s=0}^{k-1}t_{s}\xbf_{s+1}$
% and $\bar{\ybf}_{k}=T_{k}^{-1}\sum_{s=0}^{k-1}t_{s}\ybf_{s}$.
Combining~\eqref{eq:temp-04} and~\eqref{eq:temp-05}, we obtain
\begin{equation}\label{eq:mid-03}
\begin{aligned}
T_{K}\bsbra{\mathcal{L}(\bar{\xbf}_{K},\ybf)-\mathcal{L}(\xbf,\bar{\ybf}_{K})}\le &\frac{1}{2\tau_{0}}\|\xbf-\xbf_{0}\|^{2}-\frac{t_{K-1}\tau_{K-1}^{-1}}{2}\|\xbf-\xbf_{K}\|^{2}+\frac{1}{2\sigma_{0}}\|\ybf-\ybf_{0}\|^{2}.
\end{aligned}
\end{equation}
Dividing both sides by $T_{K},$ we obtain the desired result~\eqref{eq:thm2}.

Part 2. Next we show $\Ycal^*\subseteq\Ycal_{K}$. Let $\ybf^*$ be any point in $\Ycal^*$. Since \eqref{eq:mid-03} holds for any $\xbf\in\Xcal$ and $\ybf\in \cap_{0\le k\le K-1} \Ycal_k\supseteq \Ycal^*$, we can place $\xbf=\xbf^{*},\ybf=\ybf^{*}\in\Ycal^*$
in~\eqref{eq:thm2} to obtain
\[
\frac{t_{K-1}\tau_{K-1}^{-1}}{2T_{K}}\|\xbf^{*}-\xbf_{K}\|^{2}+\mathcal{L}(\bar{\xbf}_{K},\ybf^{*})-\mathcal{L}(\xbf^{*},\bar{\ybf}_{K})\leq\frac{1}{T_{K}}\Delta(\xbf^{*},\ybf^{*}).
\]
Moreover, the strong convexity of $\Lcal(\cdot, \ybf^*)$ implies
\[\mathcal{L}(\bar{\xbf}_{K},\ybf^{*}) \geq\mathcal{L}(\xbf^{*},\ybf^{*})+\frac{(\ybf^{*})^{\top}\boldsymbol{\mu}}{2}\|\bar{\xbf}_{K}-\xbf^{*}\|^{2}\ge\mathcal{L}(\xbf^{*},\bar{\ybf}_{K})+\frac{(\ybf^{*})^{\top}\boldsymbol{\mu}}{2}\|\bar{\xbf}_{K}-\xbf^{*}\|^{2}.
\]
Applying the above two inequalities yields
\begin{equation}
\frac{(\ybf^{*})^{\top}\boldsymbol{\mu}}{2}\|\bar{\xbf}_{K}-\xbf^{*}\|^{2} \leq\frac{1}{T_{K}}\Delta(\xbf^{*},\ybf^{*}),\ %
\frac{1}{2}\|\xbf_{K}-\xbf^{*}\|^{2}  \leq\frac{\tau_{K-1}\sigma_{0}}{\sigma_{K-1}}\Delta(\xbf^{*},\ybf^{*}).\label{eq:thm2-temp04}
\end{equation}
In view of \eqref{eq:thm2-temp04} and Proposition~\ref{prop:local_estimation}, we have that 
\[
(\ybf^*)^T\mubm \ge \underline{\mu}\, \norm{\ybf^*}_1 = \underline{\mu} \max\bcbra{h_1(\xbf_K, \frac{\sigma_0\tau_{K-1}\Delta_{XY}}{\sigma_{K-1}}), h_2(\bar{\xbf}_K, \frac{\Delta_{XY}}{T_K})}\coloneqq \hat{\rho}_K.
\]
Moreover, since $\Ycal^*\subseteq \Ycal_{K-1}$, we have $(\ybf^*)^T\mubm \ge \rho_{K-1}$. 
Hence we have $(\ybf^*)^T\mubm \ge \rho_{K}$ where $\rho_{K}=\max\{\hat{\rho}_K, \rho_{K-1}\}$ is the output of the \textsc{Improve} procedure.  Due to the construction of $\Ycal_K$,  we immediately see that $\ybf^*\in\Ycal_K$. This implies $\Ycal^*\subseteq \Ycal_K$ and completes our induction proof.
\endproof

% \subsection{Proof of Lemma~\ref{lem:speed}}

Next, we specify the stepsize selection in Lemma~\ref{lem:speed} and develop more concrete complexity results in Corollary~\ref{cor:mainthm}.

\begin{lemma}\label{lem:speed}
Let $\hat{\rho}_{k+1}:=\frac{\sqrt{\hat{\rho}_{k}^{2}k^{2}+(3\rho_{k+1}\hat{\rho}_{k})k}}{k+1}$
for $k\geq1$ and $\hat{\rho}_{1}=3\sqrt{\frac{\rho_{1}}{\tau_{0}}}$.
Suppose $\sigma_{k},\tau_{k}$ satisfy:
\begin{equation}\label{eq:step-size-1}
\tau_{0}^{-1}\geq L_{XY}+L_{G}^{2}\sigma_{0},\ \ \tau_{k+1}=\tau_{k}(1+\rho_{k+1}\tau_{k})^{-\frac{1}{2}},\ \ \sigma_{k+1}=\frac{\tau_{k}\sigma_{k}}{\tau_{k+1}}.
\end{equation}
Then we have
\begin{equation}\label{speed_lower}
\frac{1}{\tau_{k}^{2}}\geq\frac{\hat{\rho}_{k}^{2}}{9}k^{2}+\frac{1}{\tau_{0}^{2}},\ \ T_{k}\geq1+\frac{\tau_{0}}{6}\tilde{\rho}_{k}(k+1)k,\ \ \hat{\rho}_{k}\geq\min\{\rho_{1},\hat{\rho}_{1}\},
\end{equation}
where $\tilde{\rho}_{k}=2\sum_{s=0}^{k}\frac{\hat{\rho}_{s}s}{k(k+1)}$
for $k\geq1$. Moreover, suppose $\bar{\rho}\tau_{0}\leq2$, where
$\bar{\rho}=\bar{c}\cdot\bar{\mu}$, then we have
\begin{equation}\label{eq:gammak_upper}
\sigma_{k}^{2}\leq\sigma_{0}^{2}(k+1)^{2}.
\end{equation}
\end{lemma}
\proof
We first use induction to show that $\frac{1}{\tau_{k}^{2}}\geq\frac{\hat{\rho}_{k}^{2}}{9}k^{2}+\frac{1}{\tau_{0}^{2}}$.
It is easy to see that $\frac{1}{\tau_{k}^{2}}\geq\frac{\hat{\rho}_{k}^{2}}{9}k^{2}+\frac{1}{\tau_{0}^{2}}$
holds for $k=1$ by the definition $\hat{\rho}_{1}=3\sqrt{\rho_{1}/\tau_{0}}$
and $\tau_{1}=\tau_{0}(1+\rho_{1}\tau_{0})^{-\frac{1}{2}}$. Assume
$\frac{1}{\tau_{k}^{2}}\geq\frac{\hat{\rho}_{k}^{2}}{9}k^{2}+\frac{1}{\tau_{0}^{2}}$
holds for all $k=0,\ldots,K$, then we have
\begin{equation}
\begin{split}\frac{1}{\tau_{K+1}^{2}} & =\frac{1}{\tau_{K}^{2}}+\frac{\rho_{K+1}}{\tau_{K}}\\
 & \geq\frac{\hat{\rho}_{K}^{2}}{9}K^{2}+\frac{1}{\tau_{0}^{2}}+\rho_{K+1}\sqrt{\frac{\hat{\rho}_{K}^{2}}{9}K^{2}+\frac{1}{\tau_{0}^{2}}}\\
 & \geq\frac{\hat{\rho}_{K}^{2}}{9}K^{2}+\frac{1}{\tau_{0}^{2}}+\frac{\rho_{K+1}\hat{\rho}_{K}K}{3}\\
 & \geq\frac{\hat{\rho}_{K+1}^{2}}{9}(K+1)^{2}+\frac{1}{\tau_{0}^{2}},
\end{split}
\end{equation}
which completes our induction. It follows from $\frac{1}{\tau_{k}^{2}}\geq\frac{\hat{\rho}_{k}^{2}}{9}k^{2}+\frac{1}{\tau_{0}^{2}}$
and the relation among $T_{k},t_{k},\sigma_{k},\tau_{k}$ that, for
any $k\geq1$
\begin{equation}
\begin{aligned}
    T_{k}=&\sum_{s=0}^{k-1}t_{s}=1+\sum_{s=1}^{k-1}t_{s}\geq1+\sum_{s=1}^{k-1}\frac{\sigma_{s}}{\sigma_{0}}=1+\sum_{s=1}^{k-1}\frac{\tau_{0}}{\tau_{s}}\geq1+\tau_{0}\sum_{s=1}^{k-1}\sqrt{\frac{\hat{\rho}_{s}^{2}s^{2}}{9}+\frac{1}{\tau_{0}^{2}}}\\
    &>1+\tau_{0}\sum_{s=1}^{k-1}\frac{\hat{\rho}_{s}s}{3}=1+\frac{\tau_{0}}{6}\tilde{\rho}_{k}(k+1)k.
\end{aligned}
\end{equation}

Similarly, we use induction to prove

\begin{equation}
\hat{\rho}_{k}\geq\min\{\rho_{1},\hat{\rho}_{1}\},\forall k\geq1.\label{eq:induction}
\end{equation}

It is easy to find that $\hat{\rho}_{1}\geq\min\{\rho_{1},\hat{\rho}_{1}\}$.
We assume that $\hat{\rho}_{k}\geq\min\{\rho_{1},\hat{\rho}_{1}\},\forall k\geq1$
holds for any $k=1,\ldots,K$. Considering $\hat{\rho}_{K+1}$, we
have 
\[
\begin{split}\hat{\rho}_{K+1} & \geq\frac{1}{K+1}\sqrt{\hat{\rho}_{K}^{2}K^{2}+3\rho_{1}\hat{\rho}_{K}K}\\
 & \geq\frac{1}{K+1}\sqrt{\left(\min\left\{\rho_{1},\hat{\rho}_{1}\right\}\right)^{2}K^{2}+3\rho_{1}\cdot\min\left\{\rho_{1},\hat{\rho}_{1}\right\}K}\geq\min\left\{\rho_{1},\hat{\rho}_{1}\right\},
\end{split}
\]
which completes the induction. Moreover, we use induction to show
$\sigma_{k}^{2}\leq\sigma_{0}^{2}(k+1)^{2}$. It is obvious that the
inequality holds for $k=0$. Assume the inequality holds for all $k=0,\ldots,K,$
then we have 
\begin{equation}
\begin{split}\sigma_{K+1}^{2} & =\sigma_{K}^{2}(1+\rho_{K+1}\frac{\tau_{0}\sigma_{0}}{\sigma_{K}})\\
 & =\sigma_{K}^{2}+\rho_{K+1}\tau_{0}\sigma_{0}\sigma_{K}\\
 & \leq\sigma_{0}^{2}\left((K+1)^{2}+\rho_{K+1}\tau_{0}(K+1)\right)\\
 & \leq\sigma_{0}^{2}(K+2)^{2},
\end{split}
\end{equation}
where the last inequality use the relation $\rho_{k}\leq\bar{\rho},\forall k$,
and $\bar{\rho}\tau_{0}\leq2$.
\endproof

\subsection{Proof of Corollary~\ref{cor:mainthm}\label{proof:cor}}
% \begin{recorollary}{\ref{cor:mainthm}}
% \label{cor:appendix_mainthm}
% Suppose that $\gamma_k, \sigma_k,\tau_k,\theta_{k}, t_k$ satisfy~\eqref{eq:step-size-0},~\eqref{eq:step-size-1} and~\eqref{eq:step-size-2}.
% Then, for any $K\geq 1$, we have 
% \begin{equation}\label{eq:appendix_cor_complexity}
% \begin{aligned}
% f(\bar{\xbf}_{K})-f(\xbf^{*})
% &\leq  \frac{6}{6+\tau_{0}\trho_{K} (K+1)K}\Brbra{\frac{1}{2\tau_0}\norm{\xbf_0-\xbf^*}^2+\frac{D_Y^2}{2\sigma_0}},\\
% \|[G(\bar{\xbf}_{K})]_{+}\| & \leq \frac{6}{c^{*}(6+\tau_{0}\trho_{K} (K+1)K)}\Brbra{\frac{1}{2\tau_0}\norm{\xbf_0-\xbf^*}^2+\frac{D_Y^2}{2\sigma_0}},\\
% \frac{1}{2}\|\xbf_{K}-\xbf^{*}\|^{2}&\leq\frac{3\sigma_{0}}{\vrho_{K}^{2}\tau_{0}^{2}K^{2}+9\gamma_{0}}\Delta(\xbf^*,\ybf^*),
% \end{aligned}
% \end{equation}
% {where $c^{*}:=\brbra{f(\xbf^{*})-\min_{\xbf}f(\xbf)}/{\min_{i\in[m]}\{-g_{i}(\tilde{\xbf})\}}>0$.}
% \end{recorollary}

\proof
First, we show that the sequences $\{\tau_{k},\sigma_{k},t_{k},\rho_{k}\}$
generated by {\APDSCE} satisfy the relationship in~\eqref{eq:param-03} in Theorem~\ref{thm:gapdiminishing}. 
The first part of~\eqref{eq:param-03} can be derived  using the monotonicity of $\{\rho_k\}$ as follows:
\begin{equation*}
    \begin{aligned}
    t_{k+1}\brbra{\tau_{k+1}^{-1}-\rho_{k+1}}&=\sigma_0^{-1}\brbra{\sigma_{k+1}\tau_{k+1}-\sigma_{k+1}\rho_{k+1}} \\
    &=\sigma_0^{-1}\brbra{\sigma_{k}\tau_k \tau_{k+1}^{-2}-\sigma_{k+1}\rho_{k+1}} \\
    &=\sigma_0^{-1}\brbra{\sigma_k(1+\rho_{k+1}\tau_k)/\tau_k - \sigma_{k+1}\rho_{k+1}}\\
    &=\sigma_0^{-1}\brbra{\sigma_k / \tau_{k} + \rho_{k+1}\sigma_k - \sigma_{k+1}\rho_{k+1}}\\
    &\leq t_k\tau_k^{-1}
    \end{aligned}
\end{equation*}
% \[
% \begin{aligned}
% t_{k+1}(\tau_{k+1}^{-1}-\rho_{k+1}) 
% % & =\sigma_0^{-1}\sigma_{k+1}(\tau_{k+1}^{-1}-\rho_{k+1})\\
% % & = \sigma_0^{-1}(\gamma_{k+1}-\gamma_{k+1}\rho_{k+1}\tau_{k+1})\\
% = \sigma_0^{-1}(\gamma_{k}+\gamma_k\rho_{k+1}\tau_k-\sqrt{\gamma_{k+1}\gamma_k}\rho_{k+1}\tau_{k})
% \le  \sigma_0^{-1}\gamma_{k} = t_k\tau_k^{-1}.
% \end{aligned}
% \]
% \[
% \begin{aligned}
% t_{k+1}(\tau_{k+1}^{-1}-\rho_{k+1}) 
% & =\sigma_0^{-1}\sigma_{k+1}(\tau_{k+1}^{-1}-\rho_{k+1})\\
% % & = \sigma_0^{-1}(\gamma_{k+1}-\gamma_{k+1}\rho_{k+1}\tau_{k+1})\\
% & = \sigma_0^{-1}(\gamma_{k}+\gamma_k\rho_{k+1}\tau_k-\sqrt{\gamma_{k+1}\gamma_k}\rho_{k+1}\tau_{k})\\
% & \le  \sigma_0^{-1}\gamma_{k} = t_k\tau_k^{-1}.
% \end{aligned}
% \]
The second part of \eqref{eq:param-03} can be easily verified using the parameters setting. 

Next, we prove the last term in~\eqref{eq:param-03} by induction. Firstly, it easy to verify that for any $\sigma_{0}>0$, there exists $\tau_{0}\in(0,(L_{XY}+L_G^2 \sigma_{0})^{-1}]$ such that last term of~\eqref{eq:param-03} holds.
Hence, when $k=0$, the last term of~\eqref{eq:param-03} is directly from the first term of~\eqref{eq:cor_param}.
Suppose that the last term of~\eqref{eq:param-03} holds for $k=0,\ldots,K-1$.
From ${\sigma_{K-1}}/{\sigma_{K}}={\tau_{K}}/{\tau_{K-1}}\leq 1$,
we have
\begin{equation}\label{eq:relation_tau}
    \frac{1}{\tau_{K}}=\frac{\sigma_{K}}{\tau_{K-1}\sigma_{K-1}}\geq\frac{L_{XY}}{\sigma_{K-1}/\sigma_{K}}+L_{G}^{2}\sigma_{K} \geq L_{XY}+L_{G}^{2}\sigma_{K}.
\end{equation}

Without loss of generality, place $\xbf=\xbf^{*}$, $\ybf=\ybf^{+}:=\left(\|\ybf^{*}\|_{1}+c^{*}\right)\frac{[G(\bar{\xbf}_{K})]_{+}}{\|[G(\bar{\xbf}_{K})]_{+}\|}$
in~\eqref{eq:thm2}, and using $\norm{\ybf^*}_1\leq \bar{c}$ in Proposition~\ref{prop:bounded_set_y}. It is easy to see
$
\|\ybf^{+}\|=\|\ybf^{*}\|_{1}+c^{*}\leq \bar{c}$, 
and 
$
\norm{\ybf^{+}}_{1}\geq\norm{\ybf^{+}}=\|\ybf^*\|_1 + c^*\geq\|\ybf^*\|_1,
$
Hence, we conclude that $\ybf^{+}\in\mathcal{Y}_{k},\forall k\geq0$.

Now observe that $\mathcal{L}(\bar{\xbf}_{K},\ybf^{*})-\mathcal{L}(\xbf^{*},\ybf^{*})\geq0$,
which implies 
$
f(\bar{\xbf}_{K})+\left\langle \ybf^{*},G(\bar{\xbf}_{K})\right\rangle -f(\xbf^{*})\geq 0.
$
In view of $\left\langle \ybf^{*},G(\bar{\xbf}_{K})\right\rangle \leq\left\langle \ybf^{*},[G(\bar{\xbf}_{K})]_{+}\right\rangle \leq\|\ybf^{*}\| \cdot\|[G(\bar{\xbf}_{K})]_{+}\|$,
then we have
\begin{equation}
f(\bar{\xbf}_{K})+\|\ybf^{*}\|\cdot\|[G(\bar{\xbf}_{K})]_{+}\|-f(\xbf^{*})\geq0.\label{eq:infe1}
\end{equation}
Moreover, it follows from $\norm{\ybf^*}_1\geq \norm{\ybf^*}$ that
\begin{equation}
\begin{split}\mathcal{L}(\bar{\xbf}_{K},\ybf^{+})-\mathcal{L}(\xbf^{*},\bar{\ybf}_{K}) & \geq\mathcal{L}(\bar{\xbf}_{K},\ybf^{+})-\mathcal{L}(\xbf^{*},\ybf^{*})\\
 % & =f(\bar{\xbf}_{K})+ \left(\|\ybf^{*}\|_{1}+c^{*}\right)\|[G(\bar{\xbf}_{K})]_{+}\|-f(\xbf^{*})\\
 & \geq f(\bar{\xbf}_{K})+\left(\|\ybf^{*}\|+c^{*}\right)\|[G(\bar{\xbf}_{K})]_{+}\|-f(\xbf^{*}).
\end{split}
\label{eq:infe2}
\end{equation}
Combining~\eqref{eq:infe1},~\eqref{eq:infe2} and~\eqref{eq:thm2},
we obtain
\begin{equation}\label{eq:bound-01}
\max\bcbra{ c^{*}\|[G(\bar{\xbf}_{K})]_{+}\|,f(\bar{\xbf}_{K})-f(\xbf^{*})}
% \leq\frac{1}{T_{K}}\Delta(\xbf^{*},\ybf^{+})
\leq{\frac{1}{T_{K}}\brbra{\frac{1}{2\tau_0}\norm{\xbf_0-\xbf^*}^2+\frac{D_Y^2}{2\sigma_0}}},
\end{equation}
In view of the bound in~\eqref{speed_lower} and the relation between $\tau_{k},\sigma_{k}$, we can get  
\begin{equation}
\label{speed:extend}
\frac{\tau_{k}}{\sigma_{k}} \le \frac{3}{\vrho_{k}^{2}\tau_{0}^{2}k^{2}+9\sigma_0/\tau_0}.
\end{equation}
In view of \eqref{eq:bound-01} and \eqref{speed_lower}, we have
\[
\max\bcbra{ c^{*}\|[G(\bar{\xbf}_{K})]_{+}\|,f(\bar{\xbf}_{K})-f(\xbf^{*})}
\leq  \frac{6}{6+\tau_{0}\trho_{K} (K+1)K}{\brbra{\frac{1}{2\tau_0}\norm{\xbf_0-\xbf^*}^2+\frac{D_Y^2}{2\sigma_0}}}.
\]
Combining~\eqref{eq:thm2} and~\eqref{speed:extend} yields
$
\frac{1}{2}\|\xbf_{K}-\xbf^{*}\|^{2}
\leq{3\sigma_{0}}{\Delta(\xbf^*,\ybf^*)}/({\vrho_{K}^{2}\tau_{0}^{2}K^{2}+9\sigma_0/\tau_0}).
$
\endproof

\section{\label{sec:Proof Details in restart}
Convergence analysis of \resapdsce{}
% Proof details in Section~\ref{sec:restarted}
}

\subsection{Proof of Theorem~\ref{thm:rapdpro}\label{sec:proof_thm_rapdpro}}
% \begin{retheorem}{\ref{thm:rapdpro}}
% \label{thm:appendix_rapdpro}
% Let $\{\xbf_{0}^{s}\}_{s\geq0}$
% be the sequence generated by {\resapdsce}, then we have
% \begin{equation}
% \|\xbf_{0}^{s}-\xbf^{*}\|^{2}\leq\Delta_{s}\equiv D_{X}^{2}\cdot 2^{-s},\ \ \ \forall s\geq0.\label{appendix:eq:exp decay}
% \end{equation}
% As a consequence, {\resapdsce}
% will find a solution $\xbf_{0}^{S}$ such that $\|\xbf_{0}^{S}-\xbf^{*}\|^{2}\leq\vep$
% for any $\vep\in(0, D_{X}^{2})$ in at most $S:=\big\lceil \log_{2}\rbra{ D_{X}^{2}/\vep}\big\rceil $
% epochs. Moreover, the overall number of iterations performed by {\resapdsce}
% to find such a solution is bounded by
% \begin{equation}\label{eq:appendix_T_e}
%     T_{\vep}:= \Brbra{\frac{12}{\varpi_{1}\tau_{0}^{s}}+2}\left\lceil \log_{2}\frac{D_{X}}{\sqrt{\vep}}+1\right\rceil +\Brbra{\frac{6(\sqrt{2}+2)}{\varpi_{2}\sqrt{\tau_{0}^{s}\sigma_{0}^{s}}}}\cdot\Brbra{\frac{D_{Y}}{\sqrt{\vep}}},
% \end{equation}
% where $\varpi_1$ and $\varpi_2$ satisfy $\sum_{s=0}^{S}(\hat{\rho}_{N_{s}}^s)^{-1} = (\varpi_1)^{-1}(S+1)$ and $\sum_{s=0}^{S} {\sqrt{2}^s}/{\hat{\rho}_{N_{s}}^s}= (\varpi_2)^{-1}\sum_{s=0}^{S}\sqrt{2}^s$, respectively.
% \end{retheorem}
\proof
First, we show that the choice of $\tau_{0}^{s}=\bar{\tau},\sigma_{0}^{s}=\bar{\sigma},\forall s\geq0$
satisfy the condition~\eqref{eq:cor_param} in Corollary~\ref{cor:mainthm}:
$
(\tau_{0}^{s})^{-1}\geq \rbra{1-\nu_{0}}\rbra{\tau_{0}^{s}}^{-1}= L_{XY}+c{L_{G}^{2}}\sigma_{0}^{s}/{\delta}\geq L_{XY}+c{L_{G}^{2}}\sigma_{0}^{s}
$.

Next, we show~\eqref{eq:exp decay} holds by induction. Clearly,~\eqref{eq:exp decay}
holds for $s=0$. Assume $\|\xbf_{0}^{s}-\xbf^{*}\|^{2}\leq\Delta_{s}$
holds for $s=0,\ldots,S-1$. Then by Theorem~\ref{thm:gapdiminishing},
we have 
\begin{equation}
\|\xbf_{0}^{S}-\xbf^{*}\|^{2} \leq\frac{\sigma_{0}^{S}\tau_{N_{S}}^{S}}{\sigma_{N_{S}}^{S}}\Big(\frac{2}{\tau_{0}^{S}}\Delta_{S}+\frac{1}{\sigma_{0}^{S}}D_{Y}^{2}\Big).
\label{eq:restart-01}
\end{equation}
% \begin{equation}
% \begin{split}\|\xbf_{0}^{S}-\xbf^{*}\|^{2} & \leq\frac{\sigma_{0}^{S}\tau_{N_{S}}^{S}}{\sigma_{N_{S}}^{S}}\big(\frac{1}{\tau_{0}^{S}}\|\xbf_{0}^{S-1}-\xbf^{*}\|^{2}+\frac{1}{\sigma_{0}^{S}}\|\ybf_{0}^{S-1}-\ybf^{*}\|^{2}\big)\\
%  & \leq\frac{\sigma_{0}^{S}\tau_{N_{S}}^{S}}{\sigma_{N_{S}}^{S}}\big(\frac{2}{\tau_{0}^{S}}\Delta_{S}+\frac{1}{\sigma_{0}^{S}}D_{Y}^{2}\big).
% \end{split}
% \label{eq:restart-01}
% \end{equation}
In view of the first bound in~\eqref{speed_lower} and the relation between $\tau^s_{N_s},\sigma^s_{N_s}$, we can get  
\begin{equation}
\frac{\tau_{N_{s}}^{s}}{\sigma_{N_{s}}^{s}}\leq\frac{9}{\sigma_{0}^{s}\tau_{0}^{s}(\hat{\rho}_{N_{s}}N_{s})^{2}}.\label{eq:restart-02}
\end{equation}
Combining~\eqref{eq:restart-01} and~\eqref{eq:restart-02} yields
\begin{equation*}
\norm{\xbf_{0}^{S}-\xbf^{*}}^{2}\leq\frac{18}{(\hat{\rho}_{N_{s}}\tau_{0}^{s}N_{s})^{2}}+\frac{9D_{Y}^{2}}{\sigma_{0}^{s}\tau_{0}^{s}(\hat{\rho}_{N_{s}}N_{s})^{2}}.
\end{equation*}
Since the algorithm sets $N_{s}=\lceil\max\{6(\hat{\rho}_{N_{s}}\tau_{0}^{s})^{-1},\sqrt{2}^{s}\cdot3\sqrt{2}D_{Y}/\brbra{\hat{\rho}_{N_{s}}D_{X}\sqrt{\tau_{0}^{s}\sigma_{0}^{s}}}\}\rceil$,  it follows that
\begin{equation*}
\begin{split}\frac{18}{(\hat{\rho}_{N_{s}}\tau_{0}^{s}N_{s})^{2}} & \leq\frac{18}{(\hat{\rho}_{N_{s}}\tau_{0}^{s})^{2}}\cdot\frac{(\hat{\rho}_{N_{s}}\tau_{0}^{s})^{2}}{36}=\frac{1}{2},\\
\frac{9D_{Y}^{2}}{\sigma_{0}^{s}\tau_{0}^{s}(\hat{\rho}_{N_{s}}N_{s})^{2}} & \leq\frac{9D_{Y}^{2}}{\sigma_{0}^{s}\tau_{0}^{s}\hat{\rho}_{N_{s}}^{2}}\cdot\frac{\hat{\rho}_{N_{s}}^{2}\sigma_{0}^{s}\tau_{0}^{s}D_{X}^{2}}{18D_{Y}^{2}2^{s}}=\frac{1}{2}\cdot2^{-s}D_{X}^{2}=\frac{1}{2}\Delta_{S},
\end{split}
\end{equation*}
which implies the desired result~\eqref{eq:exp decay}.

Let the algorithm run for $S=\big\lceil \log_{2}(D_{X}^{2}/{\vep})\big\rceil$ epochs, then $\|\xbf_{0}^{S}-\xbf^{*}\|^{2}\le D_{X}^{2}\cdot2^{-S}\leq\vep$. 
% It follows from that property~\eqref{speed_lower} that the sequence $\vrho_{N_{s}}^{s}$ is lower bounded, and then we  denote $\underline{\vrho}:=\min_{s}\{\vrho_{N_{s}}^{s}\}$.
The total iteration number required by Algorithm~\ref{alg:Restart Version of Stage Two TPAPD} for attaining a solution $\xbf_0^S$ such that $\norm{\xbf_0^S-\xbf^*}^2\leq \vep$ is 
\begin{equation*}
\begin{split}\sum_{s=0}^{S}N_{s}\leq & \sum_{s=0}^{S}\Bcbra{\frac{6}{\hat{\rho}^s_{N_{s}}\tau_{0}^{s}}+\frac{3\sqrt{2}D_{Y}}{\hat{\rho}^s_{N_{s}}D_{X}\sqrt{\tau_{0}^{s}\sigma_{0}^{s}}}\sqrt{2}^{s}+1}\\
\aeq & \Brbra{\frac{6}{\varpi_{1}\tau_{0}^{s}}+1}\rbra{S+1}+{\frac{3\sqrt{2}D_{Y}}{\varpi_{2}D_{X}\sqrt{\tau_{0}^{s}\sigma_{0}^{s}}}}\sum_{s=0}^{S}\sqrt{2}^{s}\\
\leq & \Brbra{\frac{12}{\varpi_{1}\tau_{0}^{s}}+2}\left\lceil \log_{2}\frac{D_{X}}{\sqrt{\vep}}+1\right\rceil +{\frac{3\sqrt{2}D_{Y}}{\varpi_{2}D_{X}\sqrt{\tau_{0}^{s}\sigma_{0}^{s}}}}\cdot\frac{\sqrt{2}^{S+1}-1}{\sqrt{2}-1}\\
\leq & \Brbra{\frac{12}{\varpi_{1}\tau_{0}^{s}}+2}\left\lceil \log_{2}\frac{D_{X}}{\sqrt{\vep}}+1\right\rceil +{\frac{3\sqrt{2}D_{Y}(\sqrt{2}+1)}{\varpi_{2}D_{X}\sqrt{\tau_{0}^{s}\sigma_{0}^{s}}}}\cdot\brbra{\sqrt{2}^{\log_{2}(D_{X}^{2}/\vep)+2}-1}\\
\leq & \Brbra{\frac{12}{\varpi_{1}\tau_{0}^{s}}+2}\left\lceil \log_{2}\frac{D_{X}}{\sqrt{\vep}}+1\right\rceil +{\frac{6D_{Y}(\sqrt{2}+2)}{\varpi_{2}\sqrt{\tau_{0}^{s}\sigma_{0}^{s}}}}\cdot{\frac{1}{\sqrt{\vep}}},
\end{split}
\end{equation*}
where $(a)$ holds by $\sum_{s=0}^{S}(\hat{\rho}_{N_{s}}^s)^{-1} = (\varpi_1)^{-1}(S+1)$ and $\sum_{s=0}^{S} {\sqrt{2}^s}/{\hat{\rho}_{N_{s}}^s}= (\varpi_2)^{-1}\sum_{s=0}^{S}\sqrt{2}^s$.
\endproof

% \subsection{Proof of Lemma~\ref{lem:asy}}
Now, we give some proof details in dual convergence results. Let 
\begin{equation*}
    \begin{split}Q_{j}(\xbf,\ybf) & :=\frac{(\tau_{j})^{-1}-\rho_{j}}{2}\norm{\xbf-\xbf_{j}}^{2}+\frac{1}{2\sigma_{j}}\norm{\ybf-\ybf_{j}}^{2}+(\sigma_{j-1}/\sigma_{j})\left\langle \ybf_{j}-\ybf,G(\xbf_{j})-G(\xbf_{j-1})\right\rangle \\
 & \ \ +\frac{(\sigma_{j-1}/\sigma_{j})}{2/\sigma_{j-1}}\norm{G(\xbf_{j})-G(\xbf_{j-1})}^{2},
\end{split}
\end{equation*}
then we establish an important property about the solution sequence in the following lemma.
\begin{lemma}
\label{lem:asy}
Assume $\bar{\tau}^{-1}>\overline{\rho}$
and  choose  $\nu_{0}>0$ such that 
\begin{equation}
1>\inf_{j\ge 0}\{\sigma_{j-1}/\sigma_{j}\}\geq\delta+\nu_{0}.\label{eq:appendix_ini_step_assu}
\end{equation}
 Then  there exists an $\nu_{1}>0$ such that  for any $j\ge 0$ and any KKT point $(\xbf^*,\tilde{\ybf}^*)$:
\begin{equation*}
\begin{aligned}
    0& \leq t_{j}Q_{j}(\xbf^{*},\tilde{\ybf}^{*})-t_{j+1}Q_{j+1}(\xbf^{*},\tilde{\ybf}^{*}) -\nu_{1}t_{j} \Bsbra{\frac{1}{2\tau_{j}}\norm{\xbf_{j+1}-\xbf_{j}}^{2}+\frac{1}{2\sigma_{j}}\norm{\ybf_{j+1}-\ybf_{j}}^{2}},\\ 
     0 & < t_{j}Q_{j}(\xbf^{*},\tilde{\ybf}^{*}).
\end{aligned}
\end{equation*}
\end{lemma}
\proof
First, we give some results that will be used repeatedly in the following. For notation simplicity, we denote $\theta_j = \sigma_{j-1}/\sigma_{j}$. In view of Lemma~\ref{lem:speed}, and the parameter ergodic
sequence generated by {\resapdsce}, we have $\left\{(\tau_{k}^s)^{-1},\sigma_{k}^s\right\}$ is monotonically increasing sequence in $k$, $\bar{\tau}={\tau_0^{s}},\bar{\sigma}={\sigma_0^{s}}, t_0^s=1,\forall s\geq 0$, and there exist a $\nu_{3}>0$ such that $\bar{\sigma}+\nu_{3} \leq \underline{\sigma}:= \min_{s}\{\sigma_{N_s}^s\}$.
Now, for {\resapdsce}, we claim that {there exist $\nu_{1},\nu_{2}>0$ such that the following two conditions hold} 

1. For any $j\ge 0$,  we have
\begin{equation}
\min\big\{1-\delta,(\tau_{j}^{-1}-L_{XY}-L_{G}^{2}\sigma_{j})\tau_{j}\big\}\geq\nu_{1}>0,\label{eq:nu_1_ineq}
\end{equation}
and 
\begin{equation}
\begin{split}t_{j}\min\Bcbra{\tau_{j}^{-1}-\rho_{j},\frac{1}{\sigma_{j}}-\frac{\delta}{\sigma_{j-1}}} & \geq\nu_{2}>0. \end{split}
\label{eq:nu_2_ineq}
\end{equation}

2. For any $j\geq 0$,
we have
\begin{equation}
0\leq t_{j}Q_{j}(\xbf^{*},\tilde{\ybf}^{*})-t_{j+1}Q_{j+1}(\xbf^{*},\tilde{\ybf}^{*})-\nu_{1}t_{j}\big({(2\tau_{j})^{-1}}\norm{\xbf_{j+1}-\xbf_{j}}^{2}+{(2\sigma_{j})^{-1}}\norm{\ybf_{j+1}-\ybf_{j}}^{2}\big).\label{eq:res_lem}
\end{equation}

Part 1. {We first consider two subsequent points $\xbf_j$ and $\xbf_{j+1}$ within the same epoch, and assume $j\sim (s, k)$. Then, it follows from $\theta_{k}^{s}={\sigma_{k-1}^{s}}/{\sigma_{k}^{s}}$
that 
\begin{equation}
(\sigma_{k}^{s})^{-1}-\theta_{k}^{s}\delta(\sigma_{k-1}^{s})^{-1}=(\sigma_{k}^{s})^{-1}-\delta(\sigma_{k}^{s})^{-1}=\frac{1-\delta}{\sigma_{k}^{s}}\overset{\eqref{eq:appendix_ini_step_assu}}{\geq}\frac{\nu_{0}}{\sigma_{k}^{s}}.\label{mid_01}
\end{equation}
Next, we use induction to show 
\begin{equation}
\frac{1-\nu_{0}}{\tau_{k}^{s}}\geq L_{XY}+L_{G}^{2}\sigma_{k}^{s}\delta^{-1}.\label{eq:new_relationship}
\end{equation}
When $k=0$, inequality~\eqref{eq:new_relationship} degenerates as the {definition of $\tau_{0}^s,\sigma_{0}^s$}. Suppose~\eqref{eq:new_relationship}
holds for $k=0,1,\ldots,K-1$. Then, from $\theta_{K}^{s}={\sigma_{K-1}^{s}}/{\sigma_{K}^{s}}={\tau_{K}^{s}}/{\tau_{K-1}^{s}}\leq1$,
we have
\[
(1-\nu_{0})(\tau_{K}^{s})^{-1}=(1-\nu_{0})(\tau_{K-1}^{s}\theta_{K}^{s})^{-1}\geq\frac{L_{XY}}{\theta_{K}^{s}}+\frac{L_{G}^{2}\sigma_{K-1}^{s}\delta^{-1}}{\theta_{K}^{s}}\geq L_{XY}+L_{G}^{2}\sigma_{K}^{s}\delta^{-1},
\]
which completes our induction proof. Hence, combining~\eqref{mid_01}
and~\eqref{eq:new_relationship}, we have
\begin{equation}
\min\bcbra{1-\delta, \brbra{(\tau_k^s)^{-1}-L_{XY}-L_{G}^2\sigma^{s}_k/\delta}\tau_k^{s}}\geq\nu_{0},\ \ \forall k\in[N_{s}].\label{eq:step_size_mid01}
\end{equation}
 Furthermore, when switching to the next epoch $(s\to s+1)$, we have
\begin{equation}
\begin{split}\sigma_{0}^{s+1}((\sigma_{0}^{s+1})^{-1}-\theta_{0}^{s+1}\delta/\sigma_{N_{s}}^{s})\overset{(a)}{\ge}\sigma_{0}^{s+1}((\sigma_{0}^{s+1})^{-1}-(\sigma_{N_{s}}^{s})^{-1})&\overset{(b)}{\ge}1-\sigma_{0}^{s+1}\underline{\sigma}^{-1} = 1-\bar{\sigma}\underline{\sigma}^{-1}\\
((\tau_{0}^{s+1})^{-1}-L_{XY}-L_{G}^{2}\delta^{-1}\sigma_{0}^{s+1})\tau_{0}^{s+1}&\overset{(c)}{\ge}\nu_{0}\tau_{0}^{s+1}=\nu_{0}\bar{\tau},
\end{split}
\label{eq:step-size_mid_02}
\end{equation}
where $(a)$ holds by $\theta_{0}^{s}=1$, $\delta<1$, $(b)$ follows
from $(\sigma_{N_{s}}^{s})^{-1}\geq\underline{\sigma}^{-1}$. Hence,
combining~\eqref{mid_01}, ~\eqref{eq:step_size_mid01} and~\eqref{eq:step-size_mid_02},
we completes our proof of~\eqref{eq:nu_1_ineq} by setting $\nu_{1}=\min\{1-\bar{\sigma}\underline{\sigma}^{-1},\nu_{0}\bar{\tau},\nu_{0}\}$.}

Since {\resapdsce} reset the stepsize periodically and $\{t_{k}^s,(\tau_{k}^s)^{-1}\}_{k\in[N_{s}]}$
are two monotonically increasing sequences, hence 
\begin{equation}
\inf_{j\ge 0}t_{j}({\tau_{j}^{-1}}-\rho_{j})\geq t_{0}^{s}({\bar{\tau}^{-1}}-\overline{\rho})=\bar{\tau}^{-1}-\overline{\rho}.\label{eq:step_size_03}
\end{equation}

Consider $\inf_{k\in[N_{s}]}t_{k}^{s}\sigma_{k}^{s}(1-\delta{\sigma_{k}^{s}}/{\sigma_{k-1}^{s}})$.
Combining $\delta+\nu_{0}\leq\inf_{k\in [N_s]}\{\theta_{k}^{s}\}$, then 
\begin{equation}
\inf_{k\in [N_s]}t_{k}^{s}\sigma_{k}^s(1-\delta\frac{\sigma_{k}^{s}}{\sigma_{k-1}^{s}})=\inf_{k\in[N_s]}t_{k}^{s}\sigma_{k}^{s}(1-\delta/\theta_{k}^{s})\geq\nu_{0}\bar{\sigma}.\label{mid_05}
\end{equation}
 Furthermore, when switching to the next epoch $(s\to s+1)$, we have
\begin{equation}
\inf_{s\geq0}t_{0}^{s+1}\sigma_{0}^{s+1}(1-\delta\sigma_{0}^{s+1}(\sigma_{N_{s}}^{s})^{-1})=\bar{\sigma}^{2}\inf_{s\geq0}(\bar{\sigma}^{-1}-\delta(\sigma_{N_{s}}^{s})^{-1})\geq\bar{\sigma}(1-\delta)
,\label{eq:step_size_04}
\end{equation}
 where the last inequality holds by $\bar{\sigma}=\sigma_0^s\leq \sigma_{N_s}^s$.
 % $(\sigma_{0}^{s+1})^{-1}\geq(\underline{\sigma})^{-1}\geq(\sigma_{N_{s}}^{s})^{-1}$.
Hence, it follows from~\eqref{eq:step_size_03},~\eqref{mid_05} and~\eqref{eq:step_size_04}
that there exist $\nu_{2}=\min\{\bar{\tau}^{-1}-\overline{\rho},\nu_{0}\bar{\sigma},\bar{\sigma}(1-\delta)\}$
% \frac{1-1}{\underline{\sigma}}
such~\eqref{eq:nu_2_ineq} holds.

Part 2. for any $j\ge 0,$ we have
\begin{equation}
\begin{split}t_{j+1}Q_{j+1}(\xbf^{*},\tilde{\ybf}^{*}) & \leq   t_{j}\big(\frac{(\tau_{j})^{-1}}{2}\|\xbf^{*}-\xbf_{j+1}\|^{2}+\left\langle G(\xbf_{j+1})-G(\xbf_{j}),\ybf_{j+1}-\tilde{\ybf}^{*}\right\rangle\\
 &\ \  +{(2\sigma_{j})^{-1}}\|\tilde{\ybf}^{*}-\ybf_{j+1}\|^{2} +\frac{\sigma_{j}}{2}\|G(\xbf_{j+1})-G(\xbf_{j})\|^{2}\big).
\end{split}
\label{eq:mid_res_lem-03}
\end{equation}

Consider $k\in\{0,1,\ldots, N_{s}\}$. Inequality~\eqref{eq:appendix_ini_step_assu} implies~\eqref{eq:param-03} holds
% , i.e., $t_{k+1}^{s}((\tau_{k+1}^s)^{-1} - \rho_{k+1}^{s})\le t_{k}^{s}(\tau_{k}^{s})^{-1}$, $\theta_{k}^{s} 1/\sigma^{s}_{k-1}\leq(\sigma^{s}_{k})^{-1}$, $t_{k+1}^{s}(\sigma_{k+1}^{s})^{-1}\leq t_{k}(\sigma_{k}^{s})^{-1}, t_{k+1}^{s}\theta_{k+1}^{s}=t_{k}^{s}$ 
(see proof
of Corollary~\ref{cor:mainthm} in Section~\ref{proof:cor}). Hence, for  {$0\le k \le N_{s},$} we have
% \begin{equation}
% \begin{split} & t_{k}^{s}(\frac{(\tau_{k}^{s})^{-1}}{2}\|\xbf^{*}-\xbf_{k+1}^{s}\|^{2}+\frac{1}{2\sigma_{k}^{s}}\|\tilde{\ybf}^{*}-\ybf_{k+1}^{s}\|^{2}\\
%  & +\left\langle G(\xbf_{k+1}^{s})-G(\xbf_{k}^{s}),\ybf_{k+1}^{s}-\tilde{\ybf}^{*}\right\rangle +\frac{1}{2/\sigma_{k}^{s}}\|G(\xbf_{k+1}^{s})-G(\xbf_{k}^{s})\|^{2})\\
% \ge & t_{k+1}^{s}\big[\frac{(\tau_{k+1}^{s})^{-1}-\rho_{k+1}^{s}}{2}\norm{\xbf^{*}-\xbf_{k+1}^{s}}^{2}+\frac{1}{2\sigma_{k+1}^{s}}\norm{\tilde{\ybf}^{*}-\ybf_{k+1}^{s}}^{2}\\
%  & +\theta_{k+1}^{s}\left\langle \ybf_{k+1}^{s}-\tilde{\ybf}^{*},G(\xbf_{k+1}^{s})-G(\xbf_{k}^{s})\right\rangle +\frac{\theta_{k+1}^{s}}{2/\sigma_{k}^{s}}\norm{G(\xbf_{k+1}^{s})-G(\xbf_{k}^{s})}^{2}\big]\\
% = & t_{j+1}Q_{j+1}(\xbf^{*},\tilde{\ybf}^{*}),
% \end{split}
% \label{eq:mid_res_lem-01}
% \end{equation}
\begin{equation}
\begin{split} t_{j+1}Q_{j+1}(\xbf^{*},\tilde{\ybf}^{*}) &\leq  t_{k}^{s}\big(\frac{(\tau_{k}^{s})^{-1}}{2}\|\xbf^{*}-\xbf_{k+1}^{s}\|^{2} +\left\langle G(\xbf_{k+1}^{s})-G(\xbf_{k}^{s}),\ybf_{k+1}^{s}-\tilde{\ybf}^{*}\right\rangle\\
 &\ \  +\frac{1}{2\sigma_{k}^{s}}\|\tilde{\ybf}^{*}-\ybf_{k+1}^{s}\|^{2} +\frac{\sigma_{k}^{s}}{2\delta}\|G(\xbf_{k+1}^{s})-G(\xbf_{k}^{s})\|^{2}\big)
% \ge & t_{k+1}^{s}\big[\frac{(\tau_{k+1}^{s})^{-1}-\rho_{k+1}^{s}}{2}\norm{\xbf^{*}-\xbf_{k+1}^{s}}^{2}+\frac{1}{2\sigma_{k+1}^{s}}\norm{\tilde{\ybf}^{*}-\ybf_{k+1}^{s}}^{2}\\
%  & +\theta_{k+1}^{s}\left\langle \ybf_{k+1}^{s}-\tilde{\ybf}^{*},G(\xbf_{k+1}^{s})-G(\xbf_{k}^{s})\right\rangle +\frac{\theta_{k+1}^{s}}{2/\sigma_{k}^{s}}\norm{G(\xbf_{k+1}^{s})-G(\xbf_{k}^{s})}^{2}\big]\\
% = & ,
\end{split}
\label{eq:mid_res_lem-01}
\end{equation}
where $j$ corresponds to $(s,k)$. Furthermore, consider switching
to next epoch $(s\to s+1)$. 
 {Since $t_{k}^{s}(\tau_{k}^{s})^{-1}$ is an
increasing sequence in $k$, $\rho_{0}^{s+1}>0,t_0^{s+1}=1$, hence}
\begin{equation}
t_{N_{s}}^{s}(\tau_{N_{s}}^{s})^{-1}\geq t_{0}^{s+1}(\tau_{0}^{s+1})^{-1}-\rho_{0}^{s+1}t_{0}^{s+1},\forall s\geq 0.\label{eq:hats_def2}
\end{equation}
Next, we have
\begin{equation}\label{eq:switch-01}
    \frac{t_{N_{s}}^{s}}{\sigma_{N_{s}}^{s}} \overset{(a)}{=}\frac{t_{0}^{s+1}}{\sigma_{0}^{s+1}},\ t_{N_{s}}^{s}\overset{(b)}{\ge}t_{0}^{s+1}\overset{(c)}{=}t_{0}^{s+1}\theta_{0}^{s+1}, t_{N_{s}}^{s}\sigma_{N_{s}}^{s}  \overset{(b)}{\ge}t_{0}^{s+1}\sigma_{0}^{s+1}\overset{(c)}{=}t_{0}^{s+1}\sigma_{0}^{s+1}\theta_{0}^{s+1},
\end{equation}
where $(a)$ holds by the definition of $t_{k}^{s}=\frac{\sigma_{k}^{s}}{\sigma_{0}^{s}}$,
$(b)$ holds by $\{t_{k}^{s},\sigma_{k}^{s}\}$ is an increasing sequence
in $k$, and $(c)$ holds by $\theta_{0}^{s+1}=1$. Hence, by~\eqref{eq:hats_def2} and~\eqref{eq:switch-01}, we have
% \begin{equation}
% \begin{split} & t_{N_{s}}^{s}(\frac{(\tau_{N_{s}}^{s})^{-1}}{2}\|\xbf^{*}-\xbf_{0}^{s+1}\|^{2}+\frac{1}{2\sigma_{N_{s}}^{s}}\|\tilde{\ybf}^{*}-\ybf_{0}^{s+1}\|^{2}\\
%  & +\left\langle G(\xbf_{0}^{s+1})-G(\xbf_{N_{s}}^{s}),\ybf_{0}^{s+1}-\tilde{\ybf}^{*}\right\rangle +\frac{1}{2/\sigma_{N_{s}}^{s}}\|G(\xbf_{0}^{s+1})-G(\xbf_{N_{s}}^{s})\|^{2})\\
% \geq & t_{0}^{s+1}\big(\frac{(\tau_{0}^{s+1})^{-1}-\rho_{0}^{s+1}}{2}\norm{\xbf^{*}-\xbf_{0}^{s+1}}^{2}+\frac{1}{2\sigma_{0}^{s+1}}\norm{\ybf^{*}-\ybf_{0}^{s+1}}^{2}\\
%  & +\theta_{0}^{s+1}\left\langle \ybf_{0}^{s+1}-\tilde{\ybf}^{*},G(\xbf_{0}^{s+1})-G(\xbf_{N_{s}}^{s})\right\rangle +\frac{\theta_{0}^{s+1}}{2/\sigma_{0}^{s+1}}\norm{G(\xbf_{0}^{s+1})-G(\xbf_{N_{s}}^{s})}^{2}\big)\\
% = & t_{j+1}Q_{j+1}(\xbf^{*},\tilde{\ybf}^{*}),
% \end{split}
% \label{eq:mid_res_lem_02}
% \end{equation}
\begin{equation}
\begin{split}
t_{j+1}Q_{j+1}(\xbf^{*},\tilde{\ybf}^{*})&\leq   t_{N_{s}}^{s}\big(\frac{1}{2\tau_{N_{s}}^{s}}\|\xbf^{*}-\xbf_{0}^{s+1}\|^{2}+\frac{\sigma_{N_{s}}^{s}}{2}\|G(\xbf_{0}^{s+1})-G(\xbf_{N_{s}}^{s})\|^{2}\\
  &\ \ +\frac{1}{2\sigma_{N_{s}}^{s}}\|\tilde{\ybf}^{*}-\ybf_{0}^{s+1}\|^{2}+\left\langle G(\xbf_{0}^{s+1})-G(\xbf_{N_{s}}^{s}),\ybf_{0}^{s+1}-\tilde{\ybf}^{*}\right\rangle \big)
\end{split}
\label{eq:mid_res_lem_02}
\end{equation}
where $j$ corresponds to $(s,N_{s})$.
% Then we complete the proof
% of~\eqref{eq:mid_res_lem-03} by putting~\eqref{eq:mid_res_lem-01}
% and~\eqref{eq:mid_res_lem_02} together.
By putting~\eqref{eq:mid_res_lem-01}
and~\eqref{eq:mid_res_lem_02} together, we complete the proof of~\eqref{eq:mid_res_lem-03}.

Placing $(\xbf,\ybf)=(\xbf^{*},\tilde{\ybf}^{*}),(\xbf_{k+1},\ybf_{k+1})=(\xbf_{j+1},\ybf_{j+1})$
in~\eqref{eq:bounded_lagran} and multiplying $t_{j}$ on both sides, we have
\begin{equation}
\begin{split}0 & \le t_{j}[\mcal L(\xbf_{j+1},\tilde{\ybf}^{*})-\mcal L({\bf \xbf}^{*},\ybf_{j+1})]\\
 & \le t_{j}\big[\frac{\tau_{j}^{-1}-\rho_{j}}{2}\|\xbf-\xbf_{j}\|^{2}-\frac{\tau_{j}^{-1}}{2}\|\xbf-\xbf_{j+1}\|^{2}+\frac{\theta_j}{2\delta/\sigma_{j-1}}\|\qbf_{j}\|^{2}-\frac{1}{2\delta/\sigma_{j}}\|\qbf_{j+1}\|^{2}\\
 & \quad+{(2\sigma_{j})^{-1}}\left(\|\ybf-\ybf{}_{j}\|^{2}-\|\ybf-\ybf_{j+1}\|^{2}\right)+\langle\ybf-\ybf_{j+1},\qbf_{j+1}\rangle-\theta_j\langle\ybf-\ybf_{j},\qbf_{j}\rangle\\
 & \quad-\frac{\sigma_{j}^{-1}-\theta_j \delta/\sigma_{j-1}}{2}\|\ybf_{j+1}-\ybf_{j}\|^{2}+\frac{L_{G}^{2}}{2\delta/\sigma_{j}}\|\xbf_{j+1}-\xbf_{j}\|^{2}-\frac{\tau_{j}^{-1}-L_{XY}}{2}\|\xbf_{j+1}-\xbf_{j}\|^{2}\big]\\
 & \le t_{j}Q_{j}(\xbf^{*},\tilde{\ybf}^{*})-t_{j+1}Q_{j+1}(\xbf^{*},\tilde{\ybf}^{*})-\nu_{1}t_{j}[{(2\tau_{j})^{-1}}\norm{\xbf_{j+1}-\xbf_{j}}^{2}+{(2\sigma_{j})^{-1}}\norm{\ybf_{j+1}-\ybf_{j}}^{2}],
\end{split}
\label{eq:res_01}
\end{equation}
where the last inequality holds by~\eqref{eq:mid_res_lem-03} and~\eqref{eq:nu_1_ineq}.
It follows from~\eqref{eq:nu_2_ineq}, $\sigma_{j-1}/\sigma_{j}\le 1$ and $$\left\langle \ybf_{j}-\tilde{\ybf}^{*},\qbf_{j}\right\rangle\geq-\frac{\sigma_{j-1}}{2\delta}\norm{\qbf_{j}}^{2}-\frac{\delta/\sigma_{j-1}}{2}\norm{\tilde{\ybf}^{*}-\ybf_{k}}^{2}$$
that 
\begin{equation}
\begin{split}t_{j}Q_{j}(\xbf^{*},\tilde{\ybf}^{*}) 
% & \ge t_{j}\left({(2\tau_{j})^{-1}}\|\xbf^{*}-\xbf_{j}\|^{2}+{(2\sigma_{j})^{-1}}\|\tilde{\ybf}^{*}-\ybf_{j}\|^{2}-\frac{(\sigma_{j-1}/\sigma_{j})1}{2\sigma_{j-1}}\norm{\ybf_{j}-\tilde{\ybf}^{*}}^{2}\right)\\
 & \geq t_{j}\brbra{{(2\tau_{j})^{-1}}\|\xbf^{*}-\xbf_{j}\|^{2}+{(2\sigma_{j})^{-1}}\|\tilde{\ybf}^{*}-\ybf_{j}\|^{2}-\frac{\delta}{2\sigma_{j-1}}\norm{\ybf_{j}-\tilde{\ybf}^{*}}^{2}}\\
 & \geq\nu_{2}\brbra{\frac{1}{2}\norm{\xbf^{*}-\xbf_{j}}^{2}+\frac{1}{2}\norm{\ybf_{j}-\tilde{\ybf}^{*}}^{2}}>0.
\end{split}
\label{eq:res_02}
\end{equation}
Combining~\eqref{eq:res_01}
and~\eqref{eq:res_02}, we complete our proof of~\eqref{eq:res_lem}.
\endproof

\subsection{Proof of Theorem~\ref{thm:asymptotic}\label{proof_of_thm_asymptotic}}
\proof
Since  $\bcbra{(\xbf_{j},\ybf_{j})}$ located in set $\mcal X\times\mcal Y$
is a bounded sequence, it must have a convergent subsequence $\lim_{n\to\infty}(\xbf_{j_{n}},\ybf_{j_{n}})=~(\xbf^{*},\ybf^{*})$,
where $\ybf^{*}$ is the limit point. We claim  that limit point $(\xbf^{*},\ybf^{*})$ satisfies the
KKT condition. Placing $a_{j}=t_{j}Q_{j}(\xbf^{*},\tilde{\ybf}^{*})$,
$b_{j}=\nu_{1}t_{j}[{(2\tau_{j})^{-1}}\norm{\xbf_{j+1}-\xbf_{j}}^{2}+{(2\sigma_{j})^{-1}}\norm{\ybf_{j+1}-\ybf_{j}}^{2}]$
and $c_{j}=0$ in Lemma~\ref{lem:limit_converge}.
It follows from~\eqref{eq:res_lem} in Lemma~\ref{lem:asy} that
$a_{j}\geq0,b_{j}>0$. Hence, we have $\sum_{j=0}^{\infty}\|\xbf_{j+1}-\xbf_{j}\|^{2}<\infty,$
and $\sum_{j=0}^{\infty}\|\ybf_{j+1}-\ybf_{j}\|^{2}<\infty$,
which implies $\lim_{n\to\infty}\|\xbf_{j_{n}}-\xbf_{j_{n}+1}\|^{2}=0$
and $\lim_{n\to\infty}\|\ybf_{j_{n}}-\ybf_{j_{n}+1}\|^{2}=0$. 
 {There are two different cases for $\tau_{j_n}$ when $j_n \to \infty$, and we discuss the value of $B_{j_n+1}$ in~\eqref{eq:finequality} decided by $\tau_{j_n}$ in each of the two cases below.}

 Case 1: $\tau_{j_n}^{-1}<\infty$. By the definition of $B_{j_n + 1}$ in~\eqref{eq:defBk+1} and $\lim_{n\to\infty}\|\xbf_{j_{n}}-\xbf_{j_{n}+1}\|^{2}=0$, we have $  B_{j_n+1}   \le \norm{\xbf-\xbf_{j_n+1}} \cdot \norm{\xbf_{{j_n}+1} - \xbf_{j_n}} / \tau_{j_n} \overset{n \to \infty}{\longrightarrow} 0.$
% \begin{equation*}
% \begin{aligned}
%     B_{j_n+1} 
%     % &= \frac{1}{2\tau_{j_n}}\brbra{\norm{\xbf - \xbf_{j_n}}^2 - \norm{\xbf - \xbf_{j_n+1}}^2 - \norm{\xbf_{j_n+1} - \xbf_{j_n}}^2} \\
%     % &
%     \le \frac{1}{\tau_{j_n}}\brbra{\norm{\xbf-\xbf_{j_n+1}} \norm{\xbf_{{j_n}+1} - \xbf_{j_n}} } \overset{n \to \infty}{\longrightarrow} 0.
% \end{aligned}
% \end{equation*}

 {Case 2: $\tau_{j_n}^{-1}=\infty$. It follows from~\eqref{eq:gammak_upper} that $\tau_{j_n}^{-1}$ increases at order $\Theta(k)$, where $j_n\sim (s,k)$.  By~\eqref{eq:thm2}, we obtain $\norm{\xbf-\xbf_{j_n}}$ decreases at order $\mcal{O}(1/k)$ ($j_n\sim (s,k)$). Hence, combining $\lim_{n\to\infty}\|\xbf_{j_{n}}-\xbf_{j_{n}+1}\|^{2}=0$,  we have $ B_{j_n+1} \le \frac{1}{\tau_{j_n}}\brbra{\norm{\xbf-\xbf_{j_n+1}} \norm{\xbf_{{j_n}+1} - \xbf_{j_n}} } \overset{n \to \infty}{\longrightarrow} 0$.
% \begin{equation*}
% \begin{aligned}
%         B_{j_n+1} &\le \frac{1}{\tau_{j_n}}\brbra{\norm{\xbf-\xbf_{j_n+1}} \norm{\xbf_{{j_n}+1} - \xbf_{j_n}} } \overset{n \to \infty}{\longrightarrow} 0.
% \end{aligned}
% \end{equation*}
}
 {It follows from $\lim_{n\to\infty}\xbf_{j_{n}} = \xbf^{*}$, $\lim_{n\to\infty} B_{j_n+1}=0$ and~\eqref{eq:finequality}
that 
$$ f(\xbf^{*})+\langle\nabla G(\xbf^{*})\ybf^{*},\xbf^{*}\rangle\le f(\xbf)+\langle\nabla G(\xbf^{*})\ybf^{*},\xbf\rangle, \forall \xbf\in\Xcal.$$}
Hence, according to the first-order optimality condition, we have
\begin{equation}
\boldsymbol{0}\in\partial f(\xbf^{*})+\nabla G(\xbf^{*})\ybf^{*}+\mcal N_{\mcal X}(\xbf^{*}).\label{eq:mid-opt-01_new}
\end{equation}
Next, we show the complementary slackness holds for $(\xbf^*,\ybf^*)$.  {Since $\sigma_{j_n}^{-1} $ has an upper bound $\bar{\sigma}^{-1}$, $\norm{\ybf-\ybf_{j_n+1}}\le D_Y$, $\lim_{n\to\infty}\|\ybf_{j_{n}}-\ybf_{j_{n}+1}\|^{2}=0$ and the definition of $A_{j_n+1}$ in~\eqref{eq:defAk+1}, hence we obtain $   A_{j_n+1}   \leq \frac{1}{\sigma_{j_n}}\brbra{\norm{\ybf_{j_n} - \ybf_{j_n+1}}\norm{\ybf - \ybf_{j_n+1}}}\overset{n\to\infty}{\longrightarrow}0.$
% \begin{equation*}
% \begin{aligned}
%         A_{j_n+1} 
%         % &= \frac{1}{2\sigma_{j_n}}\left(\|\ybf-\ybf{}_{j_n}\|^{2}-\|\ybf-\ybf_{j_n+1}\|^{2}-\|\ybf_{j_n+1}-\ybf_{j_n}\|^{2}\right)\\
%     % &
%     \leq \frac{1}{\sigma_{j_n}}\brbra{\norm{\ybf_{j_n} - \ybf_{j_n+1}}\norm{\ybf - \ybf_{j_n+1}}}\overset{n\to\infty}{\longrightarrow}0.
% \end{aligned}
% \end{equation*}
}
 {Combining above, $\lim_{n\to \infty} \ybf_{j_n}=\ybf^*$ and~\eqref{eq:indicator}}, we
have $0\le-\langle G(\xbf^{*}),\ybf^{*}\rangle\le-\langle G(\xbf^{*}),\ybf\rangle$, $\forall\ybf\in\Ycal$.
% \[
% 0\le-\langle G(\xbf^{*}),\ybf^{*}\rangle\le-\langle G(\xbf^{*}),\ybf\rangle,\quad\forall\ybf\in\Ycal.
% \]
Moreover, due to the complementary slackness, there exists an $\hat{\ybf}^{*}\in\Ycal^{*}\subseteq\mcal Y$
such that $-\inner{G(\xbf^{*})}{\hat{\ybf}^{*}}=0$. Hence, we must
have $\langle G(\xbf^{*}),\ybf^{*}\rangle=0$, which, together with~\eqref{eq:mid-opt-01_new},
implies that $(\xbf^{*},\ybf^{*})$ is KKT point.
\endproof

% \section{\label{sec:Proof_of_msapd}
% % Proof details in Section~\ref{sec:MSAPD}
% Convergence analysis of \MSAPD{}
% }

\section{\label{sec:Proof-of-manifold}
% Proof details in Section~\ref{sec:Manifold}
Proof details for sparsity identification
}

Our proof strategy of active-set identification in {\resapdsce} is similar to those in unconstrained optimization~\citep{nutini2019activeset}. Namely, we show that the optimal sparsity pattern is identified when the iterates fall in a properly defined neighborhood dependent on $\eta$. The next lemma shows that the primal and dual sequences indeed converge to the neighborhood of the optimal primal and dual solutions, respectively,  in a finite number of iterations.
% \begin{lemma}[Informal]\label{lem:xsmall_epoch}
% There exists an epoch $\hat{S}_{0}$ such that $\forall s\geq\hat{S}_{0},$ we have
% \begin{equation}
% \|\ybf_k^{s} - \ybf^*\|\le \frac{\eta}{3\|\nabla g(\xbf^*)\|}, \ 
% \norm{\xbf_{k}^{s}-\xbf^{*}}\leq\frac{\eta}{3L_{XY}}\frac{\tau_{k}^{s}}{\tau_{k}^{s}+(2L_{XY})^{-1}},\  \forall k=0,1,\ldots N_{s}.\label{x_small-1}
% \end{equation}
% \end{lemma}

% \subsection{Proof of Lemma~\ref{lem:xsmall_epoch}}

\begin{lemma}\label{lem:xsmall_epoch}
There exists an $\hat{S}_{1}$ such that 
\begin{equation}\label{eq:def_hats1}
\norm{\xbf_{0}^{s}-\xbf^{*}}\leq\norm{\xbf_{0}^{\hat{S}_1}-\xbf^{*}}\ \ \text{and}\ \ \norm{\ybf_{0}^{s}-\ybf^{*}}\leq\norm{\ybf_{0}^{\hat{S}_1}-\ybf^{*}},\forall s\geq\hat{S}_{1},
\end{equation}
where $(\xbf^*,\ybf^*)$ is the unique solution of problem~\eqref{eq:sparselearn}.
Moreover, there exists an epoch $\hat{S}_{0}\geq\hat{S}_{1}$ such that $\forall s\geq\hat{S}_{0},$ we have
\begin{equation}
\|\ybf_k^{s} - \ybf^*\|\le \frac{\eta}{3\|\nabla g(\xbf^*)\|}, \ 
\norm{\xbf_{k}^{s}-\xbf^{*}}\leq\frac{\eta}{3L_{XY}}\frac{\tau_{k}^{s}}{\tau_{k}^{s}+(2L_{XY})^{-1}},\  \forall k=0,1,\ldots N_{s}.\label{appendix:x_small-1}
\end{equation}
\end{lemma}

\proof
{
From Theorem~\ref{thm:rapdpro} and~\ref{thm:asymptotic}, we have $\lim_{j\to\infty}(\xbf_j,\ybf_j)=(\xbf^*,\ybf^*)$, where $j$ corresponds to $(s,0)$. It implies that there exists an epoch $\hat{S}_1$ such that~\eqref{eq:def_hats1} holds.}

It follows from~\eqref{eq:mid-03} that $\norm{\xbf_{1}^{s}-\xbf^{*}}\leq\sqrt{{\sigma_{0}^{s}\tau_{0}^{s}}/{\sigma_{1}^{s}}(\norm{\xbf_{0}^{s}-\xbf^{*}}^{2}/\tau_{0}^{s}+\norm{\ybf_{0}^{s}-\ybf^{*}}^{2}/{\sigma_{0}^{s}})}$.
Hence, in order to prove $\norm{\xbf_{1}^{s}-\xbf^{*}}\leq\frac{\eta}{3L_{XY}}\cdot \frac{\tau_{k}^{s}}{\tau_{k}^{s}+(2L_{XY})^{-1}}$,
we need to prove 
\begin{equation}
\sqrt{\frac{\sigma_{0}^{s}\tau_{0}^{s}}{\sigma_{1}^{s}}(\frac{1}{\tau_{0}^{s}}\norm{\xbf_{0}^{s}-\xbf^{*}}^{2}+\frac{1}{\sigma_{0}^{s}}\norm{\ybf_{0}^{s}-\ybf^{*}}^{2})}\leq\frac{\eta}{3L_{XY}}\frac{\tau_{1}^{s}}{\tau_{1}^{s}+(2L_{XY})^{-1}}.\label{eq:suff_cond-1}
\end{equation}
From Corollary~\ref{cor:mainthm} and Theorem~\ref{thm:rapdpro},~\ref{thm:asymptotic}, we know that the left hand side of~\eqref{eq:suff_cond-1}
converges to $0$ and right hand side of~\eqref{eq:suff_cond-1}
is a positive constant. Hence, there exist a $\hat{S}_{2}$ such that~\eqref{eq:suff_cond-1} holds,
which implies~\eqref{appendix:x_small-1} holds for $k=1,s=\hat{S}_2$. Now we use induction
to prove, for $\forall k\in[N_{\hat{S}_{2}}]$, we have
\begin{equation}
\brbra{\frac{\sigma_{0}^{\hat{S}_{2}}\tau_{k}^{\hat{S}_{2}}}{\sigma_{k}^{\hat{S}_{2}}}(\frac{1}{\tau_{0}^{\hat{S}_{2}}}\norm{\xbf_{0}^{\hat{S}_{2}}-\xbf^{*}}^{2}+\frac{1}{\sigma_{0}^{\hat{S}_{2}}}\norm{\ybf_{0}^{\hat{S}_{2}}-\ybf^{*}}^{2})}^{1/2}\leq\frac{\eta}{3L_{XY}}\frac{\tau_{k}^{\hat{S}_{2}}}{\tau_{k}^{\hat{S}_{2}}+(2L_{XY})^{-1}}.\label{x_small_k-1}
\end{equation}
When $k=1$, inequality~\eqref{x_small_k-1} coincides with~\eqref{eq:suff_cond-1}
with $s=\hat{S}_{2}$. Now, assume~\eqref{x_small_k-1} holds for
$k$, we aim to prove~\eqref{x_small_k-1} holds for $k+1$. It follows from~\eqref{eq:mid-03} that
\begin{align*}
 \norm{\xbf_{k+1}^{\hat{S}_{2}}-\xbf^{*}}
 % & le\sqrt{\frac{\sigma_{0}^{\hat{S}_{2}}\tau_{k+1}^{\hat{S}_{2}}}{\sigma_{k+1}^{\hat{S}_{2}}}(\frac{1}{\tau_{0}^{\hat{S}_{2}}}\norm{\xbf_{0}^{\hat{S}_{2}}-\xbf^{*}}^{2}+\frac{1}{\sigma_{0}^{\hat{S}_{2}}}\norm{\ybf_{0}^{\hat{S}_{2}}-\ybf^{*}}^{2})}\\
 & \overset{(a)}{\leq}\sqrt{\frac{\tau_{k+1}^{\hat{S}_{2}}}{\sigma_{k+1}^{\hat{S}_{2}}}\cdot\frac{\sigma_{k}^{\hat{S}_{2}}}{\tau_{k}^{\hat{S}_{2}}}}\cdot \frac{\eta}{3L_{XY}}\cdot \frac{\tau_{k}^{\hat{S}_{2}}}{\tau_{k}^{\hat{S}_{2}}+(2L_{XY})^{-1}}
 % & =\sqrt{\frac{\sigma_{k}^{\hat{S}_{2}}}{\sigma_{k+1}^{\hat{S}_{2}}}\cdot\frac{\tau_{k}^{\hat{S}_{2}}}{\tau_{k+1}^{\hat{S}_{2}}}}\cdot \frac{\eta}{3L_{XY}}\cdot \frac{\tau_{k+1}^{\hat{S}_{2}}}{\tau_{k}^{\hat{S}_{2}}+(2L_{XY})^{-1}}\\
 % &
 \overset{(b)}{=}\frac{\eta}{3L_{XY}}\cdot \frac{\tau_{k+1}^{\hat{S}_{2}}}{\tau_{k}^{\hat{S}_{2}}+(2L_{XY})^{-1}}\\
 & \overset{(c)}{\le}\frac{\eta}{3L_{XY}}\cdot \frac{\tau_{k+1}^{\hat{S}_{2}}}{\tau_{k+1}^{\hat{S}_{2}}+(2L_{XY})^{-1}},
\end{align*}
where $(a)$ follows from induction, $(b)$ holds by $\tau_{k}^{\hat{S}_{2}}\sigma_{k}^{\hat{S}_{2}}=\tau_{k+1}^{\hat{S}_{2}}\sigma_{k+1}^{\hat{S}_{2}}$ and $(c)$ holds by $\tau_{k+1}^{\hat{S}_{2}}\leq\tau_{k}^{\hat{S}_{2}}$.
Hence, we complete our proof of~\eqref{x_small_k-1}.
From Theorem~\ref{thm:rapdpro}, we have $\norm{\xbf_{0}^{s}-\xbf^{*}}^2\leq D_{X}^2 \cdot 2^{-s}$, which implies that  there exists a $\hat{S}_{3}=\big\lceil 2 \log_{2}\big\{D_X (\frac{\eta}{3L_{XY}}\cdot\frac{\bar{\tau}}{\bar{\tau}+(2L_{XY})^{-1}})^{-1}\big\}\big\rceil$
such that $\norm{\xbf_{0}^{\hat{S}_{3}}-\xbf^{*}}\leq D_{X}\cdot \sqrt{2}^{-\hat{S}_{3}}\leq\frac{\eta}{3L_{XY}}\frac{\bar{\tau}}{\bar{\tau}+(2L_{XY})^{-1}}$,
which implies that $\norm{\xbf_{0}^{s}-\xbf^{*}}\leq D_{X}^{2}\cdot 2^{-s}\leq\frac{\eta}{3L_{XY}}\frac{\bar{\tau}}{\bar{\tau}+(2L_{XY})^{-1}}$
holds for any $s\geq\hat{S}_{3}$.

It follows from the definition
of $\hat{S}_{1}$ in~\eqref{eq:def_hats1} and stepsize will be reset at different epoch, then we have~\eqref{eq:suff_cond-1} holds for
$s\geq\max\{\hat{S}_{1},\hat{S}_{2}\}$, which implies that~\eqref{x_small_k-1}
holds with substituting $\hat{S}_{2}$ as any $s\geq\max\{\hat{S}_{1},\hat{S}_{2}\}$. 
{Furthermore, it follows from Theorem~\ref{thm:asymptotic} that $\lim_{j\to\infty}\ybf_j=\ybf^*$, where $j$ corresponds to $(s,k)$. Then there exists a $\hat{S}_4$ such that the first term in~\eqref{appendix:x_small-1} holds.}
Hence,
we can obtain that there exist a $\hat{S}_{0}=\max\{\hat{S}_{1},\hat{S}_{2},\hat{S}_{3},{\hat{S}_4}\}$
such that~\eqref{appendix:x_small-1} holds.
\endproof

It is worth noting that the primal neighborhood defined by the second term of~\eqref{appendix:x_small-1} is a bit different from the fixed neighborhood in the standard analysis~\citep{nutini2019activeset},  which involves a constant stepsize. 
As \APDSCE{} sets $\tau_k^s=\Ocal({1}/{k})$, both the point distance and neighborhood radius decay at the same $\Ocal({1}/{k})$ rate. Hence, we use a substantially different analysis to show the sparsity identification in the constrained setting.

\subsection{Proof of Proposition~\ref{prop:uniquey}\label{proof_prop_uniquey}}
% \begin{reproposition}{\ref{prop:uniquey}}
% \label{appendix:prop:uniquey}
% Under Assumptions~\ref{assu:Slater's} and~\ref{assu:dual-optim}, both the optimal solution $\xbf^*$ and the Lagrangian multiplier $\ybf^*$ satisfying the KKT condition for problem~\eqref{eq:sparselearn} are unique. 
% \end{reproposition}
\proof
The uniqueness of primal
optimal solution $\xbf^{*}$ follows from Proposition~\ref{prop:uniquex}.
The KKT condition (ensured by Slater's CQ) implies 
\begin{equation}\label{eq:kkt-2}
    \zerobf\in \partial f(\xbf^*) + \nabla g(\xbf^*) \ybf^*.
\end{equation}
According to Assumption~\ref{assu:dual-optim}, we have $\xbf^* \neq \boldsymbol{0}$, hence $\Acal^c(\xbf^*)=\{1,2,\ldots,B\}\setminus\Acal(\xbf^*)\neq\emptyset$. In view of \eqref{eq:kkt-2}, for any $i\in\Acal^c(\xbf)$, we have 
$p_i {\xbf^*_{(i)}}/{\norm{\xbf^*_{(i)}}}=-\nabla_{(i)} g(\xbf^*)\ybf^*$, which gives a unique $\ybf^*$.
\endproof

\subsection{Proof of Theorem~\ref{thm:active_set_ide}\label{sec:proof_thm_active_set_ide}}

\proof
It follows from the Lipschitz smoothness of $g(\cdot)$ and property~\eqref{appendix:x_small-1} that for any $s\ge \hat{S}_0$, we have 
\begin{equation}
\begin{aligned} \label{eq:mid-02}
& \bnorm{\left[\nabla g(\xbf_{k}^{s})\ybf_{k+1}^{s}\right]_{(i)}}-\bnorm{\left[\nabla g(\xbf^{*})\ybf_{k+1}^{s}\right]_{(i)}}\\
 & \leq \bnorm{\nabla g(\xbf_{k}^{s})\ybf_{k+1}^{s}-\nabla g(\xbf^{*})\ybf_{k+1}^{s}}\\
 & \leq L_{XY}\bnorm{\xbf_{k}^{s}-\xbf^{*}}
\leq \frac{\eta}{3}\frac{\tau_{k}^{s}}{\tau_{k}^{s}+(2L_{XY})^{-1}},\  k=0,\ldots N_s.
\end{aligned}
\end{equation}
Recall that the primal update has the following form
\[
\xbf_{k+1}^{s}=\argmin_{\xbf\in \mathcal{X}}\Bcbra{\sum_{i=1}^Bp_{i}\|\xbf_{(i)}\|+\left\langle \nabla g(\xbf_{k}^{s})\text{\ensuremath{\ybf}}_{k+1}^{s},\xbf\right\rangle +\frac{1}{2\tau_{k}^{s}}\|\xbf-\xbf_{k}^{s}\|^{2}}.
\]
Since ${\tau_{k}^s}/\rbra{\tau_{k}^s + (2L_{XY})^{-1}}$ is monotonically increasing with respect to $\tau_{k}^{s}$, {for the strictly feasible point $\tilde{\xbf}$}, we have 
\begin{equation}\label{identi-01}
\begin{aligned}
\norm{\xbf_{k+1}^s - \tilde{\xbf}}
% \leq \norm{\xbf_{k+1}^s-\xbf^*}+\norm{\xbf^* - \tilde{\xbf}}
 & \overset{(a)}{\leq} \frac{\eta}{3L_{XY}} \cdot \frac{\bar{\tau}}{\bar{\tau} + (2L_{XY})^{-1}} + \norm{\xbf^*-\tilde{\xbf}} \\
& \overset{(b)}{<}\zeta + \min\nolimits_{i\in[m]}2\sqrt{\frac{-2g_{i}(\xbf_{i}^{*})}{\mu_{i}}},    
\end{aligned}
\end{equation}
where {$(a)$ holds by~\eqref{appendix:x_small-1}, $\bar{\tau}\geq \tau_k^s$ and $(b)$ follows from the definition of $\xbf^*,\tilde{\xbf}$ and $\zeta$}. Inequality~\eqref{identi-01} implies that $\xbf_{k+1}^{s}\in \intr \mcal{X}$, and hence $\mcal{N}_{\mathcal{X}}(\xbf_{k+1}^s) = \left\{\boldsymbol{0}\right\}$. In view of the optimality condition, we have
\begin{equation}\label{eq:optimality-02}
    \Bsbra{\frac{1}{\tau_k^s}(\xbf_k^s-\xbf_{k+1}^{s}) - \nabla g(\xbf_k^s)\ybf_{k+1}^s}_{(i)} \in p_i\partial \|[\xbf_{k+1}^{s}]_{(i)}\|,\ 1\le i\le B.
\end{equation}

Our next goal is to show $[\xbf^S_{k+1}]_{(i)}=\xbf^*_{(i)}$ satisfies condition~\eqref{eq:optimality-02} for $i\in \mcal{A}(\xbf^*)$.
Placing $\xbf_{(i)}=\xbf^{*}_{(i)}$ in $ \bnorm{\big[\nabla g(\xbf_{k}^{S})\ybf_{k+1}^{S}+\frac{1}{\tau_{k}^{s}}(\xbf -\xbf_{k}^{s}) \big]_{(i)}}$, we have
\begin{equation}
\label{eq:x_star_activeset}
    \begin{split} & \bnorm{\big[\nabla g(\xbf_{k}^{s})\ybf_{k+1}^{s}+\frac{1}{\tau_{k}^{s}}(\xbf^{*}-\xbf_{k}^{s})\big]_{(i)}}\\
\leq{} & \bnorm{\left[\nabla g(\xbf_{k}^{s})\ybf_{k+1}^{s}\right]_{(i)}\|+\|\frac{1}{\tau_{k}^{s}}(\xbf_{(i)}^{*}-\xbf_{k(i)}^{s})}\\
\overset{(a)}{\le}{} & \frac{\eta}{3}\frac{\tau_{k}^{s}}{\tau_{k}^{s}+(2L_{XY})^{-1}} + \bnorm{\left[\nabla g(\xbf^{*})\ybf_{k+1}^{s}\right]_{(i)}}+\frac{\eta}{3}\frac{(L_{XY})^{-1}}{\tau_{k}^{s}+(2L_{XY})^{-1}}\\
\overset{(b)}{\leq}{} & \frac{\eta}{3}\bsbra{\frac{\tau_{k}^{s}+2(2L_{XY})^{-1}}{\tau_{k}^{s}+(2L_{XY})^{-1}} + 1 
}
+\bnorm{\left[\nabla g(\xbf^{*})\ybf^{*}\right]_{(i)}}\\
<{} & \eta+\bnorm{\left[\nabla g(\xbf^{*})\ybf^{*}\right]_{(i)}}
\overset{(c)}{\le}{}  p_{i}, \forall i\in \mcal{A}(\xbf^*).
\end{split}
\end{equation}
In above, $(a)$ follows from~\eqref{appendix:x_small-1} and~\eqref{eq:mid-02}, $(b)$ follows from 
\[\|[\nabla g(\xbf^{*})\ybf_{k+1}^{S}]_{(i)}\|-\|[\nabla g(\xbf^{*})\ybf^{*}]_{(i)}\| \le \|\ybf_{k+1}^{S} -\ybf^{*}\|\|\nabla g(\xbf^{*})\|\leq\frac{\eta}{3},
\] and $(c)$ holds by the definition of $\eta$. Combining~\eqref{eq:optimality-02} and~\eqref{eq:x_star_activeset},  we have $\mcal{A}(\xbf^*)\subseteq \mcal{A}(\xbf_{k+1}^{s}),  s \ge  \hat{S}_0, \forall k \in [N_{s}]$, which completes our proof.
\endproof
\begin{table}[htbp]
  \centering
\caption{\label{tab:dataset_description}Datasets description and parameter settings}
    \begin{tabular}{lrrrr}
    \toprule
    dataset & \multicolumn{1}{l}{Node(n)} & Edge  & \multicolumn{1}{c}{$b$} & \multicolumn{1}{c}{$\alpha$} \\
    \midrule
    bio-CE-HT & 2617  & 3K    & -0.04 & 0.4 \\
    bio-CE-LC & 1387  & 2K    & -0.05 & 0.4 \\
    econ-beaflw & 502   & 53K   & -0.01 & 0.995 \\
    DD68  & 775   & 2K    & -0.005 & 0.4 \\
    DD242 & 1284  & 3K    & -0.05 & 0.4 \\
    peking-1 & 3341  & 13.2K & -0.001 & 0.4 \\
    \bottomrule
    \end{tabular}%
\end{table}%

\begin{figure}
\begin{centering}
\begin{minipage}[t]{0.33\columnwidth}%
      \includegraphics[width=4.5cm]{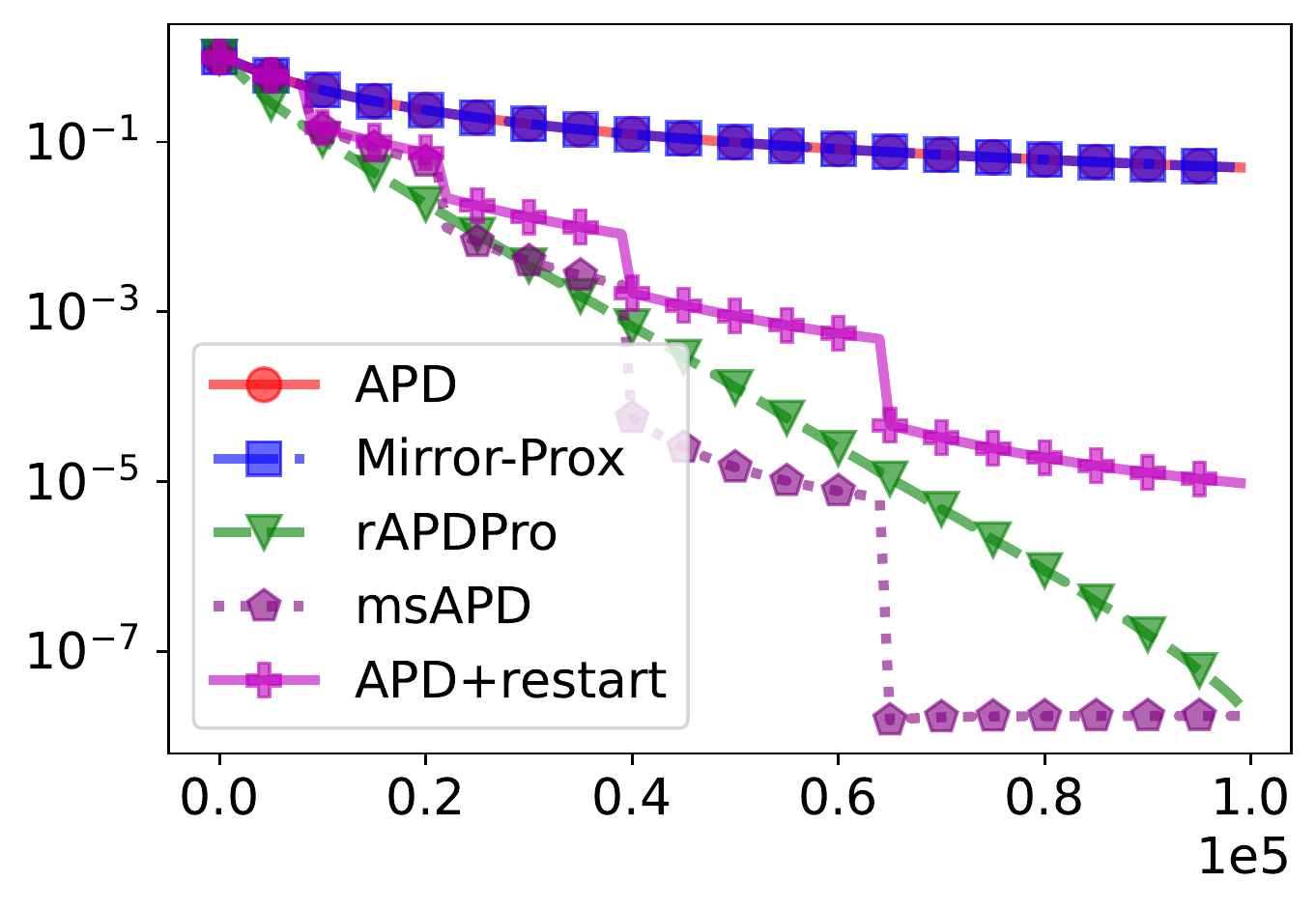}%
\end{minipage}%
\begin{minipage}[t]{0.33\columnwidth}%
\includegraphics[width=4.5cm]{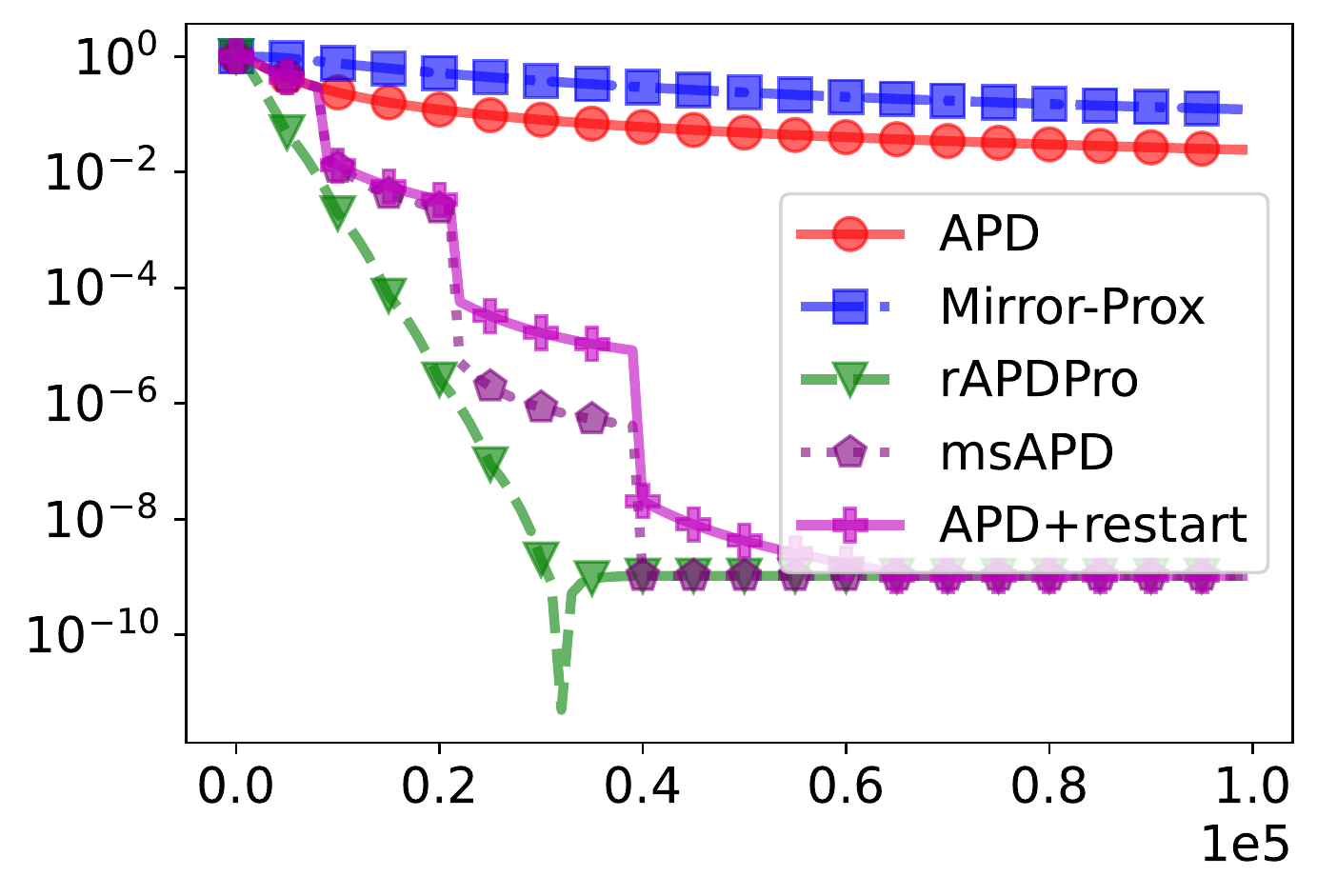}%
\end{minipage}%
\begin{minipage}[t]{0.33\columnwidth}%
\includegraphics[width=4.5cm]{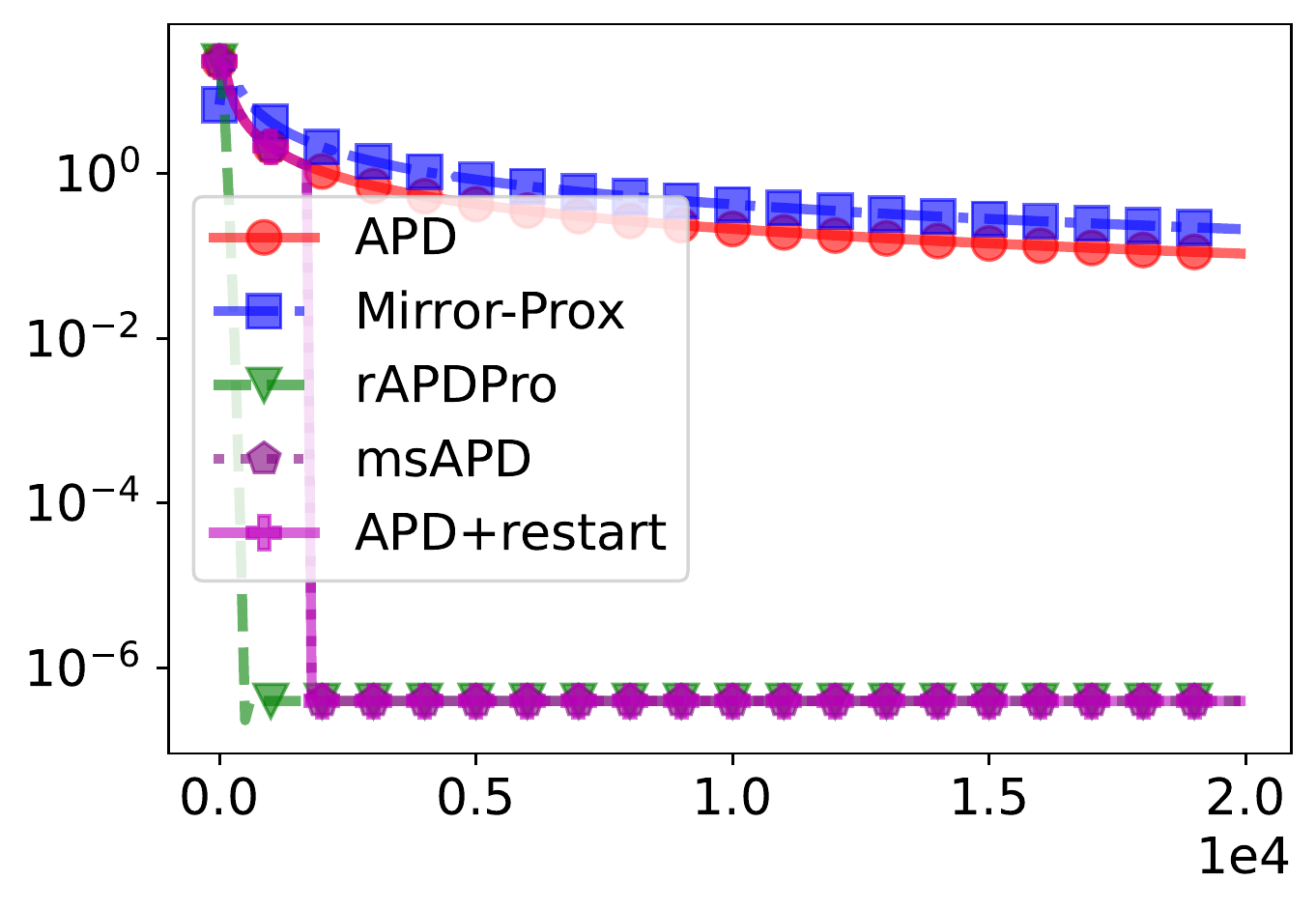}%
\end{minipage}%
\\
\begin{minipage}[t]{0.33\columnwidth}%
      \includegraphics[width=4.5cm]{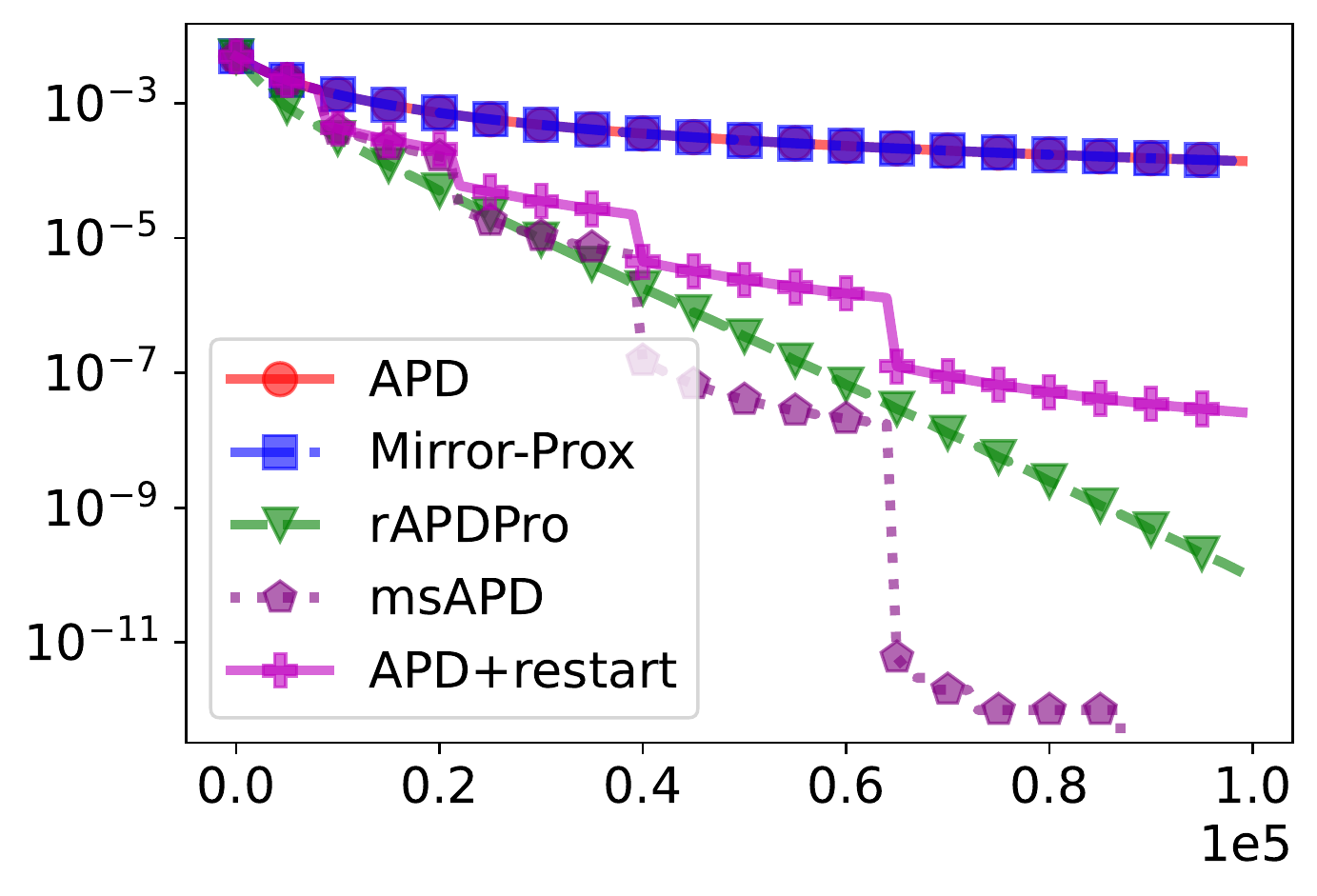}%
\end{minipage}%
\begin{minipage}[t]{0.33\columnwidth}%
\includegraphics[width=4.5cm]{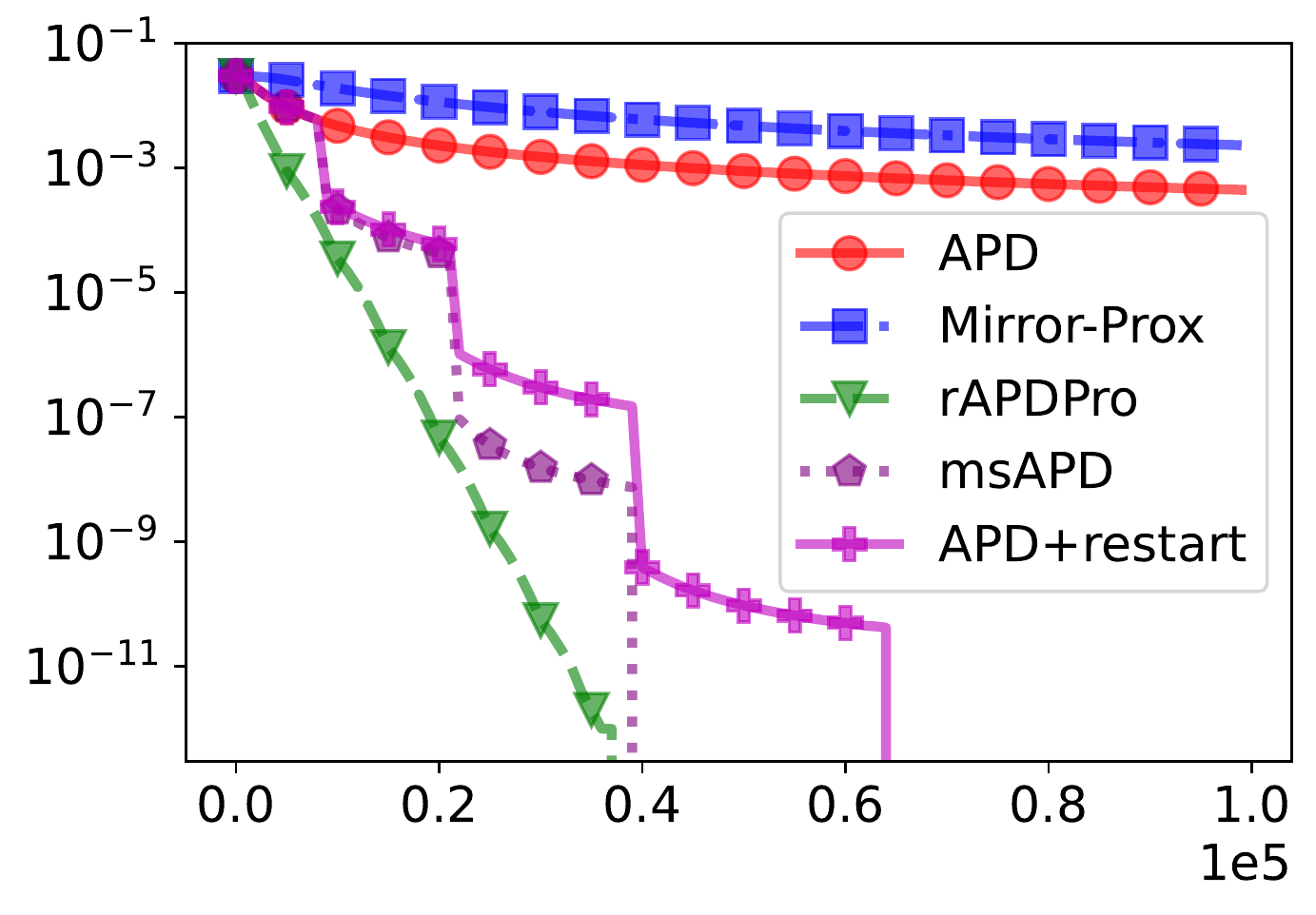}%
\end{minipage}%
\begin{minipage}[t]{0.33\columnwidth}%
\includegraphics[width=4.5cm]{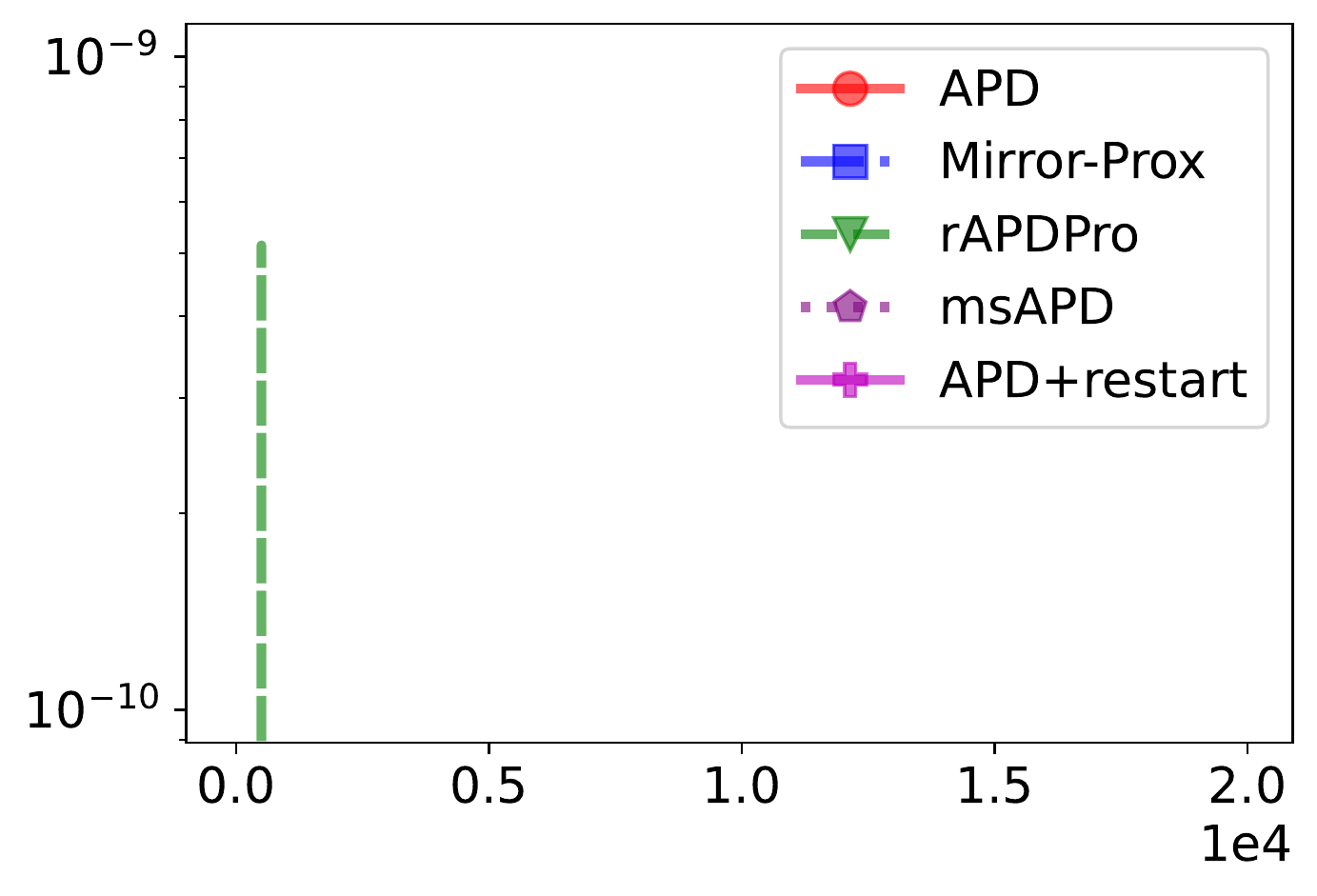}%
\end{minipage}%
\par\end{centering}
\caption{\label{fig:feasible_gap}
The first row is the results of objective convergence to optimum, where the $y$-axis reports $\log_{10}((\|D^{1/2}\xbf_{k}\|_{1}-\|D^{1/2}\xbf^{*}\|_{1})/\|D^{1/2}\xbf^{*}\|_{1})$
for {\TPAPD}, and  $\log_{10}((\|D^{1/2}\bar{\xbf}_{k}\|_{1}-\|D^{1/2}\xbf^{*}\|_{1})/\| D^{1/2}\xbf^{*}\|_{1})$
for APD, {\MSAPD} and {\mirror}. The second row is the results of feasibility violation, where $y$-axis reports the feasibility gap $\log_{10}(\max\{0, G(\protect\xbf_k)\})$
for {\TPAPD}, and $\log_{10}(\max\{0,G(\bar{\protect\xbf}_k)\})$ for
APD, APD+restart {\MSAPD} and {\mirror}.  Datasets (Left-Right order) correspond to DD68, DD242 and peking-1.}
\end{figure}

\section{\label{sec:MSAPD}A multi-stage accelerated primal-dual algorithm}
Both the previous algorithms need to solve a complicated dual problem that involves a linear cut constraint, posing a potential issue: the associated sub-problem might lack a closed-form solution.
To resolve this issue,  we present the Multi-Stage Accelerated Primal-Dual Algorithm ({\MSAPD}) in Algorithm~\ref{alg:MSAPD}, which obtains the same $\Ocal(1/\sqrt{\vep})$ complexity without introducing a new cut constraint. 
Our new method is a double-loop procedure for which an accelerated primal-dual algorithm with a pending sub-iteration number (\APDPI) is running in each stage. 
While both \APDPI{} and {\APDSCE} employ the \IMPROVE{} step to estimate the dual lower bound,   {\APDPI} only relies on the lower bound estimation to change the inner-loop iteration number adaptively, but not the stepsize selection.

\begin{algorithm}[h]
\caption{\label{alg:MSAPD}\uline{M}ulti-\uline{S}tage \uline{APD}
(\MSAPD{})}
     \begin{algorithmic}[1]
        \Require{$\bar{\mathbf{x}}^{0}\in \mathcal{X},\bar{\mathbf{y}}^{0}\in \mathcal{Y},\tilde{\sigma}, S$}
        \State{\textbf{Initialize:} $\rho_{0}^0 = 0$}
        \For{$s=0,\ldots,S$}
        \State{Compute $\tau_{0}^{s}=\brbra{L_{XY}+{L_{G}^{2}\tilde{\sigma}}\cdot 2^{\frac{s}{2}}}^{-1},\sigma_{0}^{s}=\tilde{\sigma}\cdot 2^{\frac{s}{2}}$\label{lst:line:msapd_step-size}}
        \State{$(\bar{\mathbf{x}}^{s+1},\bar{\mathbf{y}}^{s+1},\rho^{s+1}_0)\leftarrow$ \Call{APDPi}{$\tau^s_0,\sigma^s_0,\bar{\mathbf{x}}^{s},\bar{\mathbf{y}}^{s},\rho_{0}^s,s$}}
        \EndFor
    \State{\textbf{Output: $\bar{\mathbf{x}}^{S+1},\bar{\mathbf{y}}^{S+1}$}}
    \Procedure{APDPi}{$\tau_{0}^{s}, \sigma_{0}^s, \xbf_0, \ybf_{0}, \rho_{0}^s,  s$}
        \State{\textbf{Initialize:} $\left(\mathbf{x}_{-1}, \mathbf{y}_{-1}\right) \leftarrow\left(\mathbf{x}_{0}, \mathbf{y}_{0}\right),\bar{\mathbf{x}}_0=\mathbf{x}_0,k=0,N_s=\infty$, $\Delta_{XY} = \frac{1}{2\tau_0^s}D_X^2+\frac{1}{2\sigma_0^s}D_Y^2$}
        \While{$k<N_s$}
        \State{$\mathbf{z}_k \leftarrow 2G(\mathbf{x}_k)- G(\mathbf{x}_{k-1})$}
\State{$\ybf_{k+1}\leftarrow \argmin_{\ybf\in\mcal Y}\|\ybf - (\ybf_k + \sigma_k\zbf_k)\|^2 $\label{lst:line:updateyinmsapd}}
\State{$\mathbf{x}_{k+1}\leftarrow\prox_{f,\mcal{X}}(\mathbf{x}_k-\tau_{0}^{s}\nabla G(\mathbf{x}_k)\mathbf{y}_{k+1},\tau_0^s)$}
\State{
$\bar{\mathbf{x}}_{k+1}\leftarrow (k\bar{\mathbf{x}}_k + \mathbf{x}_{k+1})/(k+1),$
}
\State{$\rho_{k+1}^s\leftarrow$\Call{Improve}{$\xbf_k,  \bar{\xbf}_k, \frac{1}{2}D_{X}^2, \frac{\Delta_{XY}}{k}, \rho_{k}^s$}}
\State{Compute $N_{s}=\lceil\max\bcbra{\frac{4}{\rho_{k+1}^{s}\tau_{0}^{s}},\frac{D_{Y}^{2}}{\rho_{k+1}^{s}\sigma_{0}^{s}D_{X}^{2}}\cdot 2^{s+1}}\rceil$}
        \State{$k\leftarrow k + 1 $}
        \EndWhile
    \State{\Return{$\bar{\mathbf{x}}_{N_s},\bar{\mathbf{y}}_{N_s},\rho_{k}^s$}}
    \EndProcedure
        \end{algorithmic}
\end{algorithm}

We develop the convergence property of {\APDPI}, which paves the path to proving our main theorem.
For the convergence analysis, it suffices to verify that the initial stepsize parameter $\tau^{s}_0,\sigma^{s}_0$
satisfy assumptions in Theorem~\ref{thm:onestageconver}.
\begin{theorem}
\label{thm:onestageconver}Let $\left\{ \bar{\xbf}_{k}^s,\bar{\ybf}_{k}^s\right\}$
be the sequence generated by {\APDPI},
then  we have
\begin{equation}\label{eq:mthm}
    \mathcal{L}(\bar{\xbf}_{K}^{s},\ybf^{*})-\mathcal{L}(\xbf^{*},\bar{\ybf}_{K}^{s}) \leq\tfrac{1}{K}\Delta^{s}(\xbf^{*},\ybf^{*}),\ \ 
\tfrac{1}{2}\|\bar{\xbf}_{K}^{s}-\xbf^{*}\|^{2} \leq\tfrac{1}{(\ybf^{*})^{\top} \mubm K}\Delta^{s}(\xbf^{*},\ybf^{*}),
\end{equation}
where $\Delta^{s}(\xbf^{*},\ybf^{*})\triangleq\tfrac{1}{2\tau_{0}^{s}}\|\xbf_{0}^{s}-\xbf^{*}\|^{2}+\tfrac{1}{2\sigma_{0}^{s}}\|\ybf_{0}^{s}-\ybf^{*}\|^{2}$ and $(\xbf^*,\ybf^*)$ is a KKT point.
\end{theorem}
\proof
The stepsize $\tau_{k}^{s}=\tau_{0}^{s},\sigma_{k}^{s}=\sigma_{0}^{s}$
are unchanged at one epoch, which implies that $\rho_{k+1}= 0$, i.e., \eqref{eq:step-size-1} are satisfied. By the definition of $\tau_{0}^{s}$ and $\sigma_{0}^{s}$, we have
$
(\tau_{0}^{s})^{-1}=L_{XY}+{L_{G}^{2}\tilde{\sigma}}\sqrt{2}^{s}=L_{XY}+{L_{G}^{2}\sigma_{0}^{s}},
$
which means equality holds at the first term in~\eqref{eq:step-size-1}.

Since $g_{i}(\xbf)$ is a strongly convex function with modulus $\mu_{i}$,
then we have
\begin{equation*}
\mathcal{L}(\bar{\xbf}_{K},\ybf^{*})  \geq\mathcal{L}(\xbf^{*},\ybf^{*})+\frac{(\ybf^{*})^{\top}\boldsymbol{\mu}}{2}\|\bar{\xbf}_{K}-\xbf^{*}\|^{2},\ \ \mathcal{L}(\xbf^{*},\ybf^{*}) \geq\mathcal{L}(\xbf^{*},\bar{\ybf}_{K}).
\end{equation*}
Summing up the two inequalities above, we can get
\begin{equation}
\mathcal{L}(\bar{\xbf}_{K},\ybf^{*})-\mathcal{L}(\xbf^{*},\bar{\ybf}_{K})\geq\frac{(\ybf^{*})^{\top}\boldsymbol{\mu}}{2}\|\text{\ensuremath{\bar{\xbf}_{K}-\xbf^{*}}}\|^{2}.\label{eq:12-1-1}
\end{equation}
Combining~\eqref{eq:mthm} and~\eqref{eq:12-1-1}, we can obtain the second term of~\eqref{eq:mthm}.
\endproof

We show \MSAPD{} obtains an $\mcal O(1/\sqrt{\vep})$ convergence rate, which matches the complexity of {\APDSCE}.
\begin{theorem}
\label{thm:MSAPD}
Let $\{\bar{\xbf}_{0}^{s}\}$ be the sequence computed
by {\MSAPD}. 
Then, we have
\begin{equation}
\|\bar{\xbf}_{0}^{s}-\xbf^{*}\|^{2}\leq\Delta_{s}\equiv D_{X}^{2}\cdot 2^{-s},\quad\forall s\geq0.\label{appendix:eq:ratiodecrease-1}
\end{equation} 
For any $\vep\in(0, D_{X}^{2})$, {\MSAPD}
will find a solution $\bar{\xbf}_{0}^{s}\in\mathcal{X}$ such that
$\|\bar{\xbf}_{0}^{s}-\xbf^{*}\|^2\leq\vep$ 
in at most $\left\lceil \log_{2}{ D_{X}^{2}}/{\vep}\right\rceil $
epochs. Moreover, the overall iteration number performed by {\MSAPD}
to find such a solution is bounded by
\[
T_{\vep}=\left(\frac{8L_{XY}}{\rho_{N_{0}}^{0}}+2\right)\left\lceil\log_{2}\frac{D_{X}}{\sqrt{\vep}}+1\right\rceil+(2+\sqrt{2})\left(\tilde{\sigma}L_{G}^{2}+\frac{2D_{Y}^{2}}{\rho_{N_{0}}^{0}\tilde{\sigma}D_{X}^{2}}\right)\frac{D_{X}}{\sqrt{\vep}}.
\]
\end{theorem}
\proof
We first show that~\eqref{appendix:eq:ratiodecrease-1} holds by induction.
It is easy to verify that~\eqref{appendix:eq:ratiodecrease-1} holds for $s=0$.
Assume $\|\bar{\xbf}_{0}^{s}-\xbf^{*}\|^{2}\leq\Delta_{s}= D_{X}^{2}\cdot 2^{-s}$
holds for $s=0,\ldots,S-1$.
By Theorem~\ref{thm:onestageconver}, we have
% \begin{equation*}
% \begin{split}\norm{\bar{\xbf}_{0}^{S}-\xbf^{*}}^{2} & \le\frac{1}{(\ybf^{*})^{\top}\mubm N_{S-1}}\brbra{\frac{1}{\tau_{0}^{S-1}}\norm{\bar{\xbf}_{0}^{S-1}-\xbf^{*}}^{2}+\frac{1}{\sigma_{0}^{S-1}}\norm{\bar{\ybf}_{0}^{S-1}-\ybf^{*}}^{2}}\\
%  & \leq\frac{1}{(\ybf^{*})^{\top}\mubm N_{S-1}}\brbra{\frac{2}{\tau_{0}^{S-1}}\Delta_{S}+\frac{1}{\sigma_{0}^{S-1}}D_{Y}^{2}}.
% \end{split}
% \end{equation*}
\begin{equation*}
\norm{\bar{\xbf}_{0}^{S}-\xbf^{*}}^{2}
 \leq\frac{1}{(\ybf^{*})^{\top}\mubm N_{S-1}}\brbra{\frac{2}{\tau_{0}^{S-1}}\Delta_{S}+\frac{1}{\sigma_{0}^{S-1}}D_{Y}^{2}}.
\end{equation*}
As the algorithm sets $N_{S-1}=\big\lceil\max\bcbra{{4}/\rbra{\rho_{N_{S-1}}^{S-1}\tau_{0}^{S-1}},{2D_{Y}^{2}}/\rbra{\rho_{N_{S-1}}^{S-1}\sigma_{0}^{S-1}\Delta_{S}}}\big\rceil$,
the following inequalities hold:
\begin{equation*}
\begin{split}2\brbra{(\ybf^{*})^{\top}\mubm N_{S-1}\tau_{0}^{S-1}}^{-1} & \leq{2}\brbra{\rho_{N_{S-1}}^{S-1}N_{S-1}\tau_{0}^{S-1}}^{-1}
% \leq\frac{2}{\rho_{N_{S-1}}^{S-1}\tau_{0}^{S-1}}\cdot\frac{\rho_{N_{S-1}}^{S-1}\tau_{0}^{S-1}}{4}=
\leq \frac{1}{2},\\
{D_{Y}^{2}}\brbra{(\ybf^{*})^{\top}\mubm N_{S-1}\sigma_{0}^{S-1}}^{-1} & \leq{D_{Y}^{2}}\brbra{\rho_{N_{S-1}}^{S-1}N_{S-1}\sigma_{0}^{S-1}}^{-1}
% \leq\frac{D_{Y}^{2}}{\rho_{N_{S-1}}^{S-1}\sigma_{0}^{S-1}}\cdot\frac{\rho_{N_{S-1}}^{S-1}\sigma_{0}^{S-1}}{2D_{Y}^{2}}\Delta_{S}=
\leq \frac{1}{2}\Delta_{S}.
\end{split}
\end{equation*}
Putting these pieces together, we have $\norm{\bar{\xbf}_{S}-\xbf^{*}}^{2}\le\frac{1}{2}\Delta_{S}+\frac{1}{2}\Delta_{S}=\Delta_{S}$.
Suppose the algorithm runs for $S$ epochs to achieve the desired
accuracy $\vep$, i.e., $\|\xbf_{0}^{S}-\xbf^{*}\|^{2}\le D_{X}^{2}\cdot2^{-S}\leq\vep$.
Then the overall iteration number can be bounded by 
\begin{equation*}
\begin{split}\sum_{s=0}^{S}N_{s} & \overset{(a)}{\leq}\sum_{s=0}^{S}\Bcbra{\frac{4}{\rho_{N_{0}}^{0}\tau_{0}^{S-1}}+\frac{2D_{Y}^{2}}{\rho_{N_{0}}^{0}\sigma_{0}^{S-1}\Delta_{S}}+1}\\
 & \overset{(b)}{\leq}\sum_{s=0}^{S}\Bcbra{\Brbra{\frac{4L_{XY}}{\rho_{N_{0}}^{0}}+1}+\Brbra{\tilde{\sigma}L_{G}^{2}+\frac{2D_{Y}^{2}}{\rho_{N_{0}}^{0}\tilde{\sigma}D_{X}^{2}}}\sqrt{2}^{s}}\\
 % & \leq(\frac{4L_{XY}}{\rho_{N_{0}}^{0}}+1)(S+1)+(\frac{\tilde{\sigma}L_{G}^{2}}{1}+\frac{2D_{Y}^{2}}{\rho_{N_{0}}^{0}\tilde{\sigma}D_{X}^{2}})(2+\sqrt{2})(\sqrt{2}^{S}-\frac{\sqrt{2}}{2})\\
 & \leq\Brbra{\frac{8L_{XY}}{\rho_{N_{0}}^{0}}+2}\Big\lceil\log_{2}\frac{D_{X}}{\sqrt{\vep}}+1\Big\rceil+(2+\sqrt{2})\Brbra{\tilde{\sigma}L_{G}^{2}+\frac{2D_{Y}^{2}}{\rho_{N_{0}}^{0}\tilde{\sigma}D_{X}^{2}}}\frac{D_{X}}{\sqrt{\vep}},
\end{split}
\end{equation*}
where $(a)$ holds by $\rho_{N_{S}}^{s}\geq\rho_{N_{0}}^{0},\forall s\geq0$,
$(b)$ follows from the definition of $\tau_{0}^{s}$ and $\sigma_{0}^{s}$. 
\endproof

\begin{remark}
Theorem~\ref{thm:MSAPD} shows that \MSAPD{} obtains a worst-case complexity of  $\Ocal\brbra{\log ({D_X}/{\sqrt{\vep}})+(D_X+{D_Y^2}/{D_X})/{\sqrt{\vep}}}$, which is an upper bound of the complexity of \resapdsce{} (see Theorem~\ref{thm:rapdpro}).
The complexities of \MSAPD{} and \resapdsce{} match when $D_X=\Omega(1)D_Y$. Otherwise, \resapdsce{} appears to be much better in terms of dependence on ${D_X}/{\sqrt{\vep}}$.  On the other hand,  {\MSAPD} has a simpler subproblem, which does not involve an additional cut constraint on the dual update. 
\end{remark}

\section{ Experiment details \label{sec:appendix_exp_details}}
We examine the empirical performance for solving sparse Personalized PageRank. Let $G=(V,E)$ be a connected undirected graph with $n$ vertices. Denote the adjacency matrix of $G$ by $A$, that is, $A_{i,j}=1$ if $i\sim j$ and $0$ otherwise. Let $D=\diag(d_1,\ldots,d_n)$ be the matrix with the degrees $\{d_i\}_{i=1}^n$ in its diagonal. Then the constrained form of Personalized PageRank can be written as follows:
\begin{equation}
    \min_{\xbf\in \mbb R^n}\ \ \norm{D^{1/2}\xbf}_1\ \  \st\,\, \tfrac{1}{2}\inprod{\xbf}{Q\xbf}-\alpha \inner{\mathbf{s}}{D^{-1/2}\xbf}\leq b,
\end{equation}
where $Q=D^{-1/2}\brbra{D-\tfrac{1-\alpha}{2}(D+A)}D^{-1/2}$, $\alpha \in (0,1)$, $\mathbf{s}\in \Delta^n$ is a teleportation distribution over the nodes of the graph $G$ and $b$ is a pre-specific target level.

\paragraph{Datasets} We selected 6 small-to-median scale datasets from various domains in the Network Datasets~\citep{networkdata}. We skip large-scale networks as MOSEK struggles to achieve the optimal solution, making it unsuitable for subsequent comparison of the optimality gap. We briefly describe these datasets in Table~\ref{tab:dataset_description}. For more details, please refer to the \href{https://networkrepository.com/networks.php}{network repository}.

\paragraph{Parameter tuning} 
For all experiments, we set $r=\min_{i\in[n]}|d_i|$, $\underline{\mu} = \lambda_{\min}(Q)$ and $L_X = \lambda_{\max}(Q)$, with $\lambda_{\min}(\cdot), \lambda_{\max}(\cdot)$ denoting the smallest and largest eigenvalue, respectively. For {\MSAPD}, we have made additional parameter adjustments. Based on our observations, due to a small estimated strongly convex coefficient, {\MSAPD} could not switch to the next cycle $s$ early enough. To prevent {\MSAPD} from degrading to {\APD}, we iterate according to the predefined number of sub-iterations and manually switch to the next set of parameters. We divide $\tau$ by $\sqrt{2}$, multiply $\sigma$ by $\sqrt{2}$, and increase the number of sub-iterations in the next period by a factor of  $\sqrt{2}$. For all experiments, we tune the stepsize $\tau,\sigma,\gamma$ from $\bcbra{0.0001,0.0005,0.001,0.005,0.01}$, where $\tau, \sigma$ are the initial stepsizes of {\resapdsce}, {\MSAPD} and {\APD}, $\gamma$ is the constant stepsize of {\mirror}. All algorithms start with the primal variables initialized as zero vectors and the dual variables initialized as ones.

\paragraph{Additional experiment results}

\begin{figure}
\begin{centering}
\begin{minipage}[t]{0.33\columnwidth}%
\includegraphics[width=4.5cm]{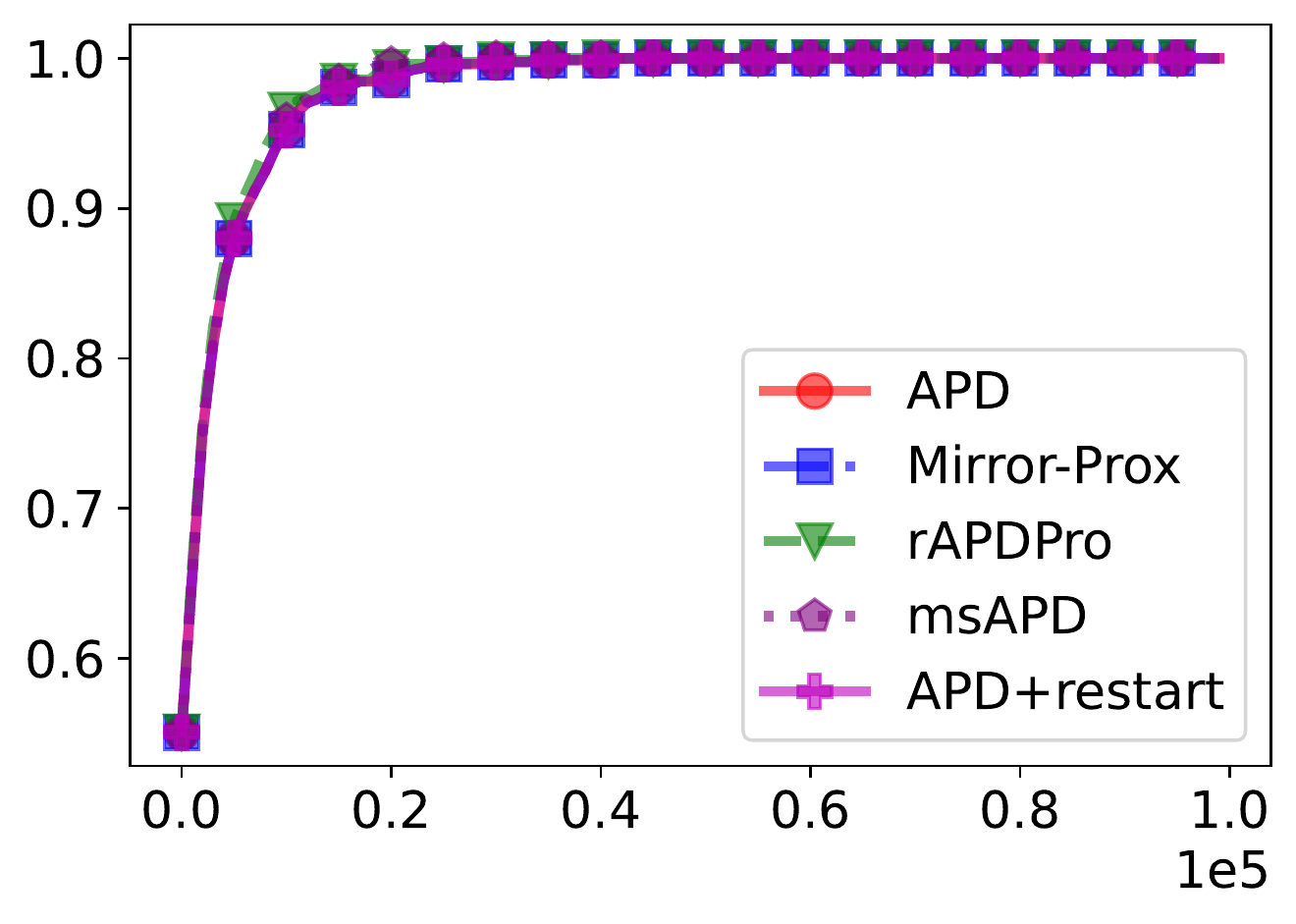}%
\end{minipage}%
\begin{minipage}[t]{0.33\columnwidth}%
\includegraphics[width=4.5cm]{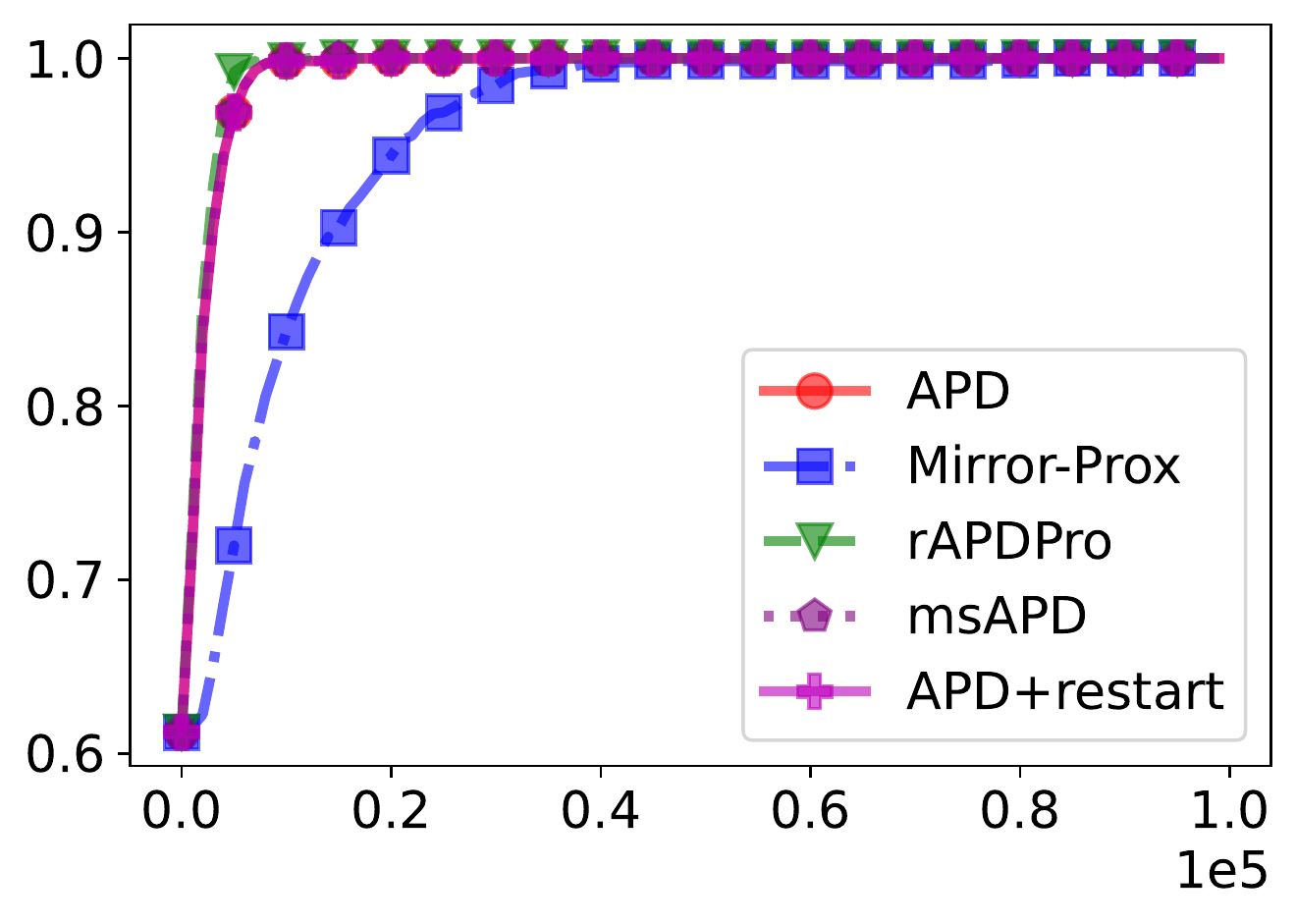}%
\end{minipage}%
\begin{minipage}[t]{0.33\columnwidth}%
\includegraphics[width=4.5cm]{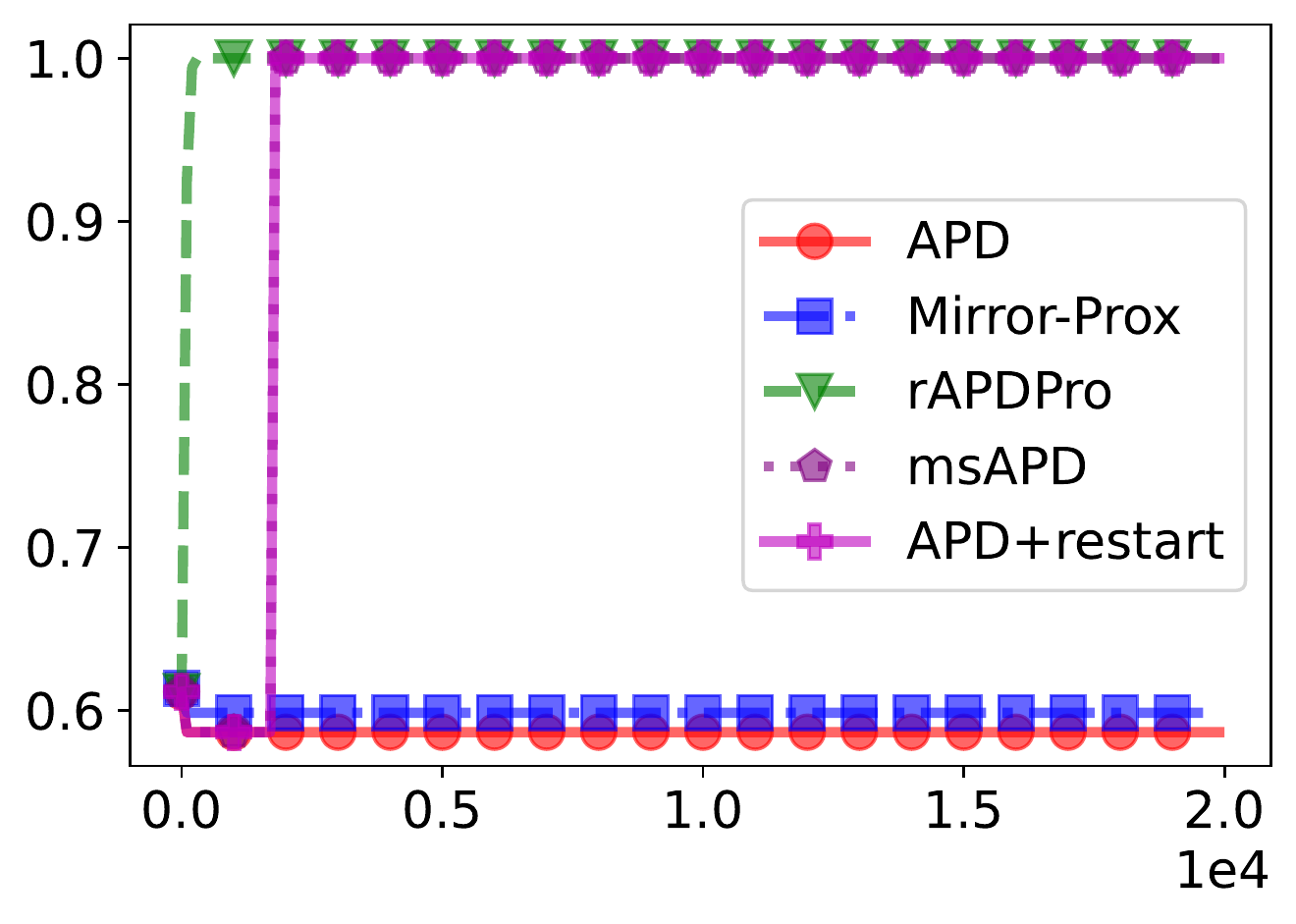}%
\end{minipage}%
\par\end{centering}
\caption{\label{fig:identification_appendix}The experimental results on active-set identification. Datasets (Left-Right order) correspond to DD68, DD242 and peking-1.
The $x$-axis reports the iteration number and  the $y$-axis reports accuracy in active-set identification.}
\end{figure}

Figure~\ref{fig:feasible_gap} and Figure~\ref{fig:identification_appendix} describe the convergence performance and active set identification results on the last three datasets:  DD68, DD242 and peking-1. Furthermore, we report the time consumption for the Personalized PageRank problem in Table~\ref{tab:time_summary}. The table indicates that, although {\resapdsce} and {\MSAPD} require moderately complex computations to determine the lower bound of the strong convexity parameter, the two methods still accelerate the algorithm's convergence and can significantly reduce the overall convergence time. 

\begin{table}[htbp]
  \centering
  \caption{Time summary when $\max\{\abs{f(\xbf)-f(\xbf^*)}/\abs{f(\xbf^*)},\max\{G(\xbf),0\}\}\leq 10^{-3}$. All experiments were conducted five times, and the results are reported as mean (standard deviation). $*$ means that upon completion of all iterations, the algorithms still fails to meet the criteria for both error measures.}
  \setlength{\tabcolsep}{1.3mm}
    \begin{tabular}{l|lllll|r}
    \toprule
    dataset & {\APD}  & APD+restart & {\resapdsce} & {\mirror} & {\MSAPD} & \multicolumn{1}{l}{mosek} \\
    \midrule
    bio-CE-HT & 187.15 (0.86)* & 115.95 (1.04) & 136.92 (0.92) & 370.50 (1.80)* & \textbf{77.21} (0.67) & 0.21 \\
    bio-CE-LC & 2.58 (0.16)* & 0.65 (0.01) & \textbf{0.44} (0.01) & 4.74 (0.33)* & 0.65 (0.03) & 0.1 \\
    econ-beaflw & 72.28 (0.59)* & 87.12 (0.43)* & \textbf{18.42} (0.44) & 116.13 (1.15)* & 66.70 (0.76) & 0.16 \\
    DD242 & 43.29 (1.20)* & 10.27 (0.39) & \textbf{6.30} (0.08) & 79.16 (0.60)* & 10.33 (0.62) & 0.16 \\
    DD68  & 36.55 (0.42)* & 19.07 (0.66) & 22.35 (0.75) & 67.73 (1.39)* & \textbf{15.69} (0.37) & 0.24 \\
    peking-1 & 122.37 (2.99)* & 11.55 (0.69) & \textbf{4.86} (0.09) & 243.45 (7.20)* & 11.24 (0.15) & 0.21 \\
    \bottomrule
    \end{tabular}%
  \label{tab:time_summary}%
\end{table}%

Nonetheless, we observe that Mosek achieves significantly faster computational efficiency for small-scale problems than our algorithm. Therefore, we test the efficiency of {\resapdsce} on some large-scale instances. For large-scale instances, we consider the following problem $\min_{\xbf \in \mbb R^n}\norm{\xbf-1}_1 \,\, \st \frac{1}{2}\xbf^\top Q_i \xbf + c_i^\top \xbf + d_i\leq 0,i=1,\ldots,m,$ where $Q_i$ are dense and positive definite matrix and generated randomly and $c_i$ are generated randomly. Furthermore, we set proper $d_i$ to make the feasible region is non-empty. When $n=5000$ and $m>10$, MOSEK crashes on our computer, which means we can not get $\xbf^*$ for calculating the optimality gap. Therefore, we report the time required for the algorithm to satisfy $\max\{\abs{f(\xbf)-f(\xbf^*)}/\abs{f(\xbf^*)},\max\{G(\xbf),0\}\}\leq 10^{-3}$ and the time taken by the algorithm to complete 10,000 iterations. On this problem, results from small datasets indicate that the performance of the 10,000-step algorithm should be sufficient to meet our specified termination criteria.
\begin{table}[htbp]
  \centering
  \caption{Comparison of computational time in seconds between {\resapdsce} and MOSEK}
    \begin{tabular}{rrr}
    \toprule
    \multicolumn{1}{l}{m} & \multicolumn{1}{l}{{\resapdsce}} & \multicolumn{1}{l}{MOSEK} \\
    \midrule
    8     & 24.612 & 50.38 \\
    10    & 53.997 & 67.99 \\
    12    & 392   & - \\
    \bottomrule
    \end{tabular}%
  \label{tab:addlabel}%
\end{table}%